\documentclass[onefignum,onetabnum]{siamonline171218}

\usepackage{lipsum}
\usepackage{amsfonts}
\usepackage{graphicx}
\usepackage{epstopdf}
\usepackage{algorithmic}

\usepackage{enumitem}
\setlist[enumerate]{leftmargin=.5in}
\setlist[itemize]{leftmargin=.5in}

\newsiamremark{remark}{Remark}
\newsiamremark{hypothesis}{Hypothesis}
\crefname{hypothesis}{Hypothesis}{Hypotheses}
\newsiamthm{claim}{Claim}

\headers{Statistical closure modeling by the ROMES method}{S. Pagani, A. Manzoni and K. Carlberg}

\title{Statistical closure modeling for reduced-order models\\
of stationary systems by the ROMES method\thanks{Submitted to the editors 9th January, 2019.
\funding{
K.\ Carlberg
was sponsored by Sandia's Advanced Simulation and Computing (ASC)
Verification and Validation (V\&V) Project/Task \#103723/05.30.02.
This paper describes objective technical results and analysis. Any subjective
views or opinions that might be expressed in the paper do not necessarily
represent the views of the U.S. Department of Energy or the United States
Government. Sandia National Laboratories is a multimission laboratory managed
and operated by National Technology \& Engineering Solutions of Sandia, LLC, a
wholly owned subsidiary of Honeywell International Inc., for the U.S.\
Department of Energy's National Nuclear Security Administration under contract
DE-NA0003525.}}}

\author{
Stefano Pagani\footnotemark[1] \and
     Andrea Manzoni\footnotemark[1] \and Kevin Carlberg\footnotemark[2]  } 

\usepackage{amsopn}

\usepackage[english]{babel}

\usepackage{hyperref}
\usepackage{graphicx}
\usepackage{subfigure}
\usepackage{framed}

\usepackage{array}
\newcolumntype{C}[1]{>{\centering\arraybackslash}p{#1}}

\usepackage{cancel}

\usepackage{algorithm,algorithmic}

\usepackage{pgfplots}
\usepackage{tikz}

\pgfplotsset{ 
  title style = {font=\small}, 
  legend style = {font=\tiny}, 
  label style = {font=\tiny}, 
  xticklabel style={
          /pgf/number format/fixed,
          /pgf/number format/precision=3
  },
  yticklabel style={
          /pgf/number format/fixed,
          /pgf/number format/precision=3
  }
}
\pgfplotsset{compat=newest}
\pgfplotsset{plot coordinates/math parser=false}
\newlength\figureheight
\newlength\figurewidth

\pgfmathdeclarefunction{lg10}{1}{ \pgfmathparse{ln(#1)/ln(10)}}

\usepackage{mathtools}

\usepackage{amsmath}
\usepackage{amsfonts}
\usepackage{amssymb}
\usepackage{dsfont}
\usepackage{stmaryrd}
\usepackage{color}
\usepackage{multirow}

\usepackage{booktabs}

\makeatletter
\def\BState{\State\hskip-\ALG@thistlm}
\makeatother

\usepackage{graphicx}
\graphicspath{{Img/}}

\DeclareMathOperator*{\argmin}{\arg\min}

\usepackage{geometry}
\usepackage{longtable}

\usepackage{listings}
\usepackage{color}

\newcommand{\parallelsum}{\mathbin{\|}}
\newcommand{\defeq}{\vcentcolon=}
\newcommand{\ndof}{{N}}
\newcommand{\nparam}{{d}}
\newcommand{\paramDomain}{\mathcal D}

\newcommand{\paramEval}{\paramDomain_\mathrm{online}}
\newcommand{\paramEvalArg}[1]{\param^\star_{#1}}
\newcommand{\paramTrain}{\paramDomain_\mathrm{POD}}
\newcommand{\paramDual}{\paramDomain_\mathrm{dual}}
\newcommand{\neval}{n_\text{online}}
\newcommand{\paramOutofplane}{\paramDomain_\textrm{out-of-plane}}
\newcommand{\paramROMES}{\mathcal D_\mathrm{ROMES}}

\newcommand{\FVUArg}[1]{\text{FVU}_{#1}}
\newcommand{\FVUQoiArgs}[2]{\text{FVU}_{\qoiFOMArg{#1}}(#2)}
\newcommand{\FVUQoiMaxArgs}[1]{\text{FVU}_{\qoiFOM,\max}(#1)}
\newcommand{\ROMEStrainingArg}[1]{\mathcal T_{#1}}
\newcommand{\ROMEStrainingArgs}[2]{\mathcal T_{#1,#2}}

\newcommand{\linesearchIt}[1]{\alpha^{(#1)}}
\newcommand{\loss}[2]{L_{#1}(#2)}
\newcommand{\lossNo}{L_{i,j}}
\newcommand{\lossArg}[3]{L_{#2,#3}(#1)}
\newcommand{\lossloglikelihoodNo}{L_{\text{likelihood},i,j}}
\newcommand{\lossloglikelihoodArg}[3]{L_{\text{likelihood},#2,#3}(#1)}
\newcommand{\lossPredictionNo}{L_{\text{C},i,j}}
\newcommand{\lossPredictionArg}[3]{L_{\text{C},#2,#3}(#1)}
\newcommand{\lossFrequencyNo}{L_{\frequency,i,j}}
\newcommand{\lossFrequencyArg}[3]{L_{\frequency,#2,#3}(#1)}

\newcommand{\lossEightyNo}{L_{0.80,i,j}}

\newcommand{\lossNinetyFiveNo}{L_{0.95,i,j}}

\newcommand{\lossKSNo}{L_{\text{KS},i,j}}
\newcommand{\lossKSArg}[3]{L_{\text{KS},#2,#3}(#1)}
\newcommand{\param}{\boldsymbol{\mu}}
\newcommand{\paramOnline}{\param^\star}
\newcommand{\validationFrequency}[2]{\omega_{#1}(#2)}
\newcommand{\validationFrequencyVal}[4]{\omega_{#1,#3}(#2;#4)}
\newcommand{\predictionInterval}[3]{C_{#1}(#2,#3)}
\newcommand{\predictionIntervalVal}[5]{C_{#1,#4}(#3,#2;#5)}
\newcommand{\Span}[1]{\mathrm{span}\{#1\}}
\newcommand{\paramDummy}{\boldsymbol{\nu}}
\newcommand{\meanPredictionInt}{\overallmeanArg{i}{\indicatorArg{i}(\param)}}
\newcommand{\meanPredictionIntVal}{\overallmeanArg{i,j}{\indicator;\hyperparams}}
\newcommand{\sigmaPredictionInt}{\overallSTDArg{i}{\indicatorArg{i}(\param)}}
\newcommand{\sigmaPredictionIntVal}{\overallSTDArg{i,j}{\indicator;\hyperparams}}
\newcommand{\erf}{\mathrm{erf}}
\newcommand{\projectorSymb}{\mathbf{P}}
\newcommand{\projectorParallel}{\projectorSymb^{\parallelsum}}
\newcommand{\projectorPerp}{\projectorSymb^\perp}

\newcommand{\projmat}{\projectorSymb}
\newcommand{\metricmat}{\mathbf{\Theta}}
\newcommand{\stateErrorRedSymb}{\delta}
\newcommand{\stateErrorRedNotation}{\boldsymbol\stateErrorRedSymb}
\newcommand{\stateErrorInplane}{\stateErrorRedNotation^{\parallelsum}}
\newcommand{\stateErrorInplaneRed}{\hat{\stateErrorRedNotation}^{\parallelsum}}
\newcommand{\stateErrorOutofplane}{\stateErrorRedNotation^{\perp}}
\newcommand{\stateErrorOutofplaneRed}{\hat{\stateErrorRedNotation}^\perp}
\newcommand{\stateErrorRed}{\hat{\stateErrorRedNotation}}
\newcommand{\stateErrorRedSingle}{\hat{\stateErrorRedSymb}}
\newcommand{\stateErrorRedArg}[1]{\hat{\stateErrorRedSymb}_{#1}}
\newcommand{\stateErrorRedAvgArg}[1]{\bar{\hat{\stateErrorRedSymb}}_{#1}}
\newcommand{\stateErrorRedVecArg}[1]{\hat{\boldsymbol\stateErrorRedSymb}_{#1}}
\newcommand{\stateErrorRedArgs}[2]{\hat{\boldsymbol\stateErrorRedSymb}_{#1,#2}}
\newcommand{\spd}[1]{\mathrm{SPD}(#1)}
\newcommand{\zerovec}{\mathbf{0}}
\newcommand{\mapping}{m}
\newcommand{\mappingArg}[1]{\mapping_{#1}}

\newcommand{\mappingApproxArg}[1]{\tilde\mapping_{#1}}
\newcommand{\mappingApproxArgsNo}[2]{\tilde\mapping_{#1,#2}}
\newcommand{\mappingApproxArgs}[4]{\tilde\mapping_{#1,#2}(#3;#4)}

\newcommand{\frequency}{\omega}
\newcommand{\indicatorFunction}[1]{\boldsymbol 1_{#1}}
\newcommand{\indicator}{\rho}
\newcommand{\indicatorArg}[1]{\rho_{#1}}
\newcommand{\indicatorVecArgs}[2]{\boldsymbol\rho_{#1,#2}}
\newcommand{\indicatorVecArg}[1]{\boldsymbol\rho_{#1}}

\newcommand{\stateErrorRedModelArg}[1]{\tilde \delta_{#1}}
\newcommand{\stateErrorRedModelArgs}[4]{\tilde \delta_{#1,#2}(#3;#4)}
\newcommand{\stateErrorInplaneRedModel}{\tilde {\boldsymbol\delta}^{\parallelsum}}
\newcommand{\stateErrorInplaneRedModelArg}[1]{\tilde \delta_{#1}^{\parallelsum}}
\newcommand{\stateErrorOutofplaneRedModel}{\tilde {\boldsymbol\delta}^\perp}
\newcommand{\stateErrorOutofplaneRedModelArg}[1]{\tilde \delta_{#1}^\perp}
\newcommand{\trialsubspace}{\mathcal V}
\newcommand{\trialsubspaceperp}{\trialsubspace^\perp}

\newcommand{\trialsubspaceperpRed}{\hat \trialsubspace^\perp}
\newcommand{\rbmat}{\mathbf{\Phi}}
\newcommand{\rbmatTot}{\bar\rbmat}
\newcommand{\rbmatTotEntry}[2]{\bar\phi_{#1#2}}
\newcommand{\testrbmat}{\mathbf{\Psi}}
\newcommand{\rbmatArg}[1]{\rbmat_{#1}}
\newcommand{\testrbmatArg}[1]{\testrbmat_{#1}}
\newcommand{\rbmatDual}{\rbmatArg{\dualSymb}}
\newcommand{\rbmatDualArg}[1]{\rbmatArg{\dualSymb,{#1}}}

\newcommand{\testrbmatDual}{\testrbmatArg{\dualSymb}}
\newcommand{\testrbmatDualArg}[1]{\testrbmatArg{\dualSymb,{#1}}}
\newcommand{\nrb}{n}
\newcommand{\nrbperp}{\nrb^\perp}
\newcommand{\nrbtot}{{\bar\nrb}}
\newcommand{\nrbArg}[1]{\nrb_{#1}}
\newcommand{\nrbDual}{\nrbArg{\dualSymb}}
\newcommand{\nrbDualArg}[1]{\nrbArg{\dualSymb,{#1}}}

\newcommand{\rbvec}[1]{\boldsymbol{\phi}_{#1}}

\newcommand{\range}[1]{\text{Ran}(#1)}

\newcommand{\res}{\mathbf{r}}
\newcommand{\jacobian}{\frac{\partial\res}{\partial\stateDummy}}

\newcommand{\stateSymb}{x}
\newcommand{\state}{\mathbf{\stateSymb}}
\newcommand{\stateModel}{\tilde\state}
\newcommand{\stateModelEntry}[1]{\tilde\stateSymb_{#1}}
\newcommand{\stateROM}{\state_\text{ROM}}
\newcommand{\stateROMEntry}[1]{\stateSymb_{\text{ROM},#1}}
\newcommand{\stateRef}{ \state_\text{ref}}

\newcommand{\stateDummy}{\mathbf{w}}
\newcommand{\stateDummyOpt}{\mathbf{y}}
\newcommand{\stateIt}[1]{\state^{(#1)}}
\newcommand{\dualSymb}{p}
\newcommand{\dual}{\mathbf{\dualSymb}}
\newcommand{\dualArg}[1]{\dual_{#1}}
\newcommand{\dualApprox}{\tilde \dual}
\newcommand{\dualApproxArg}[1]{\dualApprox_{#1}}
\newcommand{\dualRed}{\hat \dual}
\newcommand{\dualRedArg}[1]{\dualRed_{#1}}

\newcommand{\unitvec}{\mathbf{e}}
\newcommand{\unitvecArg}[1]{\unitvec_{#1}}
\newcommand{\stateRed}{\hat\state}
\newcommand{\Amat}{\mathbf{A}}

\newcommand{\rbmatperp}{\rbmat^\perp}
\newcommand{\rbmatperpT}{[\rbmat^\perp]^T}
\newcommand{\RR}[1]{\mathbb{R}^{#1}}
\newcommand{\RRstar}[1]{\mathbb{R}_\star^{#1}}

\newcommand{\xx}{\vec{x}}
\newcommand{\spaceDof}[1]{\xx_{#1}}

\newcommand{\rhs}{\mathbf{b}}
\newcommand{\mat}{\mathbf{A}}
\newcommand{\expectation}[1]{\mathbb{E}[#1]}
\newcommand{\expectationParam}[1]{\frac{1}{\card{\paramEval}}\sum_{\param\in\paramEval}\left(#1\right)}
\newcommand{\errorROM}{e_\state}
\newcommand{\errorROMQoiArg}[1]{e_{\qoiArg{#1}}}

\newcommand{\qoiMatrixExperiment}{\boldsymbol \gamma}
\newcommand{\errorROMESParallel}{\tilde e_\state^{\parallelsum}}
\newcommand{\errorROMESParallelPerp}{\tilde e_\state^{\parallelsum+\perp}}

\newcommand{\errorROMESParallelPerpQoiArg}[1]{\tilde e_{\qoiFOMArg{#1}}^{\parallelsum+\perp}}

\newcommand{\errorROMParallel}{e_\state^{\parallelsum}}

\newcommand{\errorROMParallelPerp}{e_\state^{\parallelsum+\perp}}

\newcommand{\vondermonde}{\mathbf{H}}
\newcommand{\nbasis}{n_h}
\newcommand{\basisArg}[1]{h_{#1}}
\newcommand{\vondermondeROMES}{[\mathbf{1}\ \indicatorVecArgs{i}{j}]}
\newcommand{\vondermondeROMESSingle}{[\mathbf{1}\ \indicatorVecArg{i}]}
\newcommand{\card}[1]{|#1|}
\newcommand{\qoiFunc}{\mathbf{s}}
\newcommand{\qoiFuncArg}[1]{{s}_{#1}}
\newcommand{\qoiFOM}{\mathbf{q}}
\newcommand{\qoiArg}[1]{q_{#1}}
\newcommand{\qoiFOMArg}[1]{{q}_{#1}}
\newcommand{\qoiFOMAvgArg}[1]{{\bar q}_{#1}}
\newcommand{\qoiFOMAvgMax}{{\bar q_{\max}}}
\newcommand{\qoiModel}{\mathbf{\tilde q}}
\newcommand{\qoiModelArg}[1]{{\tilde q}_{#1}}
\newcommand{\qoiROM}{\mathbf{q}_\mathrm{ROM}}
\newcommand{\qoiROMArg}[1]{q_{\mathrm{ROM},#1}}
\newcommand{\nqoi}{s}
\newcommand{\stateError}{{\boldsymbol \delta}_\state}
\newcommand{\stateErrorModel}{\tilde{\boldsymbol \delta}_\state}
\newcommand{\stateErrorModelEntry}[1]{\tilde{ \delta}_{\state,#1}}
\newcommand{\qoiError}{\boldsymbol \delta_\qoiFOM}
\newcommand{\qoiErrorArg}[1]{\delta_{\qoiArg{#1}}}
\newcommand{\qoiErrorModel}{\tilde{\boldsymbol \delta}_\qoiFOM}

\newcommand{\trainingSet}{\mathcal T}
\newcommand{\ntrain}{n_\text{train}}

\newcommand{\feature}{\boldsymbol{x}}
\newcommand{\featureTrainSet}{\mathbf{x}}

\newcommand{\nfeature}{n_{\feature}}
\newcommand{\featureTrainArg}[1]{\feature_{#1}}
\newcommand{\featurePredict}{\boldsymbol{x}^\star}
\newcommand{\featurePredictSet}{\featurePredict}
\newcommand{\featureAllSet}{\underline{\mathbf{x}}}
\newcommand{\featureAllArg}[1]{\underline{\boldsymbol{x}}_{#1}}
\newcommand{\response}{y}
\newcommand{\responseTrainArg}[1]{\response_{#1}}
\newcommand{\responseTrainVec}{\mathbf{y}}
\newcommand{\responseTrainVecROMES}{\stateErrorRedArgs{i}{j}}
\newcommand{\responseTrainVecROMESSingle}{\stateErrorRedVecArg{i}}

\newcommand{\variance}{\sigma^2}
\newcommand{\varianceArg}[1]{\sigma_{#1}^2}
\newcommand{\kernelMat}{\mathbf{K}}
\newcommand{\kernelMatFunc}[2]{\mathbf{K}(#1,#2)}
\newcommand{\kernelFuncNo}{\kappa}
\newcommand{\kernelFunc}[2]{\kernelFuncNo(#1,#2)}
\newcommand{\identity}{\mathbf{I}}
\newcommand{\width}{\ell}
\newcommand{\widthArg}[1]{\width_{#1}}
\newcommand{\signalsd}{\gamma}
\newcommand{\signalsdArg}[1]{\signalsd_{#1}}
\newcommand{\hyperparams}{\boldsymbol{\theta}}
\newcommand{\hyperparamsSet}{\boldsymbol{\Theta}}
\newcommand{\hyperparamsArg}[1]{\hyperparams_{#1}}
\newcommand{\hyperparamsKernel}{\hyperparams_{\kernelFuncNo}}
\newcommand{\hyperparamsKernelArg}[1]{\hyperparams_{\kernelFuncNo,#1}}

\newcommand{\overallmean}[1]{\nu(#1)}
\newcommand{\overallmeanArg}[2]{\nu_{#1}(#2)}
\newcommand{\overallmeanArgNo}[1]{\nu_{#1}}
\newcommand{\overallmeanDef}[1]{\kernelMatFunc{#1}{\featureTrainSet}(\kernelMatFunc{\featureTrainSet}{\featureTrainSet}+\variance\identity)^{-1}\responseTrainVec}
\newcommand{\overallVar}[1]{\bar\sigma^2(#1)}
\newcommand{\overallVarArg}[2]{\bar\sigma_{#1}^2(#2)}
\newcommand{\overallSTDArg}[2]{\bar\sigma_{#1}(#2)}
\newcommand{\overallVarArgNo}[1]{\bar\sigma_{#1}^2}
\newcommand{\overallVarDef}[1]{\kernelMatFunc{#1}{#1} -
	\kernelMatFunc{#1}{\featureTrainSet}(\kernelMatFunc{\featureTrainSet}{\featureTrainSet}
+ \variance\identity)^{-1}\kernelMatFunc{\featureTrainSet}{#1} + \variance}

\newcommand{\generalSetOneSymb}{w}
\newcommand{\generalSetOne}{\mathbf \generalSetOneSymb}
\newcommand{\generalSetOneVec}{\boldsymbol{\generalSetOneSymb}}
\newcommand{\generalSetOneVecArg}[1]{\boldsymbol{\generalSetOneSymb}_{#1}}
\newcommand{\ngeneralSetOne}{n_{\generalSetOneSymb}}
\newcommand{\generalSetTwoSymb}{z}
\newcommand{\generalSetTwo}{\mathbf \generalSetTwoSymb}
\newcommand{\generalSetTwoVec}{\boldsymbol{\generalSetTwoSymb}}
\newcommand{\generalSetTwoVecArg}[1]{\boldsymbol{\generalSetTwoSymb}_{#1}}
\newcommand{\ngeneralSetTwo}{n_{\generalSetTwoSymb}}

\newcommand{\normal}[2]{\mathcal N \left(#1,#2 \right)}
\newcommand{\samplespace}{\Omega}

\newcommand{\spatialDomain}{\Omega_{\xx}}
\newcommand{\spatialDomainArg}[1]{\Omega_{\xx,#1}}
\newcommand{\diffusionConstant}{\kappa_\text{d}}
\renewcommand{\div}{\text{div}}

\newcommand{\regweights}{\boldsymbol{\beta}}
\newcommand{\regweightsArg}[1]{\regweights_{#1}}
\newcommand{\order}{O}
\usepackage{needspace}

\begin{document}

\maketitle

\renewcommand{\thefootnote}{\fnsymbol{footnote}}
\footnotetext[1]{MOX, Dipartimento di Matematica, Politecnico di Milano, P.za Leonardo da Vinci 32, I-20133 Milano, Italy,
            \email{stefano.pagani@polimi.it}, \email{andrea.manzoni1@polimi.it}, \url{https://stefanopagani.github.io/}. }%
\footnotetext[2]{
Extreme-scale Data Science and Analytics Department, Sandia National
Laboratories, Livermore, CA 94550.X
\email{ktcarlb@sandia.gov}.}

\begin{abstract}
This work proposes a technique for constructing a statistical closure model
	for reduced-order models (ROMs) applied to stationary systems
	modeled as parameterized systems
of algebraic equations. The proposed technique extends the reduced-order-model
error surrogates
	(ROMES) method \cite{kevin:romes} to closure modeling. The original ROMES
	method applied Gaussian-process regression to construct a
statistical model that maps cheaply computable error indicators
(e.g., residual norm, dual-weighted residuals) to a random variable
for either (1) the norm of the state error or (2) the error in a scalar-valued
quantity of interest. Rather than target these two types of errors, this work
	proposes to construct a statistical model for the \textit{state error}
	itself; it achieves this by constructing statistical models for the
	generalized coordinates characterizing both the in-plane error (i.e., the error
	in the trial subspace) and a low-dimensional approximation of the
	out-of-plane error. The former can be considered a statistical closure
	model, as it quantifies the error in the ROM generalized coordinates.
Because any quantity of interest can be computed as a functional of the state,
the proposed approach enables any quantity-of-interest error to be
	statistically quantified \textit{a posteriori}, as the state-error model can
	be propagated through the associated quantity-of-interest functional.
	Numerical experiments performed on both linear and nonlinear stationary
	systems illustrate the ability of the technique (1) to improve (expected) ROM
	prediction accuracy by an order of magnitude, (2) to statistically quantify
	the error in arbitrary quantities of interest, and (3)
	to realize a more cost-effective methodology for reducing the error than a
	ROM-only approach in the case of nonlinear systems.

	\end{abstract}

 \begin{keywords}
 model reduction; error surrogate; error modeling; closure modeling;
	 Gaussian-process regression;
 	uncertainty propagation; supervised machine learning
 \end{keywords}

\section{Introduction}

Computational models of stationary systems modeled as parameterized systems
of algebraic equations (e.g., those arising from the spatial discretization of a
partial-differential-equations problem) are being increasingly used in complex
decision-making scenarios. However, such scenarios are often \textit{many
query} or \textit{real time} in nature.  For example, uncertainty
propagation often requires hundreds or thousands of solutions to the
parameterized system in order to adequately characterize uncertainties;
\textit{in-situ} structural health monitoring requires such solutions to be
computed in near real time.  As a result, employing truly high-fidelity models
characterized by large-scale systems of algebraic equations (e.g., arising from a
fine spatial discretization) is often computationally intractable.
To mitigate this computational burden, analysts often replace such
computationally expensive high-fidelity `full-order models' (FOMs) with
computationally inexpensive surrogate models, which can be categorized as {\em
(i)} \textit{data fits}, which construct a regression model (e.g., via
polynomial interpolation) that directly approximates the mapping from
(parameter) inputs to (quantity-of-interest) outputs; {\em (ii)}
\textit{lower-fidelity models}, which introduce modeling simplifications
(e.g., coarsened mesh, neglected physics); and  {\em (iii)}
\textit{reduced-order models (ROMs)} constructed by performing a projection
process on the equations governing the high-fidelity model to reduce the state-space dimensionality.
Although typically more intrusive to implement, ROMs often
yield more accurate approximations than data fits, and usually
generate more significant computational gains than lower-fidelity models. For
this reason, this work considers ROMs as the surrogate of interest.  See,
e.g., Refs.~\cite{benner2015survey,quarteroni2016RBook,hesthaven2016certified}
for reviews on reduced-order-modeling techniques.

 %
To rigorously apply ROMs within a decision-making scenario, their
error with respect to the FOM must be quantified and properly accounted for in
the ultimate prediction or assessment. In uncertainty quantification (UQ)
applications, for example, the epistemic uncertainty\footnote{The ROM error
can be considered a source of epistemic uncertainty, as it can be reduced
by employing either the original high-fidelity model or a
higher-fidelity surrogate model.} introduced by the surrogate model should be
statistically quantified
\cite{kevin:romes,manzoni2016accurate,pagani2017}; on the other hand,
risk-averse scenarios may demand a deterministic bound on a quantity of
interest to ensure it does not exceed a specified threshold.

To this end, a variety of approaches have been proposed to quantify the error
introduced by reduced-order models. 
\begin{enumerate}
 \item \textit{Error indicators}. Error indicators are quantities that
	 are informative of the error, yet are relatively inexpensive to compute.
		One example is the residual norm, i.e., the norm of the FOM residual
		evaluated at the ROM solution. This quantity can be used as an error
		indicator to guide greedy methods for parameter-space sampling
\cite{bui2008parametric,bui2008model,hinze2012residual,amsallem2015design,
wu2015adaptive,yano2018lp} and when employing ROMs within a trust-region
setting \cite{zahr2015progressive,zahrCarlbergKouri}.
	 Alternatively, dual-weighted residuals employed in adjoint
error estimation provide a first-order approximation of the error in a
scalar-valued quantity of interest; they are often used for error estimation
and adaptive mesh refinement,
	as well
	as in nonlinear model reduction
\cite{meyer2003efficient,carlberg2014adaptive}.  Unfortunately, error
indicators alone are not easily amenable to UQ, as they
generate a prediction of the error that is both deterministic and is often
significantly biased.
	\item \textit{\textit{A posteriori} error bounds}. These approaches
		derive deterministic bounds for the norm of either the state error or
(quantity-of-interest) output error; they typically require evaluating the FOM
residual at the ROM solution, as well as stability/continuity constant
bounds, and dual quantities related to the quantity of interest. Such
approaches aim to derive bounds that are rigorous,
sharp, and inexpensive to compute \cite{veroy2005certified}.
However, these objectives are
often competing, as improving bound sharpness can significantly increase the
computational cost \cite{huynh2007successive,huynh2010natural}. Heuristic
		strategies to speed up error bound evaluation have been proposed
		\cite{LMR_M2AN_11,manzoni2015heuristic,wirtzDeim}, which yield error
		estimates rather than strict bounds.
Alternatively, Ref.~\cite{hain2018hierarchical} proposes a hierarchical error estimator
that can compute sharper estimates without requiring the computation of these constants, at the
cost of solving a higher-dimensional ROM.
Similarly to error indicators, deterministic error bounds are not directly
useful for UQ applications, where a probability
distribution for the ROM error is more amenable to quantifying the ROM-induced epistemic
uncertainty.
\item \textit{Error models}. These methods directly construct a regression model for the
ROM error, i.e., they construct an approximation of the mapping from
chosen regression-model inputs (or features) to a prediction of the ROM error.
Nearly all approaches in this category employ the parameter inputs as the
regression-model inputs
\cite{gano2005hybrid,NE12,march2012provably,eldred:soc,pagani2017reduced}.
These approaches are effective when the ROM error exhibits a low variability in
the parameter space and the parameter space is low dimensional. However, the
ROM error is often a highly oscillatory function of the inputs and the
parameter space dimension is often high-dimensional in many practical
settings, which can cause the approach to fail \cite{NE12,kevin:romes}.  The
reduced-order-model error surrogates (ROMES) method
\cite{kevin:romes} addresses this problem in the case of
stationary systems. Rather than employing parameter inputs as
regression-model inputs, the ROMES method instead employs the aforementioned error
indicators and rigorous error bounds for this purpose.  Because these
quantities are cheaply computable, low-dimensional, and are often highly
informative of the ROM error, the resulting error model is typically
computationally inexpensive to evaluate, exhibits low variance, and can be
sufficiently trained and validated using a relatively small amount of training data. Further,
because the approach employs Gaussian-process regression, its prediction
corresponds to a Gaussian random variable
for the ROM error that can be readily integrated into
UQ analyses; the variance of this random variable can be
interpreted as the ROM-induced epistemic uncertainty. Ref.~\cite{TrehanCD17}
extended this work to dynamical systems; rather than requiring the user to hand
select a small number of error indicators, this approach employs
high-dimensional regression models from machine learning (e.g., LASSO, random
forests) to enable a large number of candidate error indicators to be used as
inputs of the error model.  Ref.~\cite{frenoCarlberg} also extended the method
in several ways: (1) it enabled
the quantity-of-interest errors incurred by \textit{any
approximate solution} to be
quantified, (2) it proposed a much wider range of inexpensive-to-compute
residual-based features (e.g., gappy POD approximation of the residual), and
(3) it applied a wide range of regression methods of varying capacity (e.g.,
support vector regression, artificial neural networks) within a
model-selection framework.
\end{enumerate}

Due to its ability to generate inexpensive-to-evaluate, low-variance,
statistical error models that can be trained with relatively small amounts of
training data, this work
considers the ROMES method for constructing error models. We focus in
particular on addressing one major shortcoming of the approach: it requires
constructing of a \textit{separate error model} for each quantity of interest.
In many engineering applications, the analyst is often interested in field
quantities (e.g., the solution field itself, the pressure field); constructing
an error model for each element of the associated discrete error
vector---whose dimension is the same as that of the full-order model---is
computationally intractable.  Alternatively, in exploratory contexts, the
analyst may not have \textit{a priori} knowledge of which quantities will be
of interest;
while ROMES models for each quantity of interest could in principle be
constructed \textit{a posteriori} in this case, this violates the natural
offline--online decomposition leveraged by model reduction.

To address these shortcomings of the ROMES method, this work proposes
to construct
a statistical model for the \textit{state error} itself. To
avoid the need to construct a model for each element of the high-dimensional
state error
vector, the approach decomposes the state error into the \textit{in-plane
error} (i.e., the
component belonging to the low-dimensional trial subspace) and the
\textit{out-of-plane error} (i.e., the component orthogonal to
the trial subspace).
Because the orthogonal complement of the trial subspace is high dimensional,
the method employs a low-dimensional subspace---which is orthogonal to the
trial subspace---to represent the out-of-plane error.
Then, the method constructs a statistical error model for each
generalized coordinate characterizing the low-dimensional representations of
the in-plane and out-of-plane errors.
The error model for the in-plane error can be considered a statistical closure
model, as it aims to model the error in the preserved state variables (i.e., the
generalized coordinates of the ROM solution) due to omitting the remaining
variables from the formulation.
The resulting error model can then be employed to statistically quantify
the error in the entire state. Further, it can be used to generate an error
model for any quantity of interest (including field quantities) \textit{a posteriori} by propagating the
state error through the associated quantity-of-interest functional.
Numerical experiments demonstrate the ability of the method
(1) to improve (expected) ROM prediction accuracy by
an order of magnitude, (2) to statistically quantify the error in arbitrary
quantities of interest \textit{a posteriori}, and (3) to realize a more
cost-effective methodology for reducing the error than a ROM-only approach in
the case of nonlinear stationary systems.

We note that many existing works have proposed closure models for
reduced-order models; see, e.g.,
Refs~\cite{wang2012proper,iliescu2013variational,rebollo,pan2018data,xie2018data}.
However, these methods are all applied to dynamical systems (typically in the
context of fluid dynamics), and none of these techniques constructs a
statistical model, which is essential for uncertainty
quantification.



The paper is structured as follows.
Section~\ref{sec:problem_formulation} formulates the problem by presenting the full-order model,
the reduced-order model, and the state and quantity-of-interest errors associated with
parameterized systems of algebraic equations.
Section \ref{sec:ROMESstate} describes the decomposition of the state error
into in-plane and out-of-plane components, as well as computable
dual-weighted-residuals that approximate the error in the associated generalized
coordinates.
Section \ref{sec:GPstate} describes the proposed approach, i.e., the
proposed statistical model (Section \ref{sec:statisticalmodel}), the proposed error indicators
(Section \ref{sec:errorindicator}), a summary of Gaussian-process
regression (Section \ref{sec:GP}), the application of the method to construct statistical error
models for the state and
quantities of interest (Section \ref{seq:stateqoimodels}), and the
offline--online decomposition of the approach (Section \ref{sec:offlineonline}).
Finally, Section \ref{sec:NumResults} presents numerical experiments that
assess the proposed method on both linear and nonlinear stationary
systems, focusing particularly on model validation, the expected accuracy of
the error models, and the computational efficiency of the proposed technique.
\section{Problem formulation} \label{sec:problem_formulation}

This section presents the formulations of the FOM and ROM (with attendant state
and quantity-of-interest errors) in the context of stationary
systems.

\subsection{Full-order model}

In this work, the FOM corresponds to a stationary system modeled as a
parameterized system of algebraic equations
\begin{equation} \label{eq:discretizedPB}
\res( \state ; \param ) = \mathbf{0},
\end{equation}
where
$\res:(\stateDummy;\paramDummy)\mapsto \res(\stateDummy;\paramDummy)$ with
$\res : \RR{\ndof} \times \paramDomain \rightarrow \RR{\ndof}$
denotes the residual operator,
$\param \in \paramDomain \subset \RR{\nparam}$ denotes the system parameters,
and
$\state \equiv \state(\param)\in\RR{\ndof}$ denotes the state implicitly defined as the solution to
\eqref{eq:discretizedPB} given parameters $\param$.
%
If the residual operator $\res$ is linear in its first argument, then it takes
the form
\begin{equation}\label{eq:reslinear}
\res:( \stateDummy ; \paramDummy ) \mapsto \rhs (\paramDummy)  -
\mat(\paramDummy) \stateDummy
\end{equation}
where $ \rhs: \paramDomain  \rightarrow \RR{\ndof}$ and
$\mat: \paramDomain  \rightarrow \RRstar{\ndof \times \ndof}$ denote the
(parameterized) right-hand-side vector and the system matrix, respectively, and
 $\RRstar{m\times n}$ denotes the set of full-column rank $m
\times n$ matrices (the non-compact Stiefel manifold).
If instead the residual is nonlinear in its first argument, then
Eq.~\eqref{eq:discretizedPB} can be solved iteratively, e.g.,  via globalized Newton's method by executing the following
iterations:
given initial guess $\stateIt{0}$,
solve the linear system
$$
\jacobian( \stateIt{k-1} ; \param ) \delta {\state}^{(k)} =  -\res( \stateIt{k-1} ; \param ),
$$
and set
$$
\state^{(k)} = \state^{(k-1)} + \linesearchIt{k}\delta {\state}^{(k)},
$$
for $k=1,\ldots,K$,
where $\linesearchIt{k}\in\RR{}$ denotes a step length that can be computed
to ensure global convergence (e.g., by satisfying the strong Wolfe conditions), and
$K$ is determined by the satisfaction of a convergence criterion.
%
Many practical scenarios in science and engineering are characterized by the
following attributes:
 \begin{enumerate}
		  \item
				The FOM is \textit{high-dimensional},
i.e., $\ndof$ is large. This arises when the FOM corresponds to the fine
		 spatial discretization of a stationary partial-differential-equations
		 problem, for example.
\item The primary goal of the analysis is compute
	\textit{quantities of interest} that are functionals of the state, i.e.,
	for a given parameter instance $\param$, the goal is to compute $\qoiFOM(\param)$ with
	$\qoiFOM:\paramDummy\mapsto\qoiFunc(\state(\paramDummy);\paramDummy)$ and
	$ \qoiFunc:\RR{\ndof}\times\paramDomain\rightarrow\RR{\nqoi}$ denoting the
	quantity-of-interest functional.

\item The scenario is \textit{many-query} in nature, i.e., it requires the computation
	of $\qoiFOM(\param)$ for
		 $\param\in
	\paramEval\equiv\{\paramEvalArg{i}\}_{i=1}^{\neval}\subseteq\paramDomain$ with
	$\neval$ large. This arises in parameter studies,
	UQ applications, and design-optimization settings,
	for example.
				 \end{enumerate}
In such cases, simply solving Eq.~\eqref{eq:discretizedPB} for $\param\in
\paramEval$ and subsequently computing the quantities of interest is usually
computationally intractable, and a surrogate model
is required to reduce the computational cost. As discussed in the introduction,
this work focuses on applying reduced-order models for this purpose.



\subsection{Reduced-order model}\label{sec:ROM}
Reduced-order models reduce the dimensionality of the FOM governing equations
\eqref{eq:discretizedPB} via projection.
In particular, they seek approximate solutions
$\stateROM\approx\state$
in an $\nrb$-dimensional  affine trial
subspace (with $\nrb\ll\ndof$), i.e.,
\begin{equation}\label{eq:romSol}
 \stateROM(\param)= \stateRef(\param) + \rbmat  \stateRed(\param)\in\stateRef(\param) +
 \trialsubspace,
\end{equation}
with $\stateROM:\paramDomain\rightarrow\RR{\ndof}$. Here,
$\stateRef:\paramDomain\rightarrow\RR{\ndof}$ denotes a reference state (e.g.,
the mean of the snapshots in the case of proper orthogonal decomposition); the
trial-basis matrix
$\rbmat\equiv[\rbvec{1}\ \cdots\ \rbvec{\nrb}] \in \RRstar{\ndof \times
\nrb}$ may
be
constructed by a variety of means (e.g., the reduced-basis method
\cite{rozza2007reduced,hesthaven2016certified,quarteroni2016RBook}, proper orthogonal
decomposition \cite{POD});
$\trialsubspace\defeq\range{\rbmat}$ denotes the linear part of the affine trial subspace, where $\range{\Amat}$ denotes the range
of matrix $\Amat$; and
$\stateRed:\paramDomain\rightarrow\RR{\nrb}$ denotes the generalized
coordinates of the ROM solution.

Model-reduction approaches compute the approximate solution by substituting $\state\leftarrow \stateROM$ in
Eq.~\eqref{eq:discretizedPB} and enforcing orthogonality of the residual to
an $\nrb$-dimensional linear test subspace, which yields the ROM governing equations
\begin{equation}\label{eq:romEq}
\testrbmat(\stateRed;\param)^T\res(\stateRef + \rbmat  \stateRed;\param) = \mathbf{0},
\end{equation}
where $\testrbmat:\RR{\nrb}\times\paramDomain\rightarrow \RRstar{\ndof \times
\nrb}$ denotes the test-basis matrix that generally may depend on the
generalized coordinates and parameters. Common choices for the test basis
include Galerkin projection, which employs $\testrbmat = \rbmat$, and
least-squares Petrov--Galerkin (LSPG) projection
\cite{bui2008model,LeGresleyThesis,CarlbergGappy}, which employs
$\testrbmat(\stateRed;\param) = \jacobian(\stateRef +
\stateRed;\param)\rbmat$ and thus associates Eq.~\eqref{eq:romEq} with the
necessary optimality conditions for a minimum-residual problem; see
Ref.~\cite{carlbergGalDiscOpt} for a detailed comparison of the two approaches
in the context of nonlinear dynamical systems.
When the residual operator is nonlinear in its first argument or
nonaffine
 in (functions of) its second argument, additional `hyper-reduction' techniques
must be employed in order to ensure that the ROM equations~\eqref{eq:romEq} can
be solved with a computational cost that is independent of the FOM
dimension $\ndof$. Such techniques include the empirical
 interpolation method (EIM) \cite{Barrault2004,maday2007general,GreplMaday:EIM}, its discrete
 variant  DEIM \cite{Chaturantabut2010,AHS_CECAM},
 gappy POD \cite{willcox2006unsteady,carlberg2013gnat}, and missing point estimation
 \cite{astrid2008missing}.

Finally, we denote the ROM-predicted output as $
\qoiROM:\paramDummy\mapsto\qoiFunc(\stateROM(\paramDummy);\paramDummy) $.
%
%
%
%
Critically, because the ROM solution approximates the FOM solution, the ROM
will generally introduce both a state error and a quantity-of-interest error
\begin{align}\label{eq:stateError}
		\stateError&:\paramDummy\mapsto\state(\paramDummy) -
		\stateROM(\paramDummy) \quad\text{and}\quad\qoiError :\paramDummy\mapsto\qoiFOM(\paramDummy) - \qoiROM(\paramDummy),
\end{align}
respectively,
with
$\stateError:\paramDomain\rightarrow\RR{\ndof}
$
and $\qoiError:\paramDomain\rightarrow\RR{\nqoi}.$

\section{State error decomposition and approximation}\label{sec:ROMESstate}
The objective of this work is to compute statistical models of the state error
$\stateError$ and quantity-of-interest error $\qoiError$ in a manner that does
not require identifying the quantities of interest \textit{a priori}. We now
present the mathematical framework that will be leveraged by the proposed
method, which is presented in Section \ref{sec:GPstate}. In particular, this
section (1) decomposes the state error
$\stateError$ into in-plane and out-of-plane errors
(Section \ref{sec:decomposition}),
(2) identifies low-dimensional subspaces for each of these error components
(Section \ref{sec:lowdimapprox}), (3) derives first-order estimates of
the generalized coordinates characterizing these error
components (Section \ref{sec:firstorderapprox}), and
(4) applies model reduction to inexpensively approximate these estimates
(Section \ref{sec:romDWR}).

\subsection{State error decomposition}\label{sec:decomposition}

We begin by defining the in-plane projector, which is the linear operator that
computes the orthogonal projection onto the linear part
 $ \trialsubspace$
of the affine trial subspace, as
\begin{equation} \label{eq:projector}
	\projectorParallel\defeq
	\rbmat(\rbmat^T\metricmat\rbmat)^{-1}\rbmat^T\metricmat.
\end{equation}
Here, $\metricmat\in\spd{\ndof}$ is a matrix defining an inner product
(e.g., the discrete counterpart to
the inner product characterizing a Sobolev space)
with $(\stateDummy,\stateDummyOpt)_\metricmat\defeq
\stateDummy^T\metricmat\stateDummyOpt$ and
$\|\stateDummy\|_\metricmat\defeq\sqrt{(\stateDummy,\stateDummy)_\metricmat}$,
and $\spd{N}$ denotes the set of $N\times N$
symmetric-positive-definite matrices. The in-plane projector $\projectorParallel$ inherits
standard properties of orthogonal projectors, i.e., optimality
\begin{equation}
	\projectorParallel \stateDummy  =
	\underset{\stateDummyOpt\in
	\trialsubspace}{\argmin}\| \stateDummy-\stateDummyOpt\|_\metricmat;
\end{equation}
orthogonality
$(\stateDummy-\projectorParallel \stateDummy ,\stateDummyOpt)_\metricmat=0$,
$\forall\stateDummy\in\RR{\ndof}$,
$\forall\stateDummyOpt\in\trialsubspace$;
and idempotency,
$(\projectorParallel)^2 =
\projectorParallel $.

\begin{figure}[h!]
\centering
\includegraphics[width=0.6\textwidth]{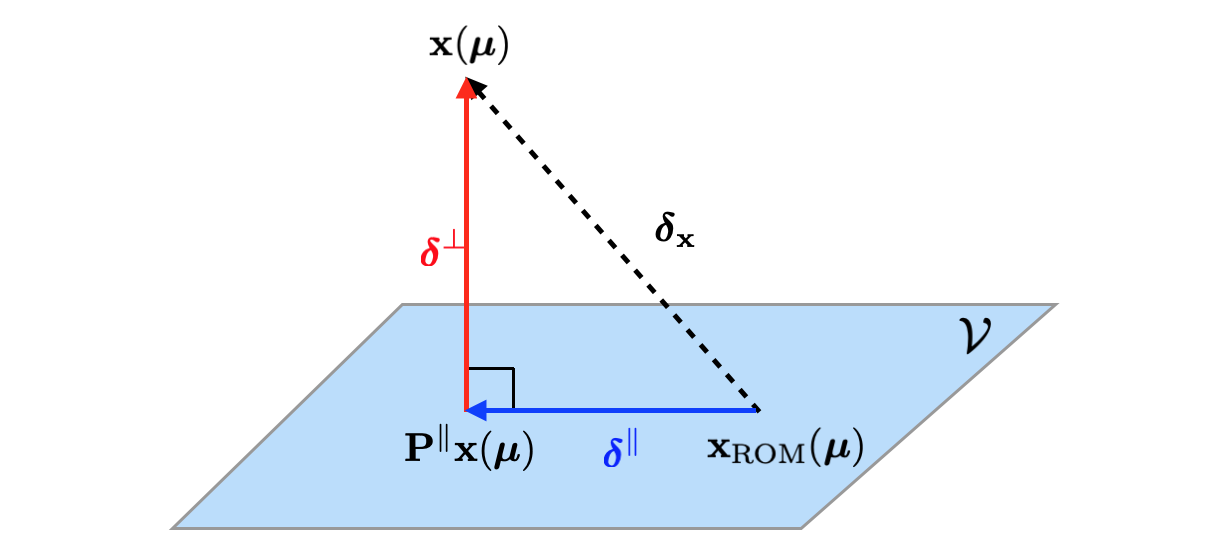}
 \caption{Graphical depiction of the decomposition of the state error
$\stateError\in\RR{\ndof}$ into the in-plane error $\stateErrorInplane\in\trialsubspace$ and the out-of-plane
error $\stateErrorOutofplane\in\trialsubspaceperp$. }
 \label{fig:projection_new}
\end{figure}

The projector enables the state error to be decomposed as
\begin{equation}\label{eq:stateErrorDecomp}
\stateError = \stateErrorInplane +
\stateErrorOutofplane,
\end{equation}
where the \textit{in-plane error} lies within the
linear part $\trialsubspace$ of the affine trial subspace and is defined as
\begin{align} \label{eq:inplaneerror}
\begin{split}
	\stateErrorInplane&:\paramDummy\mapsto\projectorParallel \stateError(\paramDummy) \\
&:\paramDomain\rightarrow\trialsubspace,
\end{split}
\end{align}
and the \textit{out-of-plane error} is orthogonal to the linear subspace $\trialsubspace$  and is
defined as
\begin{align} \label{eq:outofplaneerror}
\begin{split}
	\stateErrorOutofplane&:\paramDummy\mapsto\stateError(\paramDummy)-\projectorParallel \stateError(\paramDummy) \\
&:\paramDomain\rightarrow\trialsubspaceperp.
\end{split}
\end{align}
We note that $(\stateErrorOutofplane(\param),\stateDummyOpt)_\metricmat=0$,
$\forall\stateDummyOpt\in\trialsubspace$, $\forall\param\in\paramDomain$, and
$(\stateErrorOutofplane(\param),\stateErrorInplane(\param))_\metricmat=0$, $\forall\param\in\paramDomain$.
Figure \ref{fig:projection_new} depicts this error decomposition graphically.

We note that the in-plane error $\stateErrorInplane(\param)$ can be
interpreted as the closure error, as the in-plane
error expresses the error in the `preserved variables' (i.e., the solution
component in the affine trial subspace $\stateRef +
 \trialsubspace$) incurred by solving equations that omit the
`neglected variables' (i.e., the solution component in $\trialsubspaceperp$).
\begin{remark}[Necessary conditions for zero in-plane error]
	The in-plane (i.e., closure) error $\stateErrorInplane(\param)$ is zero if the
	residual $\res$ is linear in its first argument such that
	\eqref{eq:reslinear} holds and either: (1) Galerkin projection is employed
	(i.e., $\testrbmat = \rbmat$), the system matrix is symmetric and positive
	definite (i.e., $\Amat(\param)\in\spd{\ndof}$), and the chosen metric is
	equal to the system matrix (i.e., $\metricmat = \Amat(\param)$, or (2)
	least-squares Petrov--Galerkin (LSPG) projection is employed (i.e.,
	$\testrbmat(\stateRed;\param) = \Amat(\param)\rbmat$) and the chosen metric is equal to the associated
	normal-equations matrix (i.e., $\metricmat = \Amat(\param)^T\Amat(\param)$).
\end{remark}

\subsection{Low-dimensional representations of the state-error
components}\label{sec:lowdimapprox}

We aim to construct statistical models for both the in-plane error
$\stateErrorInplane$ and
out-of-plane error $\stateErrorOutofplane$; the former can be considered a
statistical closure model. However, the offline cost of
constructing a model for each of the $\ndof$ elements of these error vectors
is computationally costly, and the online complexity of evaluating these models is
$\ndof$-dependent. To mitigate this cost, we instead aim to construct
a statistical model for \textit{generalized coordinates} representing these
errors in low-dimensional subspaces.

This is a straightforward task for the in-plane error, as
$\stateErrorInplane(\param)\in\trialsubspace$ with
$\dim(\trialsubspace) = \nrb\ll\ndof $ by construction, and thus
\begin{equation}\label{eq:inplaneerrorDecomp}
	\stateErrorInplane(\param) =
	\rbmat\stateErrorInplaneRed(\param),
\end{equation}
where $\stateErrorInplaneRed:\paramDomain\rightarrow\RR{\nrb}$ denote the
generalized coordinates of the in-plane error that---from the definition of the in-plane projector
\eqref{eq:projector},
in-plane error \eqref{eq:inplaneerror}, and associated decomposition
\eqref{eq:inplaneerrorDecomp}---satisfies
\begin{equation}\label{eq:stateErrorInplaneRedExpression}
\stateErrorInplaneRed:\paramDummy\mapsto
(\rbmat^T\metricmat\rbmat)^{-1}\rbmat^T\metricmat\stateError(\paramDummy).
\end{equation}

In contrast, the out-of-plane error satisfies
$\stateErrorOutofplane(\param)\in\trialsubspaceperp$ with
$\dim(\trialsubspaceperp) = \ndof-\nrb$. Because $\nrb\ll\ndof$, the subspace
$\trialsubspaceperp$ is high-dimensional. To address this, we approximate the
out-of-plane error
in an $\nrbperp$-dimensional (with $\nrbperp\ll\ndof$) linear subspace as
\begin{equation} \label{eq:stateErrorOutofplaneApprox}
\stateErrorOutofplane(\param)\approx
\rbmatperp\stateErrorOutofplaneRed(\param)\in\trialsubspaceperpRed\subset\trialsubspaceperp
\end{equation}
with
\begin{equation} \label{eq:stateErrorOutofplaneRed}
	\stateErrorOutofplaneRed:\paramDummy
	\mapsto\projectorPerp\stateErrorOutofplane(\paramDummy),
\end{equation}
where
\begin{equation} \label{eq:projectorPerp}
	\projectorPerp\defeq
\rbmatperp(\rbmatperpT\metricmat\rbmatperp)^{-1}\rbmatperpT\metricmat
\end{equation}
is the orthogonal projection onto the linear subspace $\trialsubspaceperpRed$;
$\rbmatperp\in\RRstar{\ndof\times\nrbperp}$ denotes the out-of-plane error
basis matrix such that $\trialsubspaceperpRed\defeq\range{\rbmatperp}$; and
$\stateErrorOutofplaneRed:\paramDomain\rightarrow\RR{\nrbperp}$ denotes
the generalized coordinates of the out-of-plane error that---from the
definition of the  out-of-plane error \eqref{eq:outofplaneerror},
the associated approximation \eqref{eq:stateErrorOutofplaneApprox},
generalized-coordinate definition \eqref{eq:stateErrorOutofplaneRed},
and out-of-plane projector
\eqref{eq:projectorPerp}---satisfies
\begin{equation} \label{eq:stateErrorOutofplaneredExpression}
\stateErrorOutofplaneRed(\param) =
(\rbmatperpT\metricmat\rbmatperp)^{-1}\rbmatperpT\metricmat\stateError(\param),
\end{equation}
where
we have used $\rbmat^T\metricmat\rbmatperp = \zerovec$.

Comparing Eqs.~\eqref{eq:stateErrorInplaneRedExpression} and
\eqref{eq:stateErrorOutofplaneredExpression} and using
$\rbmat^T\metricmat\rbmatperp = \zerovec$ allows the definition
of the (in-plane and out-of-plane) error generalized coordinates
\begin{align}\label{eq:stateErrorRedDef}
	\stateErrorRed:\paramDummy\mapsto
	\begin{bmatrix}
		\stateErrorInplaneRed(\paramDummy)\\
\stateErrorOutofplaneRed(\paramDummy)
	\end{bmatrix}
	=
\projmat\stateError(\param)
\end{align}
where
\begin{equation}
	\projmat\defeq(\rbmatTot^T\metricmat\rbmatTot)^{-1}\rbmatTot^T\metricmat\in\RR{\nrbtot\times\ndof},\quad
	\rbmatTot\defeq\begin{bmatrix}
		\rbmat & \rbmatperp
	\end{bmatrix}\in\RRstar{\ndof\times\nrbtot}
\end{equation}
and $\nrbtot\defeq\nrb+\nrbperp$.

\begin{remark}[Out-of-plane basis matrix
	construction]\label{rem:outofplanemat}
The out-of-plane basis matrix $\rbmatperp$ can be constructed by a variety of
means. For example, if $\rbmat$ corresponds to a truncated proper orthogonal
decomposition (POD) basis, then $\rbmatperp$ can be set to the (discarded)
$\nrb+1$ to $\nrbtot$
POD
modes; this idea has also been
employed in the context of ROM error estimation \cite{wirtzDeim}. Alternatively, the basis
can be constructed by computing the projection error of FOM solutions over a
parameter set $\param\in\paramOutofplane\subset\paramDomain$ such that
$\trialsubspaceperpRed\subseteq\Span{\state(\param) - \stateRef -
\projectorParallel \state(\param) }_{\param\in\paramOutofplane}$.
\end{remark}

\subsection{Dual-weighted-residual error estimation}\label{sec:firstorderapprox}

We now derive first-order approximations for the in-plane-error generalized
coordinates $\stateErrorInplaneRed$ and out-of-plane-error general coordinates
$\stateErrorOutofplaneRed$.
		Assuming the residual is twice continuously differentiable, we can
		approximate the residual of the FOM solution to first order about the
		residual of the ROM solution as
\begin{equation}
	\zerovec = \res(\state(\param);\param) = \res(\stateROM(\param);\param) +
	\jacobian (\stateROM(\param);\param) \stateError(\param)
	+ \order( \|\stateError(\param)\|^2  ),\quad \text{as}\
	\|\stateError(\param)\|\rightarrow 0.
\end{equation}
Note that the high-order term is zero if the residual is linear in its first
argument, i.e., if \eqref{eq:reslinear} holds. Solving for the state error yields
\begin{equation}\label{eq:stateErrorFirstOrder}
	\stateError(\param)= -\left[\jacobian (\stateROM(\param);\param)\right]^{-1}\res(\stateROM(\param);\param) +
	\order( \|\stateError(\param)\|^2  ),\quad \text{as}\
	\|\stateError(\param)\|\rightarrow 0.
\end{equation}
Substituting
Eq.~\eqref{eq:stateErrorFirstOrder} in Eq.~\eqref{eq:stateErrorRedDef}
yields
\begin{equation}
\stateErrorRed(\param) =
-\projmat\left[\jacobian
(\stateROM(\param);\param)\right]^{-1}\res(\stateROM(\param);\param)
+
	\order( \|\stateError(\param)\|^2  ),\quad \text{as}\
	\|\stateError(\param)\|\rightarrow 0.
\end{equation}

Defining the $i$th dual
$\dualArg{i}:\paramDomain\rightarrow\RR{\ndof}$, $i=1,\ldots,\nrbtot$
as the solution to the $\ndof$-dimensional system of linear equations
\begin{equation} \label{eq:dualinplane}
\left[\jacobian
(\stateROM(\param);\param)\right]^T\dualArg{i}(\param)= -
\projmat^T\unitvecArg{i},\quad i=1,\ldots,\nrbtot,
\end{equation}
where $\unitvecArg{i}\in\{0,1\}^{\ndof}$ denotes the $i$th canonical unit
vector, we can express  the $i$th
error generalized coordinate as
\begin{equation}\label{eq:stateErrorInplaneRedExpressionTwo}
	\stateErrorRedArg{i}(\param)=
	\dualArg{i}(\param)^T\res(\stateROM(\param);\param)+
	\order( \|\stateError(\param)\|^2  ),\quad \text{as}\
	\|\stateError(\param)\|\rightarrow 0,\quad i=1,\ldots,\nrbtot.
\end{equation}

\subsection{Reduced-order model approximation to dual-weighted-residual error
estimates}\label{sec:romDWR}

Dual problems \eqref{eq:dualinplane}
are linear, even if the original problem
\eqref{eq:discretizedPB} is nonlinear; however, their dimension $\ndof$
remains the same as that of the full-order model. Thus, employing the
associated dual vectors for \textit{a posteriori} error modeling as suggested
by Eq.~\eqref{eq:stateErrorInplaneRedExpressionTwo} is computationally expensive.

To mitigate this cost, we propose to approximate these duals via
model reduction in analogue to the approach described in Section \ref{sec:ROM}
for approximating the state. First, we approximate the duals as
$\dualArg{i}\approx\dualApproxArg{i}$, $i=1,\ldots,\nrbtot$, where
\begin{align}\label{eq:dualApproxArgDef}
	&\dualApproxArg{i}(\param)=\rbmatDualArg{i} \dualRedArg{i}(\param) ,\quad
	i=1,\ldots,\nrbtot
\end{align}
where $\rbmatDualArg{i}\in\RRstar{\ndof\times\nrbDualArg{i}}$ denote the dual trial-basis matrices,
$\dualRedArg{i}:\paramDomain\rightarrow\RR{\nrbDualArg{i}}$, and $\nrbDualArg{i}\ll\ndof$ for
$i=1,\ldots, \nrbtot$.
As in the case of the trial-basis matrix $\rbmat$, the trial-basis matrices $\rbmatDualArg{i}$, $i=1,\ldots,\nrbtot$ can be constructed by
a variety
of means,
e.g., the  reduced-basis method, POD.
We then substitute
$\dualArg{i}\leftarrow\dualApproxArg{i}$
in Eqs.~\eqref{eq:dualinplane}
and enforce orthogonality of the residual to
the range of
associated test basis matrices
$\testrbmatDualArg{i}(\param)\in\RRstar{\ndof\times\nrbDualArg{i}}$, $i=1,\ldots,\nrbtot$
to obtain the ROM systems of equations
\begin{align}
	\label{eq:dualinplaneROM}
	&[\testrbmatDualArg{i}(\param)]^T\left[\jacobian
(\stateROM(\param);\param)\right]^T\rbmatDualArg{i}\dualRedArg{i}(\param)= -
[\testrbmatDualArg{i}(\param)]^T\metricmat\rbmat(\rbmat^T\metricmat\rbmat)^{-1}\unitvecArg{i}, 
\end{align}
for $i=1,\ldots,\nrbtot$, whose solutions define the generalized coordinates $\dualRedArg{i}$, $i=1,\ldots,\nrbtot$.

As before, a Galerkin projection corresponds to
$\testrbmatDualArg{i} = \rbmatDualArg{i}$, $i=1,\ldots,\nrbtot$,
while an LSPG projection corresponds to
$\testrbmatDualArg{i}(\param) = \left[\jacobian
(\stateROM(\param);\param)\right]^T\rbmatDualArg{i}$
$i=1,\ldots,\nrbtot$. Again, if the residual operator
is nonlinear in the state or
nonaffine in functions of the parameter inputs, then hyper-reduction is required to
ensure the cost of assembling the linear systems in Eqs.~\eqref{eq:dualinplaneROM}
does not scale with the dimension $\ndof$.

Now, substituting $\dualArg{i}\leftarrow\dualApproxArg{i}$
Eqs.~\eqref{eq:stateErrorInplaneRedExpressionTwo}
and ignoring high-order terms yields
cheaply computable approximations to the in-plane and out-of-plane error generalized
coordinates
\begin{align}
	\label{eq:stateErrorInplaneRedExpressionTwoCheap}
	&\stateErrorRedArg{i}(\param) \approx
	[\dualApproxArg{i}(\param)]^T\res(\stateROM(\param);\param),
	\quad
	i=1,\ldots,\nrbtot.
\end{align}
Note that the approximation is induced by the use of a model reduction to
approximate the duals, as well as truncation error in the case of
nonlinear FOM equations \eqref{eq:discretizedPB}.
The next section describes how this approximation to the error generalized
coordinates can be used to construct a statistical model of the state error.

\begin{remark}[Unique v.\ shared dual bases]\label{rem:uniqueShared}
The strategy outlined above for applying model reduction to the dual problems
requires the construction of $\nrbtot$ dual trial-basis matrices and dual
test-basis matrices $\rbmatDualArg{i}$ and $\testrbmatDualArg{i}$,
$i=1,\ldots,\nrbtot$, respectively.
If each of these basis matrices is unique, then the cost of solving
Eqs.~\eqref{eq:dualinplaneROM} is approximately $2/3\sum_{i=1}^{\nrbtot}
\nrbDualArg{i}^3+ 2\sum_{i=1}^{\nrbtot}
\nrbDualArg{i}^2$; this cost is small if each dual trial-basis matrix dimension
$\nrbDualArg{i}$ is small. However, each of the basis matrices must be trained
independently; in the event of limited training, these basis matrices
individually may be too low-dimensional to generate accurate dual
approximations $\dualApproxArg{i}$, which can lead to large approximation errors.

Alternatively, one may employ a single `shared' dual trial-basis matrix
$\rbmatDual\in\RRstar{\ndof\times\nrbDual}$ and test-basis matrix $\testrbmatDual\in\RRstar{\ndof\times\nrbDual}$ such that
$\rbmatDualArg{i}=\rbmatDual$ and $\testrbmatDualArg{i}=\testrbmatDual$,
$i=1,\ldots,\nrbtot$. In this case---because each of the linear systems
\eqref{eq:dualinplaneROM} is characterized by the same system matrix---the cost
of solving the resulting systems is approximately $2/3\nrbDual^3 + 2\nrbtot\nrbDual^2$. In many cases, this cost is significant,
as the dimension $\nrbDual$ is typically large, due to the fact that the basis is constructed from jointly training all duals; in the worst case, if $\nrbDual = \sum_{i=1}^\nrbtot\nrbDualArg{i}$, then the cost is
$2/3(\sum_{i=1}^\nrbtot\nrbDualArg{i})^3 + 2\nrbtot(\sum_{i=1}^\nrbtot\nrbDualArg{i})^2$.
On the other hand, this approach often requires less training
to compute a trial-basis matrix with good approximation properties, as
information across all dual solutions informs the basis.
\end{remark}


\section{ROMES error models}\label{sec:GPstate}

We now leverage the framework
presented in Section \ref{sec:ROMESstate} to describe the application of the
ROMES method \cite{kevin:romes} to construct statistical models of the
in-plane and out-of-plane error generalized coordinates using indicators
corresponding to the approximated dual-weighted residuals.
Section
\ref{sec:statisticalmodel} describes the formulation for the statistical
model, Section \ref{sec:errorindicator} describes the error indicator (i.e.,
feature) we employ, Section \ref{sec:GP} provides an overview of
Gaussian-process regression, which is the technique we employ to construct the
statistical model, Section \ref{seq:stateqoimodels} describes the application
of
ROMES error models to obtain statistical models for the state and
quantities-of-interest errors, and Section \ref{sec:offlineonline} describes
the offline/online computational strategy employed to realize the method in
practice.

\subsection{Statistical model}\label{sec:statisticalmodel}

Our objective is to construct a low-dimensional, statistical model of the
high-dimensional, deterministic, and generally unknown ROM error. The
probability distribution of the random variable representing the ROM error reflects the
epistemic uncertainty about its value.
Define a probability space $(\samplespace, \mathcal F, P)$.  We aim to approximate the
high-dimensional deterministic mappings
\begin{align}
\begin{split}
	&\mappingArg{i}:\paramDummy\mapsto\stateErrorRedArg{i}(\paramDummy),\quad
i=1,\ldots,\nrbtot
\end{split}
\end{align}
with $\mappingArg{i}:\paramDomain\rightarrow\RR{}$  and
$\nparam$ possibly large, by univariate stochastic mappings
\begin{align}\label{eq:ROMESmodels}
\begin{split}
	&\mappingApproxArg{i}:\indicatorArg{i}(\paramDummy)\mapsto\stateErrorRedModelArg{i}(\paramDummy),
\quad i=1,\ldots,\nrbtot
\end{split}
\end{align}
respectively, where
$\indicatorArg{i}:\paramDomain\rightarrow\RR{}$, $i=1,\ldots,\nrbtot$  denote
\textit{error indicators} and
$\stateErrorRedModelArg{i}(\paramDummy):\samplespace\rightarrow\RR{}$,
$i=1,\ldots,\nrbtot$ denote random variables for the error generalized
coordinates. The stochastic mapping should satisfy the following desiderata
(see Refs.~\cite{kevin:romes,frenoCarlberg}):
\begin{enumerate}
	\item \label{prop:cheap}the error indicators are
\textit{cheaply computable} given any
$\param\in\paramDomain$;
\item \label{prop:lowvar}the stochastic mappings exhibit
\textit{low variance}, i.e.,
$\expectation{(\mappingApproxArg{i}(\indicatorArg{i}(\param))-\expectation{\mappingApproxArg{i}(\indicatorArg{i}(\param))})^2}$
is `small' for all $\param\in\paramDomain$ (this ensures the ROM-induced
epistemic uncertainty is small); and
\item\label{prop:validated} the stochatic mappings are
\textit{validated}, i.e., the (empirical) distribution of test data is `close'
to the (reference) distribution prescribed by the stochastic mappings
(using, e.g., prediction intervals, the Komolgorov--Smirnov test).
\end{enumerate}

We now describe choices of error indicators and stochastic-mapping methods
that lead to statistical models satisfying the above conditions.

\subsection{Error indicators}\label{sec:errorindicator}

The error indicator should be selected so that it is both cheaply computable
(Condition \ref{prop:cheap} above) and can lead to a low-variance stochastic
mapping (Condition \ref{prop:lowvar} above); the latter condition implies that
the error indicator should be informative of the error such that the mean of
the stochastic mapping can explain most of the variance in the observed error.

Inspired by the analysis of Section \ref{sec:ROMESstate}, and
Eq.~\eqref{eq:stateErrorInplaneRedExpressionTwoCheap} in particular, we
propose employing the approximated dual-weighted residual as an error
indicator, i.e.,
\begin{equation}\label{eq:indicatorDef}
	\indicatorArg{i}(\param)=
	[\dualApproxArg{i}(\param)]^T\res(\stateROM(\param);\param),
	\quad
	i=1,\ldots,\nrbtot.
\end{equation}
From Eq.~\eqref{eq:stateErrorInplaneRedExpressionTwoCheap}, we can see that
$\stateErrorRedArg{i}(\param)\approx\indicatorArg{i}(\param)$, where the
approximation arises both to the use of model reduction to approximate the
duals and truncation error when the residual is nonlinear in the state.

\subsection{Gaussian-process regression}\label{sec:GP}

As in Ref.~\cite{kevin:romes}, we propose to construct the stochastic mappings
$\mappingApproxArg{i}$, $i=1,\ldots,\nrbtot$ using Gaussian process (GP)
kernel regression \cite{rasmussen}, which is a supervised machine learning
method. We first provide a brief review of this technique. A GP is a collection of random variables such that
any finite number of them has a joint Gaussian distribution. GP kernel
regression computes this GP by Bayesian inference using a kernel function
and training data
$\trainingSet =
\{(\featureTrainArg{i},\responseTrainArg{i})\}_{i=1}^{\ntrain}$, where
$\featureTrainArg{i}\in\RR{\nfeature}$ and
$\responseTrainArg{i}\in\RR{}$ denote the $i$th instance of the
features and response, respectively. We consider a single prediction
point characterized by features $\featurePredict\in\RR{\nfeature}$, as we
treat all predictions as arising from independent samples of the GP.
First, the approach sets the prior distribution to
\begin{equation}\label{eq:prior}
	\tilde y_\text{prior}(\featureAllSet)\sim
	\normal
		{{\vondermonde \regweights}}{\kernelMatFunc{\featureAllSet}{\featureAllSet} +
		\variance\identity}.
\end{equation}
Here, $\featureAllSet\defeq[\featureAllArg{1}\
\cdots\ \featureAllArg{\ntrain+1}]^T\in\RR{(\ntrain+1)\times\nfeature}$ with
$\featureAllArg{i} = \featureTrainArg{i}$, $i=1,\ldots,\ntrain$ and
$\featureAllArg{\ntrain+1} = \featurePredict$;
element $(i,j)$ of the matrix
$\vondermonde\in\RR{(\ntrain+1)\times\nbasis}$ is
$[\vondermonde]_{ij}\defeq\basisArg{j}(\featureAllArg{i})$
with $\basisArg{j}:\RR{\nfeature}\rightarrow\RR{}$, $j=1,\ldots,\nbasis$ denoting the considered basis
functions (e.g., polynomials);
$\regweights\in\RR{\nbasis}$ denotes the basis-expansion coefficients;
and $\variance$ denotes the additive
noise arising from the non-uniqueness of the mapping from the features to the response.
Element $(i,j)$ of the kernel matrix
$\kernelMatFunc{\generalSetOne}{\generalSetTwo}$
with $\generalSetOne\defeq[\generalSetOneVecArg{1}\ \cdots\
\generalSetOneVecArg{\ngeneralSetOne}]$
and $\generalSetTwo\defeq[\generalSetTwoVecArg{1}\ \cdots\
\generalSetTwoVecArg{\ngeneralSetTwo}]$
is
\begin{equation}
	[\kernelMatFunc{\generalSetOne}{\generalSetTwo}]_{ij} \defeq
	\kernelFunc{\generalSetOneVecArg{i}}{\generalSetTwoVecArg{j}}, \quad
	i=1,\ldots,\ngeneralSetOne,\ j=1,\ldots,\ngeneralSetTwo
\end{equation}
and many choices of the kernel function $\kernelFuncNo$ exist. The kernel
function is typically characterized by its own hyperparameters
$\hyperparamsKernel$, e.g., the
length scale in the case of the squared exponential kernel. Given the
training data $\trainingSet$ and fixed values of the coefficients
$\regweights$, noise variance $\variance$, and
kernel hyperparameters $\hyperparamsKernel$,
the prediction corresponds to a random variable with
posterior distribution
\begin{equation}
	\tilde y(\featurePredictSet;\hyperparams)\sim
\normal
{\overallmean{\featurePredictSet}}{\overallVar{\featurePredictSet}}
\end{equation}
with
\begin{align}
	&\overallmean{\featurePredictSet} \defeq
	\overallmeanDef{\featurePredictSet} + {[\basisArg{1}(\featurePredict)\
	\cdots\ \basisArg{\nbasis}(\featurePredict)] \regweights}\\
	&\overallVar{\featurePredictSet} \defeq
	\overallVarDef{\featurePredictSet},
\end{align}
where $\responseTrainVec\defeq\left[\responseTrainArg{1}\ \cdots\
\responseTrainArg{\ntrain}\right]^T$ and
$\featureTrainSet\defeq[\featureTrainArg{1}\ \cdots\
\featureTrainArg{\ntrain}]$. Indeed, computing the posterior
using
the training data is a simple operation derived from
conditioning a joint Gaussian distribution\footnote{In this
respect, consider the fundamental result
\[
\begin{pmatrix}
{\bf v} \\ {\bf w}
\end{pmatrix} \sim
\mathcal{N}
\left(
\begin{pmatrix}
\boldsymbol{\nu}_v \\ \boldsymbol{\nu}_w
\end{pmatrix}
 ,
\begin{pmatrix}
 \boldsymbol{\Sigma}_{vv} &  \boldsymbol{\Sigma}_{vw} \\  \boldsymbol{\Sigma}_{vw}^T &  \boldsymbol{\Sigma}_{ww}
\end{pmatrix}
\right) \  \  \Rightarrow \ \
{\bf v} | {\bf w} \sim \mathcal{N} (\boldsymbol{\nu}_{v|w},
\boldsymbol{\Sigma}_{v|w} ),
\]
where $\boldsymbol{\nu}_{v|w} = \boldsymbol{\nu}_v + \boldsymbol{\Sigma}_{vw}
\boldsymbol{\Sigma}_{ww}^{-1} ({\bf w} - \boldsymbol{\nu}_w  )$  and
$\boldsymbol{\Sigma}_{v|w} = \boldsymbol{\Sigma}_{vv} -
\boldsymbol{\Sigma}_{vw}  \boldsymbol{\Sigma}_{ww}^{-1}
\boldsymbol{\Sigma}_{vw}^T$.}.  The parameters $\regweights$, $\variance$, and
hyperparameters characterizing the kernel $\hyperparamsKernel$ can be set in a
variety of ways, e.g., via maximum likelihood estimation (as in
Ref.~\cite{kevin:romes}), cross-validation.

\subsection{Gaussian-process ingredients and cross-validation for ROMES}\label{sec:GPROMES}

In this work, we specify the GP ingredients as follows.  Following the
ROMES method \cite{kevin:romes} and the presentation of Section
\ref{sec:statisticalmodel},  we propose to apply GP regression to
(independently) construct each of the mappings $\mappingApproxArg{i}$,
$i=1,\ldots,\nrbtot$, wherein the feature corresponds to the prescribed error
indicator (i.e., $\feature = \indicatorArg{i}$ with $\nfeature=1$), and the
response corresponds to the error generalized coordinate (i.e., $\response =
\stateErrorRedArg{i}$).
Eqs.~\eqref{eq:stateErrorInplaneRedExpressionTwoCheap} and
\eqref{eq:indicatorDef} illustrate that the relationship between the features and
the response is approximately linear in this case; it is exactly linear if the
approximated dual is exact (i.e., $\dualArg{i}=\dualApproxArg{i}$) and
the when the residual operator $\res$ is linear in its first argument (i.e.,
Eq.~\eqref{eq:reslinear} holds). Thus, we select the basis functions in the
mean of the prior distribution \eqref{eq:prior} to enable linear responses, i.e.,
$\basisArg{1}:\indicatorArg{i}\mapsto
1$, $\basisArg{2}:\indicatorArg{i}\mapsto
\indicatorArg{i}$ with $\nbasis = 2$.

For training the ROMES models, we employ training data comprising
indicator--error pairs computed at ROMES-training parameter instances
$\paramROMES\subset\paramDomain$ with $\card{\paramROMES} = \ntrain$,
i.e., the training data for model $\mappingApproxArg{i}$ corresponds to
\begin{align}
	&\ROMEStrainingArg{i}\defeq\{(\indicatorArg{i}(\param),\stateErrorRedArg{i}(\param))\}_{\param\in\paramROMES},\quad
	i=1,\ldots,\nrbtot.
\end{align}
Recall that the error generalized coordinates $\stateErrorRed(\param)$ can be
computed from Eq.~\eqref{eq:stateErrorRedDef}; this expression requires
computing the state error $\stateError(\param)$, which in turn requires
computing both the FOM state $\state(\param)$ and ROM state
$\stateROM(\param)$.

Given training data $\ROMEStrainingArg{i}$, we train
model $\mappingApproxArg{i}$ as follows: we determine the hyperparmeters
$\hyperparamsArg{i}\equiv(\varianceArg{i},\hyperparamsKernelArg{i})$ using $K$-fold cross validation \textit{with
specialized loss functions} that target different interpretations of
statistical validation (condition \ref{prop:validated} in Section
\ref{sec:statisticalmodel}), and we determine the coefficients
$\regweightsArg{i}$
using maximum likelihood estimation.
In particular, we first
separate the ROMES training data $\ROMEStrainingArg{i}$ into $K$
		non-overlapping subsets $\ROMEStrainingArgs{i}{j}$, $j=1,\ldots, K$ such
		that $\ROMEStrainingArg{i} = \cup_{j=1}^K\ROMEStrainingArgs{i}{j}$ and
		$\ROMEStrainingArgs{i}{j}\cap\ROMEStrainingArgs{i}{k}=\emptyset$,
		$j\neq k$. We then define a set of candidate hyperparameter values
		$\hyperparamsSet$.
For each candidate value of the hyperparameters
$\hyperparams\equiv(\variance,\hyperparamsKernel)\in\hyperparamsSet$, we compute the values of basis-expansion
coefficients using maximum likelihood estimation as
\begin{equation} \label{eq:regMLE}
	\regweights_{i,j}(\hyperparams) =
	\left(\vondermondeROMES^T(\kernelMatFunc{\indicatorVecArgs{i}{j}}{\indicatorVecArgs{i}{j}}+\variance\identity)^{-1}\vondermondeROMES\right)^{-1}
	\vondermondeROMES^T(\kernelMatFunc{\indicatorVecArgs{i}{j}}{\indicatorVecArgs{i}{j}}+\variance\identity)^{-1}\responseTrainVecROMES,
\end{equation}
where $\indicatorVecArgs{i}{j}\in\RR{\card{\ROMEStrainingArg{i}\setminus\ROMEStrainingArgs{i}{j}}}$ and
$\stateErrorRedArgs{i}{j}\in\RR{\card{\ROMEStrainingArg{i}\setminus\ROMEStrainingArgs{i}{j}}}$ denote the
vectorized features and responses associated with training set
$\ROMEStrainingArg{i}\setminus\ROMEStrainingArgs{i}{j}$, and $\mathbf{1}$ denotes a
vector of ones. Note that we have suppressed the dependence of the kernel
matrix $\kernelMat$ on the hyperparameters $\hyperparamsKernel$ for
notational simplicity. The values $\hyperparams$
and
	$\regweights_{i,j}(\hyperparams)$ define a candidate ROMES model
$\mappingApproxArgsNo{i}{j}$ characterized by
\begin{equation}\label{eq:stateErrorRedModelArgDefCross}
	\mappingApproxArgs{i}{j}{\indicatorArg{i}(\param)}{\hyperparams}=
	\stateErrorRedModelArgs{i}{j}{\param}{\hyperparams}
	\sim\normal{
	\overallmeanArg{i,j}{\indicatorArg{i}(\param);\hyperparams}
	}{\overallVarArg{i,j}{\indicatorArg{i}(\param);\hyperparams}},\quad
	i=1,\ldots,\nrbtot
\end{equation}
with
\begin{align}
	& \overallmeanArg{i,j}{\indicator;\hyperparams}\defeq
	\kernelMatFunc{\indicator}{\indicatorVecArgs{i}{j}}
(\kernelMatFunc{\indicatorVecArgs{i}{j}}{\indicatorVecArgs{i}{j}}+\variance\identity)^{-1}\responseTrainVecROMES
+
	[1\ \indicator]\regweights_{i,j}(\hyperparams)
	\\
	&\overallVarArg{i,j}{\indicator;\hyperparams} \defeq
	\kernelMatFunc{\indicator}{\indicator} -
\kernelMatFunc{\indicator}{\indicatorVecArgs{i}{j}}
(\kernelMatFunc{\indicatorVecArgs{i}{j}}{\indicatorVecArgs{i}{j}}+\variance\identity)^{-1}
\kernelMatFunc{\indicatorVecArgs{i}{j}}{\indicator} + \variance.
\end{align}

Subsequently, each candidate value of the
hyperparameters $\hyperparams\in\hyperparamsSet$ is assigned a loss
$\loss{i}{\hyperparams}$ with
\begin{equation} \label{eq:lossi}
\loss{i}{\hyperparams} \equiv
\frac{1}{K}\sum_{j=1}^K\lossArg{\hyperparams}{i}{j},
\end{equation}
where $\lossArg{\hyperparams}{i}{j}$ denotes the loss for the $i$th ROMES
model on the $j$th
validation set corresponding to hyperparameters $\hyperparams$.
We then set the hyperparameters for the $i$th ROMES model to be the minimizer
of the associated loss over the validation set, i.e.,
\begin{equation}\label{eq:hyperparamsOpt}
\hyperparamsArg{i} =
\underset{\hyperparams\in\hyperparamsSet}{\arg\min}\loss{i}{\hyperparams}.
\end{equation}

One benefit of this cross-validation approach is that it admits flexibility in
selecting the loss function $\lossArg{\hyperparams}{i}{j}$, which determines
hyperparameter selection. Because one of our objectives is to achieve statistical
validation (condition \ref{prop:validated} in Section
\ref{sec:statisticalmodel}), we can define this loss function to align with
different notions of statistical validation; this is particularly important
when the errors do not exhibit a Gaussian distribution, as this case
precludes the ability to
achieve statistical validation in every possible
metric. We thus propose the following loss functions:
\begin{enumerate}
	\item the negative log-likelihood
		\begin{equation}\label{eq:loglikelihoodloss}
	 \lossloglikelihoodArg{\hyperparams}{i}{j} \defeq
			\frac{\card{\ROMEStrainingArgs{i}{j}}}{2}\ln(2\pi) +
			\frac{1}{2}\sum_{\indicator\in\ROMEStrainingArgs{i}{j}}\ln(\overallSTDArg{i,j}{\indicator;\hyperparams}^2)
			+
		\frac{1}{2}\sum_{(\indicator,\stateErrorRedSingle)\in\ROMEStrainingArgs{i}{j}}
			\frac{(\stateErrorRedSingle-\overallmeanArg{i,j}{\indicator;\hyperparams})^2}{\overallSTDArg{i,j}{\indicator;\hyperparams}^2}
			  \end{equation}
	\item the matching of a
		$\frequency$-prediction interval with $\frequency\in(0,1)$, i.e.,
		\begin{equation} \label{eq:lossFrequency}
\lossFrequencyArg{\hyperparams}{i}{j} =
		(\frequency-\validationFrequencyVal{i}{\frequency}{j}{\hyperparams})^2;
	\end{equation}
		where the validation frequency is
\begin{equation}
	\validationFrequencyVal{i}{\omega}{j}{\hyperparams}\defeq
	\frac{\card{\{(\indicator,\stateErrorRedSingle)\in\ROMEStrainingArgs{i}{j} |\
	\stateErrorRedSingle\in\predictionIntervalVal{i}{\omega}{\indicator}{j}{\hyperparams}\}}
	}
	{\card{\ROMEStrainingArgs{i}{j}}},
\end{equation}
		with prediction interval
\begin{align}
\begin{split}
	\predictionIntervalVal{i}{\omega}{\indicator}{j}{\hyperparams} \defeq [&
		\meanPredictionIntVal - \sqrt{2}\sigmaPredictionIntVal\erf^{-1}(\omega)
,
		\meanPredictionIntVal + \sqrt{2}\sigmaPredictionIntVal\erf^{-1}(\omega)
		].
\end{split}
\end{align}
\item the Komolgorov--Smirnov statistic, i.e.,
	$\lossKSArg{\hyperparams}{i}{j}$, which measures the maximum discrepancy between the
		cumulative distribution function (CDF) of the standard Gaussian distribution $\normal{0}{1}$ and
		the empirical CDF of the standardized data $\{\stateErrorRedSingle -
		\overallmeanArg{i,j}{\indicator;\hyperparams})/\overallSTDArg{i,j}{\indicator;\hyperparams}\}_{(\indicator,\stateErrorRedSingle)\in\ROMEStrainingArgs{i}{j}}$.
\end{enumerate}
We also consider a linear combination of these proposals to be employed as the
loss function $\lossNo$, e.g., a linear combination of the losses
$\lossFrequencyNo$ for different values of $\frequency$.

After the hyperparameters $\hyperparamsArg{i}$ have been computed according to
Eq.~\eqref{eq:hyperparamsOpt} the associated basis-expansion coefficients are
computed via maximum likelihood estimation on the full training set
$\ROMEStrainingArg{i}$ as
\begin{equation} \label{eq:regMLEfinal}
	\regweights_{i} =
	\left(\vondermondeROMESSingle^T(\kernelMatFunc{\indicatorVecArg{i}}{\indicatorVecArg{i}}+\varianceArg{i}\identity)^{-1}\vondermondeROMESSingle\right)^{-1}
	\vondermondeROMESSingle^T(\kernelMatFunc{\indicatorVecArg{i}}{\indicatorVecArg{i}}+\varianceArg{i}\identity)^{-1}\responseTrainVecROMESSingle,
\end{equation}
where $\indicatorVecArg{i}\in\RR{\card{\ROMEStrainingArg{i}}}$ and
$\stateErrorRedVecArg{i}\in\RR{\card{\ROMEStrainingArg{i}}}$ denote the
vectorized features and responses associated with training set
$\ROMEStrainingArg{i}$.

Given the hyperparameters $\hyperparamsArg{i}$ and basis-expansion coefficients
$\regweights_{i}$, the statistical models for the error generalized
coordinates at
arbitrary prediction parameter instances $\param\in\paramEval$
 are
\begin{equation}\label{eq:stateErrorRedModelArgDef}
	\stateErrorRedModelArg{i}(\param)=\mappingApproxArg{i}(\indicatorArg{i}(\param))\sim\normal{
	\overallmeanArg{i}{\indicatorArg{i}(\param)}}{\overallVarArg{i}{\indicatorArg{i}(\param)}},\quad
	i=1,\ldots,\nrbtot,
\end{equation}
where $\overallmeanArgNo{i}$ and $\overallVarArgNo{i}$ denote the mean and
variance associated with the GP model for $i$th error generalized coordinate,
defined as
\begin{align}
	& \overallmeanArg{i}{\indicatorArg{i}(\param)}\defeq
	\kernelMatFunc{\indicatorArg{i}(\param)}{\indicatorVecArg{i}}
	(\kernelMatFunc{\indicatorVecArg{i}}{\indicatorVecArg{i}}+\varianceArg{i}\identity)^{-1}\responseTrainVecROMESSingle
+
	[1\ \indicatorArg{i}(\param)]\regweights_{i}
	\\
	&\overallVarArg{i}{\indicatorArg{i}(\param)} \defeq
	\kernelMatFunc{\indicatorArg{i}(\param)}{\indicatorArg{i}(\param)} -
\kernelMatFunc{\indicatorArg{i}(\param)}{\indicatorVecArg{i}}
	(\kernelMatFunc{\indicatorVecArg{i}}{\indicatorVecArg{i}}+\varianceArg{i}\identity)^{-1}
	\kernelMatFunc{\indicatorVecArg{i}}{\indicatorArg{i}(\param)} +
	\varianceArg{i}.
\end{align}

\subsection{State and quantity-of-interest statistical
models}\label{seq:stateqoimodels}

Recall from Eqs.~\eqref{eq:stateErrorDecomp}, \eqref{eq:inplaneerrorDecomp}, and
\eqref{eq:stateErrorOutofplaneApprox} that we can approximate the state error
as
\begin{align} \label{eq:stateErrorDecompApprox}
	\stateError(\param) \approx
	\rbmat\stateErrorInplaneRed(\param) +
	\rbmatperp\stateErrorOutofplaneRed(\param),
\end{align}
where the approximation arises from the low-dimensional approximation of the
out-of-plane error from expression \eqref{eq:stateErrorOutofplaneApprox}.
Replacing the error generalized coordinates with their statistical
models such that
$\stateErrorInplaneRed\leftarrow\stateErrorInplaneRedModel$
and
$\stateErrorOutofplaneRed\leftarrow\stateErrorOutofplaneRedModel$
with
$
\stateErrorInplaneRedModelArg{i}\defeq
\stateErrorRedModelArg{i}
$, $i=1,\ldots,\nrb$ and
$
\stateErrorOutofplaneRedModelArg{i}\defeq
\stateErrorRedModelArg{i+\nrb}
$, $i=1,\ldots,\nrbperp$ in expression \eqref{eq:stateErrorDecompApprox} yields a statistical model for the state of the form
\begin{align} \label{eq:stateErrorModelDef}
	\stateErrorModel:\paramDummy\mapsto
	\rbmat\stateErrorInplaneRedModel(\paramDummy) +
	\rbmatperp\stateErrorOutofplaneRedModel(\paramDummy),
\end{align}
where $\stateErrorModel(\paramDummy):\samplespace\rightarrow\RR{\ndof}$ is an
$\ndof$-vector
of Gaussian random variables whose $i$th entry has a probability distribution
\begin{equation}
	\stateErrorModelEntry{i}(\paramDummy)\sim \normal
	{\sum_{j=1}^{\nrbtot}\rbmatTotEntry{i}{j}\overallmeanArg{j}{\indicatorArg{j}(\paramDummy)}}{
	\sum_{j=1}^{\nrbtot}\rbmatTotEntry{i}{j}\overallVarArg{j}{\indicatorArg{j}(\paramDummy)}}.
\end{equation}
Substituting $\stateError\leftarrow\stateErrorModel$ in the
definition of the state error \eqref{eq:stateError} yields a statistical model
for the state, which comprises deterministic and stochastic components, i.e.,
\begin{align} \label{eq:stateModelDef}
	\stateModel:\paramDummy \mapsto \underbrace{\stateROM(\paramDummy)}_\text{deterministic} +
	\underbrace{\stateErrorModel(\paramDummy)}_\text{stochastic}
\end{align}
such that $\stateModel(\paramDummy):\samplespace\rightarrow\RR{\ndof}$ is
also an
$\ndof$-vector
of Gaussian random variables whose $i$th entry is distributed as
\begin{equation}
	\stateModelEntry{i}(\paramDummy)\sim\normal{
	\stateROMEntry{i}(\paramDummy) +
\sum_{j=1}^{\nrbtot}\rbmatTotEntry{i}{j}\overallmeanArg{j}{\indicatorArg{j}(\paramDummy)}}{
\sum_{j=1}^{\nrbtot}\rbmatTotEntry{i}{j}\overallVarArg{j}{\indicatorArg{j}(\paramDummy)}}.
\end{equation}
Substituting this state model into the quantity-of-interest functional yields
the corresponding statistical model for the quantity of interest
\begin{align} \label{eq:qoiModelDef}
	\qoiModel:\paramDummy  \mapsto \qoiFunc(\stateModel(\paramDummy);\paramDummy)
\end{align}
and associated quantity-of-interest error model
\begin{align} \label{eq:qoiErrorModelDef}
	\qoiErrorModel:\paramDummy\mapsto \underbrace{\qoiModel(\paramDummy)}_\text{stochastic}  -
	\underbrace{\qoiROM(\paramDummy)}_\text{deterministic},
\end{align}
where
$\qoiModel(\paramDummy),\qoiErrorModel(\paramDummy):\samplespace\rightarrow\RR{\nqoi}$
are
$\nqoi$-vectors of random variables, which are Gaussian if the
quantity-of-interest functional
$\qoiFunc$ is linear in its first argument.

\subsection{Offline/online decomposition}\label{sec:offlineonline}
Algorithms \ref{alg:offline} and  \ref{alg:online} describe the steps required
for the offline and the online stages of the proposed method, respectively.
If hyper-reduction is applied to the ROM governing equations \eqref{eq:romEq}
or ROM dual system \eqref{eq:dualinplaneROM}, then hyper-reduction steps can be integrated in
the standard way, i.e., through additional data collection during the offline
stage.
We now highlight several attributes of the proposed method.

\begin{algorithm}[h!tbp]
\caption{Offline stage}
\begin{algorithmic}[1]\label{alg:offline}
\REQUIRE Training parameter sets $\paramDual,
\paramROMES\subset\paramDomain$
\ENSURE Reduced basis matrices $\rbmat$; $\testrbmat$;
	$\rbmatDualArg{i}$, $i=1,\ldots,\nrbtot$; $\testrbmatDualArg{i}$, $i=1,\ldots,\nrbtot$; and
ROMES models $\mappingApproxArg{i}$, $i=1,\ldots,\nrbtot$
\STATE \label{step:basis} Construct trial reduced basis matrix $\rbmat$ and
test reduced basis matrix $\testrbmat$. For example, the
trial reduced basis matrix $\rbmat$ can be constructed
via
POD using states computed by solving
the FOM equations \eqref{eq:discretizedPB} at training points $\param\in\paramTrain$ and
the test reduced basis matrix $\testrbmat$ can be subsequently defined by
employing Galerkin or LSPG projection.
	\STATE \label{step:rbmatperp}Construct the out-of-plane error basis matrix $\rbmatperp$, e.g., from
the discarded POD modes or by computing the projection error of the FOM
solutions for a parameter set
(see Remark \ref{rem:outofplanemat}).
\STATE \label{step:dualbasis} Construct the dual trial-basis matrices
	$\rbmatDualArg{i}$, $i=1,\ldots,\nrbtot$ via POD using dual vectors computed
from solving the ROM equations \eqref{eq:romEq} for $i=1,\ldots,\nrbtot$ and subsequently the dual
equations \eqref{eq:dualinplane} for $\param\in\paramDual$. Subsequently, set
	the dual test-basis matrices $\testrbmatDualArg{i}$, $i=1,\ldots,\nrbtot$ via
Galerkin or LSPG projection.
\STATE\label{step:ROMEStraining} Compute ROMES training data $\ROMEStrainingArg{i}$,
$	i=1,\ldots,\nrbtot$ by computing
$\indicatorArg{i}(\param)$, $\stateErrorRedArg{i}(\param)$,
$\param\in\paramROMES$.
Computing the state error generalized coordinates $\stateErrorRedArg{i}$ requires
solving FOM equations \eqref{eq:discretizedPB} and ROM
equations \eqref{eq:romEq} for $\param\in\paramROMES$, and subsequently
projecting the state error via Eq.~\eqref{eq:stateErrorRedDef}. Computing the
error indicators $\indicatorArg{i}$ requires additionally solving the dual ROMs
\eqref{eq:dualinplaneROM} and computing the indicators via
Eq.~\eqref{eq:indicatorDef}.
\STATE \label{step:ROMESconstruction} Construct ROMES models
$\mappingApproxArg{i}$, $i=1,\ldots,\nrbtot$ by Gaussian-process regression
(Section \ref{sec:GPROMES}).
\end{algorithmic}
\end{algorithm}
\begin{algorithm}[h!btp]
\caption{Online stage}
\begin{algorithmic}[1]\label{alg:online}
\REQUIRE Online parameter instance
$\paramOnline\in\paramEval\subset\paramDomain$;
reduced basis matrices $\rbmat$; $\testrbmat$;
	$\rbmatDualArg{i}$, $i=1,\ldots,\nrbtot$; $\testrbmatDualArg{i}$,
	$i=1,\ldots,\nrbtot$; and ROMES models $\mappingApproxArg{i}$,
	$i=1,\ldots,\nrbtot$
\ENSURE Statistical models
for the state
$\stateModel(\paramOnline)$, state error
$\stateErrorModel(\paramOnline)$, quantity of
interest $\qoiModel(\paramOnline)$, and quantity-of-interest error
$\qoiErrorModel(\paramOnline)$
\STATE \label{step:ROMsol}Compute ROM state $\stateROM(\paramOnline)$ by solving ROM equations
\eqref{eq:romEq} for $\param = \paramOnline$.
\STATE\label{step:dual} Compute approximate dual solutions
$\dualApproxArg{i}(\paramOnline)$, $i=1,\ldots,\nrbtot$ by solving Eqs.~\eqref{eq:dualinplaneROM} and
evaluating Eqs.~\eqref{eq:dualApproxArgDef}.
\STATE Compute dual-weighted-residual error indicators
$\indicatorArg{i}(\paramOnline)$, $i=1,\ldots,\nrbtot$
via Eqs.~\eqref{eq:indicatorDef}.
\STATE Evaluate ROMES models for the error generalized coordinates
$\stateErrorRedModelArg{i}(\paramOnline)=\mappingApproxArg{i}(\paramOnline)$, $i=1,\ldots,\nrbtot$
via Eqs.~\eqref{eq:stateErrorRedModelArgDef}.
\STATE \label{step:statmodels}Compute statistical models for the state
$\stateModel(\paramOnline)$ (Eq.~\eqref{eq:stateModelDef}), state error
$\stateErrorModel(\paramOnline)$ (Eq.~\eqref{eq:stateErrorModelDef}), quantity of
interest $\qoiModel(\paramOnline)$ (Eq.~\eqref{eq:qoiModelDef}), and quantity-of-interest error
$\qoiErrorModel(\paramOnline)$ (Eq.~\eqref{eq:qoiErrorModelDef})
\end{algorithmic}
\end{algorithm}
\begin{remark}[Offline stage: training cost]
The primary cost of Steps \ref{step:dualbasis}--\ref{step:ROMESconstruction}
	of Algorithm \ref{alg:offline}
	during the offline stage is incurred by the need to
solve the
FOM equations \eqref{eq:discretizedPB} for $\param\in\paramROMES$,
the ROM equations \eqref{eq:romEq} for $\param\in\paramDual\cup\paramROMES$,
and the dual ROMs \eqref{eq:dualinplaneROM} for $\param\in\paramDual$.
If POD is employed in step \ref{step:basis} of Algorithm \ref{alg:offline} to compute the trial reduced basis
matrix, then the
FOM equations \eqref{eq:discretizedPB} must be also solved for
$\param\in\paramTrain$. To reduce the training burden, these sets can overlap,
although this risks sacrificing generalizability of the statistical models.
For example, if POD is used to compute $\rbmat$ in Step \ref{step:basis} of Algorithm \ref{alg:offline} and
one employs $\paramROMES\subseteq\paramTrain$, then no additional FOM
solves are required; however, the ROMES training
data $\ROMEStrainingArg{i}$,
$	i=1,\ldots,\nrbtot$ will likely include only small-magnitude errors (i.e.,
$\|\stateErrorRedArg{i}(\param)\|$ small for $\param\in\paramROMES$), as the ROM is typically accurate over
the training set $\paramTrain$. Thus, the ROMES models may be not
generalize to parameter
instances corresponding to large-magnitude ROM errors. Similarly,
one could employ $\paramROMES\subseteq\paramDual$ to avoid additional ROM
solves; however, this would lead to a training set with extremely accurate
error indicators (i.e., $\indicatorArg{i}(\param)$ accurately represents of
the error $\stateErrorRedArg{i}(\param)$ for $\param\in\paramROMES$), as the
dual ROM is typically accurate over the associated training
set $\paramDual$. Thus, the ROMES models may not generalize to parameter
instances corresponding to large-magnitude dual-ROM errors.
\end{remark}

\begin{remark}[Offline stage: specifying quantities of interest not
	required]\label{rem:qoiNotNeeded}
	Unlike the original ROMES method \cite{kevin:romes}, the proposed approach
	does not require prescribing the quantities of interest in the offline
	stage. This is apparent from Algorithm \ref{alg:offline}: no steps require
	specification of these quantities, and Step \ref{step:ROMESconstruction} of
	Algorithm \ref{alg:offline} constructs ROMES models
	for the generalized coordinates only. Instead, the quantities of interest
	must be prescribed only in the online stage. This facilitates exploratory
	scenarios wherein the analyst may not know \textit{a priori} which quantities
	are of interest; further, it enables the statistical models of
	high-dimensional quantities of interest characterized by $\nqoi$ `large'
	(e.g., field quantities) to be
	efficiently computed, as
the quantity-of-interest models $\qoiModel$ and $\qoiErrorModel$ can be obtained
by substituting the low-dimensional state-error model $\stateErrorModel$ in the
quantity-of-interest functional $\qoiFunc$ in Step \ref{step:statmodels} of
Algorithm \ref{alg:online}.
\end{remark}

\begin{remark}[Online stage: comparison with a `ROM-only' approach]\label{rem:statmodel}
Note that typical `ROM-only' approaches execute only Step
\ref{step:ROMsol} of Algorithm \ref{alg:online} and directly employ
$\stateROM(\paramOnline)$ as an approximation to $\state(\paramOnline)$ and
$\qoiROM(\paramOnline)$ as an approximation to $\qoiFOM(\paramOnline)$;
some approaches additionally compute bounds for the approximation errors
$\|\state(\paramOnline)-\stateROM(\paramOnline)\|$
and
$\|\qoiFOM(\paramOnline)-\qoiROM(\paramOnline)\|$.

 In contrast, the proposed
approach computes
statistical models for the FOM state $\stateModel(\paramOnline)$ and quantity
of interest $\qoiModel(\paramOnline)$
	via Steps \ref{step:dual}--\ref{step:statmodels} of Algorithm \ref{alg:online}, which can be directly employed as
statistical models for those quantities that model their epistemic
	uncertainty. However, this benefit incurs an additional computational cost
	through the execution of
Steps \ref{step:dual}--\ref{step:statmodels}. Thus, even if the proposed
	approach is able to generate more accurate predictions of the state and
	quantity of interest for a fixed ROM dimension $\nrb$, it does so at an
	increased cost. It is unclear \textit{a priori} whether this approach is
	more computationally efficient than a `ROM-only' approach, as a larger ROM
	dimension $\nrb$ could be employed in the ROM-only approach to match the
	computational cost of the proposed method while also increasing its
	accuracy.

	However, we note that the additional cost of the proposed method
	is dominated by the dual ROM solves performed in Step \ref{step:dual}.
	These dual ROM equations \eqref{eq:dualinplaneROM} are always linear in their first argument, even when the ROM
	equations \eqref{eq:romEq} are nonlinear in their first argument. Thus, we expect the proposed method to
	perform favorably relative to a ROM-only approach when the governing
	equations are nonlinear.
\end{remark}

\section{Numerical Results} \label{sec:NumResults}
This section assesses the ability of the proposed method to construct
accurate error models on two model problems.
Experiments in Section \ref{sec:5_1} consider the Sobolev $H^1(\spatialDomain)$ inner product, whereas those in Section \ref{sec:5_2} consider the Euclidean
inner-product such that
$\metricmat=\identity$ in the definition of the projectors $\projectorParallel$
and $\projectorPerp$.
All timings are obtained by performing calculations on an Intel(R) Core i7-8700K CPU with 64 Gb DDR4 2666 MHz RAM using \textsc{Matlab}(R).
The \textsc{Matlab}(R) code used to generate these results has been released and is freely available at \url{https://stefanopagani.github.io/ChROME/}.


\subsection{Test case 1: linear diffusion}\label{sec:5_1}
\begin{figure}[h!]
\centering
\includegraphics[width=0.8\textwidth]{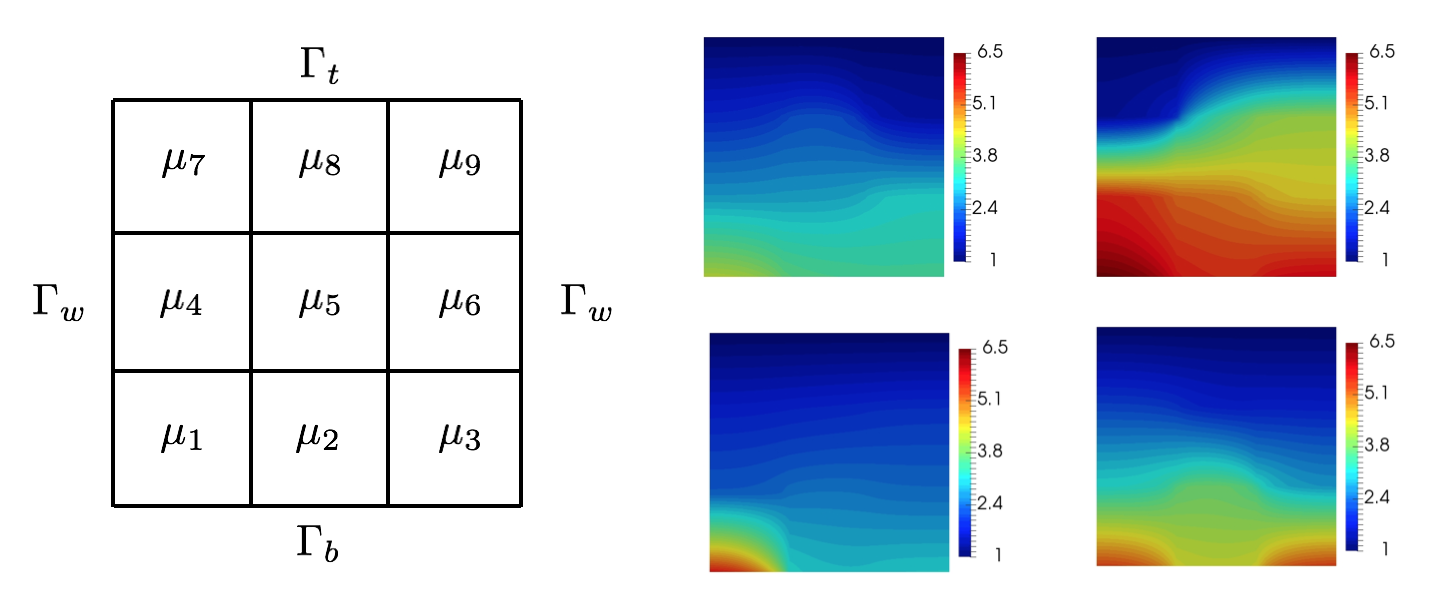}
 \caption{\textit{Test case 1}. Schematic representation of the
	 computational domain and finite-element
 solutions for different values of the system parameters $\param$.}
 \label{fig:domaintest1}
\end{figure}

We first assess the method on a problem characterized by a residual that
is linear in its first argument.
%
We consider  the following diffusion problem:
\begin{equation} \label{TestPDE}
    \begin{array}{c c}
      \begin{cases}
         \nabla \cdot( \diffusionConstant(\xx,\param) \nabla u) = 0, &
				 \xx\in\spatialDomain=[ 0, 1 ]^2\\
        \diffusionConstant(\xx,\param) \nabla u \cdot \mathbf{n} = 0, &  \xx\in \Gamma_w \\
    \diffusionConstant(\xx,\param) \nabla u \cdot \mathbf{n} = 1, &  \xx\in \Gamma_b \\
              u = 0, & \xx\in \Gamma_t,
    \end{cases}
  \end{array}
\end{equation}
where $\partial \spatialDomain=\Gamma_w \cup \Gamma_b \cup \Gamma_t$ (see
Fig.\
\ref{fig:domaintest1}, left). Here, $\diffusionConstant(\xx,\param)$ denotes
the parametrized diffusion coefficient, set to
\[
	\diffusionConstant(\xx,\param)=  0.01 \indicatorFunction{\spatialDomainArg{0}}(\xx) + \sum_{i=1}^9
	\mu_i \indicatorFunction{\spatialDomainArg{i}} (\xx) ,
\]
where $\indicatorFunction{A}(\xx)$ is the indicator function of the set $A$; we set
the parameter domain to
$\paramDomain = [0.01,1]^9$ with dimension $\nparam = 9$.
We discretize the spatial domain using the finite-element method on a
computational mesh given by  $5270$ triangular elements and quadratic finite
elements. This yields FOM governing equations of the form
\eqref{eq:discretizedPB} with $\ndof=2726$ degrees of freedom,
where $\res$ linear in its first argument such that Eq.~\eqref{eq:reslinear}
holds with $\mat$ symmetric and positive definite.

We execute the offline stage using Algorithm \ref{alg:offline} as follows. 
The training parameter sets---which comprise the algorithm inputs---are
constructed by drawing uniform random samples from the parameter domain
$\paramDomain$. We set
$\card{\paramDual}= 800$, while
$\card{\paramROMES}$ depends on the particular experiment.
In Step \ref{step:basis}, we apply POD to FOM solutions
computed at parameter instances
$\paramTrain$ with $\card{\paramTrain}=500$
drawn uniformly at random from
the parameter domain $\paramDomain$.
We  employ Galerkin
projection such that $\testrbmat = \rbmat$; the reduced-subspace dimension $\nrb$
depends on the particular experiment.
Step \ref{step:rbmatperp} constructs the basis
matrix
$\rbmatperp$ from the discarded POD modes; the out-of-plane subspace dimension $\nrbperp$ also
depends on the particular experiment.
In Step \ref{step:dualbasis}, we construct a
single shared trial dual basis matrix $\rbmatDual$ (i.e., $\rbmatDualArg{i} =
\rbmatDual$, $i=1,\ldots,\nrbtot$)
by combining snapshots from from all $\nrbtot$ dual solves executed at
parameter instances $\param\in\paramDual$ (see
Remark \ref{rem:uniqueShared}).
We also employ Galerkin projection for the dual problem such
that  $\testrbmatDualArg{i}=\testrbmatDual = \rbmatDual$,
$i=1,\ldots,\nrbtot$; the dual-basis dimension $\nrbDual$ also depends on the
particular experiment.
For constructing the ROMES models via Gaussian-process regression in Step \ref{step:ROMESconstruction},
we apply the procedure described in Section \ref{sec:GPROMES}.
We adopt the squared-exponential kernel function, which is
defined as \begin{equation}\label{eq:squaredExpKernel}
	\kernelFuncNo:(\generalSetOneVec,\generalSetTwoVec)\mapsto \signalsd
	\exp\left(\frac{\|\generalSetOneVec-\generalSetTwoVec\|_2^2}{2\width}\right)
\end{equation} and is characterized by hyperparameters $\hyperparamsKernel=(\signalsd,
\width)\in\RR{2}$.
We also consider several loss functions $\lossNo$ for hyperparameter selection as
described in Section \ref{sec:GPROMES}.

For the online stage, we execute Algorithm \ref{alg:online} for all parameter
instances in $\paramEval$, which comprises $\card{\paramEval} = 1500$
values drawn uniformly at random from $\paramDomain$. The remaining inputs to Algorithm
\ref{alg:online} result from the  outputs of Algorithm \ref{alg:offline}.

\subsubsection{ROMES model validation} \label{sec:GPvalTest1}

We first consider statistical validation of the ROMES models, i.e., Condition
\ref{prop:validated} in Section \ref{sec:statisticalmodel}.
We set the reduced subspace dimension to $\nrb=2$, the
out-of-plane subspace dimension to $\nrbperp = 0$, and the
(shared) dual-basis dimension to $\nrbDual=10$.
When constructing
the ROMES models in Step \ref{step:ROMESconstruction} according to the
description in Section \ref{sec:GPROMES}, we
define the set of candidate hyperparameter values $\hyperparamsSet$ via
uniform full-factorial sampling in each hyperparameter dimension characterized
by 12 equispaced values within the limits $\variance\in [0.01 \sigma_{t} ,
0.25 \sigma_{t}]$, $\signalsd\in [ 0.1 \sigma_{t}  , \sigma_{t} ]$, and
$\width\in [ 0.001 \sigma_{t} , 0.1 \sigma_{t} ]$, with $\sigma_{t}$ denoting
the standard deviation of the data set $\{
	\stateErrorRedArg{i}(\param))\}_{\param\in\paramROMES}$.

\begin{figure}[h!t]
 \centerline{
  \includegraphics{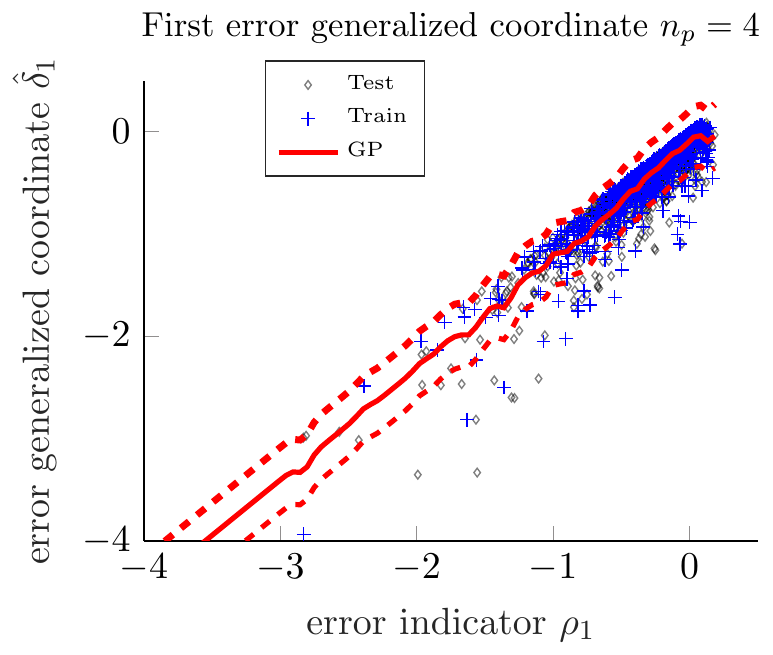}
  \includegraphics{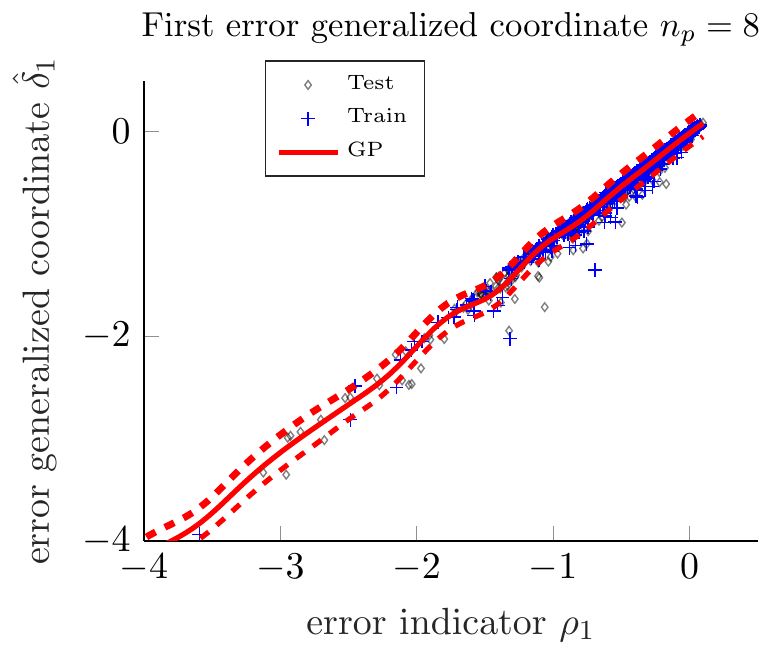}
  }
  \caption{\label{fig:GPMSE}  \textit{Test case 1}. ROMES models
	 constructed for the first two error generalized coordinates.
	 The solid line
	 represents the GP mean; the dashed lines represent the limits of the 99\%
	 prediction interval; the grey diamonds represent data related to online
	 points $\param\in\paramEval$, while the blue crosses represent training
	 data related to training
	 points $\param\in\paramROMES$.
	 We have employed $\nrb=2$, $\nrbperp = 0$, $\nrbDual=10$ and
	 we have selected hyperparameters according to Eq.~\eqref{eq:hyperparamsOpt} with
	 $\lossArg{\hyperparams}{i}{j} = \lossloglikelihoodArg{\hyperparams}{i}{j}$, and a
	 training set with $\card{\paramROMES} = 1000$.}
 \end{figure}

We first employ the negative log-likelihood loss function
$\lossArg{\hyperparams}{i}{j} = \lossloglikelihoodArg{\hyperparams}{i}{j}$
defined in Eq.~\eqref{eq:loglikelihoodloss} for hyperparameter selection.  Figure \ref{fig:GPMSE} reports
the resulting ROMES models constructed for the first error generalized
coordinate using a training set with $\card{\paramROMES} = 1000$ with two
values for the
dual-subspace dimension $\nrbDual$.
We note that for $\nrbDual=4$, the data appear to be somewhat skewed and the resulting
Gaussian process exhibits moderate variance.
By increasing the dual-subspace
dimension to $\nrbDual=8$, which incurs a larger computational cost due to
the increase dimension of the dual ROM equations \eqref{eq:dualinplaneROM}, the
feature becomes higher quality and thus leads to a lower-variance Gaussian
process.
Indeed, the ROMES model with $\nrbDual=8$ appears to qualitatively capture the relationship between the
error indicator and error generalized coordinate well; we now
investigate this further.


We assess the effect of the number of
training-parameter instances $\card{\paramROMES}$ on prediction accuracy, as
measured by (1) the fraction of variance unexplained (FVU)
\begin{equation} \label{eq:FVUArgDef}
	\FVUArg{i}\defeq
	\frac{\sum_{\param\in\paramEval}(\stateErrorRedArg{i}(\param)-\meanPredictionInt
		)^2}
		{\sum_{\param\in\paramEval}
		(
\stateErrorRedArg{i}(\param)-
\stateErrorRedAvgArg{i}
		)^2}
\end{equation}
		where
$\stateErrorRedAvgArg{i}$ denotes the mean value of
the error $\stateErrorRedArg{i}(\param)$ for $\param\in\paramEval$; (2)
the validation frequency
\begin{equation} \label{eq:validationFrequencyDef}
	\validationFrequency{i}{\omega}\defeq
	\frac{\card{\{\param\in\paramEval\ |\
	\stateErrorRedArg{i}(\param)\in\predictionInterval{i}{\omega}{\param}\}}
	}
	{\card{\paramEval}},
\end{equation}
where $\predictionInterval{i}{\omega}{\param}$ denotes the $\omega$-prediction
interval associated with ROMES model $\stateErrorRedModelArg{i}$, i.e.,
\begin{equation}
	\predictionInterval{i}{\omega}{\param} \defeq [
		\meanPredictionInt - \sqrt{2}\sigmaPredictionInt\erf^{-1}(\omega)
,
		\meanPredictionInt + \sqrt{2}\sigmaPredictionInt\erf^{-1}(\omega)
		],
\end{equation}
where $\overallmeanArgNo{i}$ and $\overallVarArgNo{i}$ denote the mean and
variance associated with the $i$th ROMES model (see
Eq.~\eqref{eq:stateErrorRedModelArgDef}); and (3) the Komolgorov--Smirnov (KS)
statistic, which quantifies the maximum discrepancy between the cumulative
distribution function (CDF) of the standard Gaussian distribution
$\normal{0}{1}$ and the empirical CDF of the standardized samples of the error
generalized coordinates
$\{\stateErrorRedArg{i}(\param) -
\overallmeanArg{i}{\indicatorArg{i}(\param)})/\overallSTDArg{i}{\indicatorArg{i}(\param)}\}_{\param\in\paramEval}$.

\begin{figure}[h!t]
 \centerline{
 \includegraphics{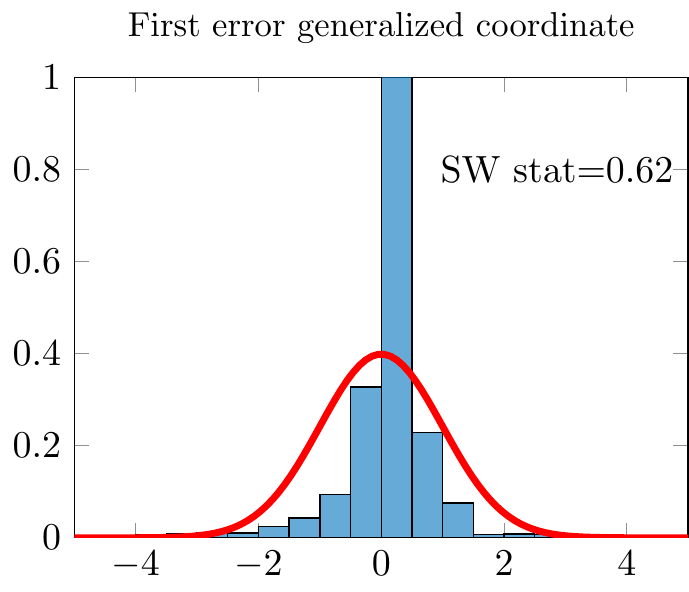}
 \includegraphics{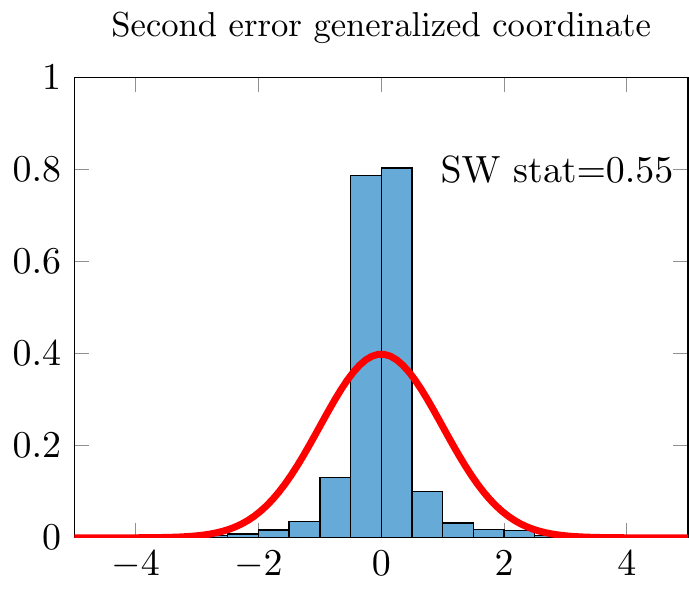}
 }
 \caption{\label{fig:histGPMSE}  \textit{Test case 1}. Histogram of the
	standardized data
	$\{\stateErrorRedArg{i}(\param) -
	\overallmeanArg{i}{\indicatorArg{i}(\param)})/\overallSTDArg{i}{\indicatorArg{i}(\param)}\}_{\param\in\paramEval}$,
	$i=1,2$ (blue bar plot)
	as compared to the PDF of the standard Gaussian distribution $\normal{0}{1}$
	(red curve).
	We have employed $\nrb=2$, $\nrbperp = 0$, $\nrbDual=10$, and
	 have selected hyperparameters according to Eq.~\eqref{eq:hyperparamsOpt} with
	 $\lossArg{\hyperparams}{i}{j} = \lossloglikelihoodArg{\hyperparams}{i}{j}$.
	 The number of training-parameter instances is $\card{\paramROMES}=1000$.
	} \vspace{0.5cm}
\end{figure}

While the FVU quantifies the ability of the ROMES model
$\stateErrorRedModelArg{i}$
to accurately model the error generalized coordinate $\stateErrorRedArg{i}$ in expectation, the
validation frequency and KS statistic
assess the statistical properties of the model, i.e., its ability to
accurately reflect the underlying data distribution.
Table  \ref{Fig:ValidationFrequencyInplane} reports these results, which show
that employing
$\card{\paramROMES}=1000$ is sufficient for the FVU to have reasonably
stabilized; thus, subsequent experiments in this section set
$\card{\paramROMES}=1000$. 
However, the converged prediction levels are not all correct; for example,
$\validationFrequency{1}{0.8} = 0.93$ even though
this value should be 0.8.
We observe that one likely source of this lack of statistical validation
arises from the fact that the true error does not exhibit Gaussian behavior as
reported in Figure \ref{fig:histGPMSE}; indeed, these data do not pass the
Shapiro--Wilk (SW) normality test, as they yield a SW statistic of $0.62$ for the first error generalized coordinate and of $0.55$ for the second. This
implies that it will not be possible to achieve statistical validation in
every possible metric if we employ Gaussian-process regression; rather, we may
only be able to satisfy a subset of statistical-validation criteria. This
motivates the need for tailored loss functions for hyperparameter selection as
described in Section \ref{sec:GPROMES},  as such loss functions enable the
method to target specific statistical-validation criteria.


 \begin{table}[h!tb]
 \centering
 \begin{tabular}{@{} c | c c c c | c c c c @{}}
   \toprule
		error index $i$& \multicolumn{4}{c |}{$1$} & \multicolumn{4}{c}{$2$} \\
  $\card{\paramROMES}$   & $100$ & $400$ & $700$ & $1000$ & $100$ & $400$ & $700$ & $1000$ \\
  \midrule
		FVU &  0.0128  &  0.0133  &  0.0131  &  0.0124  &  0.0113  &  0.0115  &  0.0112  &  0.0104 \\
	         $\validationFrequency{i}{0.8}$  & 0.8941  &  0.8961  &  0.9127  &  0.9300 & 0.8075  &  0.8901  &  0.9167  &  0.9307 \\
					 $\validationFrequency{i}{0.9}$  & 0.9267  &  0.9314  &  0.9400  &  0.9534 & 0.8601  &  0.9234  &  0.9394  &  0.9500\\
					 $\validationFrequency{i}{0.95}$ & 0.9414  &  0.9494  &  0.9587  &  0.9634 & 0.8894  &  0.9407  &  0.9527  &  0.9614\\
					 $\validationFrequency{i}{0.99}$ & 0.9560  &  0.9680  &  0.9727  &  0.9753 & 0.9234  &  0.9614  &  0.9720  &  0.9747 \\
					 KS statistic                    &  0.2169  &  0.2185  &  0.2324  &  0.2504  &  0.0985  &  0.1942  &  0.2189  &  0.2348  \\
   \bottomrule
 \end{tabular}
   \caption[ValidationFrequencyInplane]{  \label{Fig:ValidationFrequencyInplane}
		 \textit{Test case 1}. Convergence of error measures associated with the ROMES models
		 constructed for the first two error generalized coordinates as the number
		 of training-parameter instances $\card{\paramROMES}$ increases.
		 We have employed $\nrb=2$,
	  $\nrbperp = 0$, $\nrbDual=10$ and $\card{\paramEval}=1500$, and
	 have selected hyperparameters according to Eq.~\eqref{eq:hyperparamsOpt} with
	 $\lossArg{\hyperparams}{i}{j} = \lossloglikelihoodArg{\hyperparams}{i}{j}$.
	 }
 \end{table}


To this end, we now adopt several different strategies for defining the loss
function $\lossArg{\hyperparams}{i}{j}$ employed for hyperparameter selection
(see Section \ref{sec:GPROMES}).  In particular, we select the hyperparameters
$\hyperparamsArg{i}\equiv(\varianceArg{i},\signalsdArg{i}, \widthArg{i})$
characterizing the $i$th ROMES model by employing five different loss
functions
$\lossArg{\hyperparams}{i}{j}$: (1) the negative log-likelihood loss
$\lossloglikelihoodNo$ (Eq.~\eqref{eq:loglikelihoodloss}), (2) the loss based on
matching the 0.80-prediction interval $\lossEightyNo$
(Eq.~\eqref{eq:lossFrequency} with $\frequency=0.80$),
(3) the loss based on
matching the 0.95-prediction interval $\lossNinetyFiveNo$
(Eq.~\eqref{eq:lossFrequency} with $\frequency=0.95$), (4) the loss based on a
linear combination of $\frequency$-prediction interval losses
\begin{equation} \label{eq:lossPredictionNo}
	\lossPredictionNo \defeq
	\sum_{\frequency\in\{0.80,0.90,0.95,0.99\}}\lossFrequencyNo,
\end{equation}
 and (5) the loss based on the KS statistic $\lossKSNo$.
Table \ref{Fig:ValidationFrequencyInplaneTwo} reports these results for
$\card{\paramROMES}=1000$. 
We note that the loss function $\lossNo$  employed for hyperparameter
selection has a significant effect on the performance of the resulting ROMES
models according to  different statistical-validation criteria. In
particular, the loss function can be selected to optimize improve performance
with respect to particular criteria.
For example,
$\validationFrequency{1}{0.8} = 0.7941$
when $\lossEightyNo$ is adopted, while $\validationFrequency{1}{0.8} = 0.9300$
when  $\lossloglikelihoodNo$ is adopted. This is an important practical result
implying that the user should select the loss function to coincide with
desired statistical validation criterion.

Due to its favorable performance over a range of statistical-valation
criteria, we employ a loss function of $\lossNo=\lossPredictionNo$ in the remaining
experiments within Section \ref{sec:5_1}.
 \begin{table}[h!tb]
   \vspace{0.5cm}
 \centering
 \begin{tabular}{@{} c | c  c c c c }
   \toprule
		error generalized coordinate index $i$& \multicolumn{5}{c }{$1$} \\
	 loss function $\lossArg{\hyperparams}{i}{j}$ & $\lossloglikelihoodNo$ & $\lossPredictionNo$ &
	 $\lossEightyNo$ & $\lossNinetyFiveNo$ & $\lossKSNo$ \\
  \midrule
		FVU & 0.0124  &  0.0126  &  0.0123  &  0.0124  &  0.0136 \\
	 $\validationFrequency{i}{0.8}$   &  0.9300  &  0.8361  &  0.7941  &  0.9320  &  0.8368  \\
					 $\validationFrequency{i}{0.9}$   & 0.9534  &  0.8854  &  0.8381  &  0.9547  &  0.8881 \\
					 $\validationFrequency{i}{0.95}$  & 0.9634  &  0.9194  &  0.8761  &  0.9640  &  0.9194  \\
					 $\validationFrequency{i}{0.99}$  & 0.9753  &  0.9507  &  0.9294  &  0.9747  &  0.9487 \\
					 Komolgorov--Smirnov statistic  &  0.2504  &  0.2294  &  0.2350  &  0.2489  &  0.2288  \\
   \bottomrule
		error generalized coordinate index $i$& \multicolumn{5}{c }{$2$} \\
	 loss function $\lossArg{\hyperparams}{i}{j}$ & $\lossloglikelihoodNo$ & $\lossPredictionNo$ &
	 $\lossEightyNo$ & $\lossNinetyFiveNo$ & $\lossKSNo$\\
  \midrule
		FVU  &  0.0104  &  0.0101  &  0.0114  &  0.0111  &  0.0098 \\
	 $\validationFrequency{i}{0.8}$   & 0.9307  &  0.8661  &  0.8008  &  0.9427 &   0.7209   \\
					 $\validationFrequency{i}{0.9}$   &  0.9500  &  0.8967  &  0.8534  &   0.9594  &  0.7648 \\
					 $\validationFrequency{i}{0.95}$  &  0.9614  &  0.9154  &  0.8754 &  0.9700  &  0.7995  \\
					 $\validationFrequency{i}{0.99}$  & 0.9747  &  0.9407  &  0.9107  &  0.9767  &  0.8454 \\
					 Komolgorov--Smirnov statistic  & 0.2348  &  0.1731  &  0.1196  &  0.2554  &  0.0821 \\
   \bottomrule
 \end{tabular}
   \caption[ValidationFrequencyInplane]{  \label{Fig:ValidationFrequencyInplaneTwo}
 	 \textit{Test case 1}. Statistical-validation criteria evaluated on $\paramEval$ (with $\card{\paramEval}=1500$) for ROMES models when
	 different loss functions $\lossNo$ are employed for hyperparameter
	 selection according to Eqs.~\eqref{eq:lossi} and \eqref{eq:hyperparamsOpt}
	 in Section \ref{sec:GPROMES}.
We have employed $\nrb=2$, $\nrbperp = 0$,
	 $\nrbDual=10$ and have selected hyperparameters
	 according to Eq.~\eqref{eq:hyperparamsOpt} with
	 the specified loss function $\lossArg{\hyperparams}{i}{j}$, $10$ K-fold subdivisions and a search grid made by $12^3$ values.
	 The number of training-parameter instances is $\card{\paramROMES}=1000$.
	}
 \end{table}

\subsubsection{In-plane and out-of-plane error
approximation}\label{sec:inplaneoutofplane}
We now assess the ability of the ROMES method to accurately
approximate the in-plane error $\stateErrorInplane$ and the
out-of-plane error $\stateErrorOutofplane$.
To assess the ability of the method to approximate the
former, we compare the
mean relative ROM error
\begin{equation} \label{eq:ROMerror}
\errorROM \defeq \expectationParam{ \frac{ \| \stateError(\param) \|_2 }{ \| \state(\param) \|_2 } },
\end{equation}
with the mean relative ROM error after applying the in-plane ROMES correction
\begin{equation} \label{eq:ROMESerror}
\errorROMESParallel \defeq \expectationParam{ \frac{ \| \stateError(\param) -
	\rbmat  \expectation{\stateErrorInplaneRedModel(\param)}\|_2 }{ \|
	\state(\param) \|_2 } },
\end{equation}
and
the mean relative projection error
\begin{equation} \label{eq:Perror}
\errorROMParallel\defeq\expectationParam{ \frac{ \| \stateError(\param)-\stateErrorInplane(\param)\|_2 }{ \| \state(\param) \|_2 }  }.
\end{equation}
We note that the latter represents the minimum value achievable by the ROM
error with the in-plane ROMES correction.

\begin{figure}[h!t]
  \vspace{-0.1cm}
 \centerline{
  \includegraphics{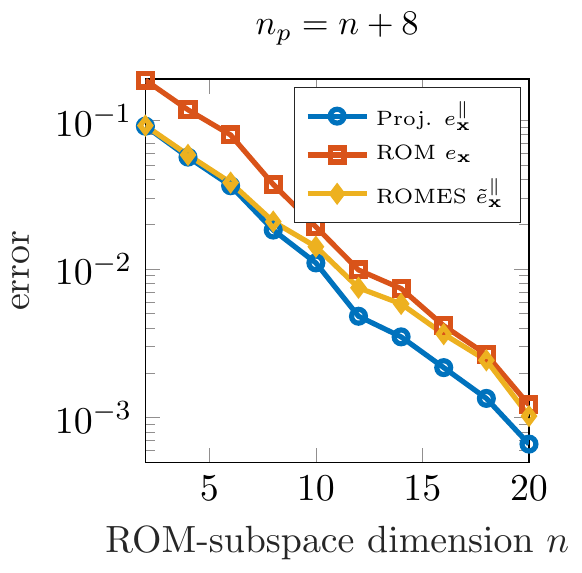}
  \includegraphics{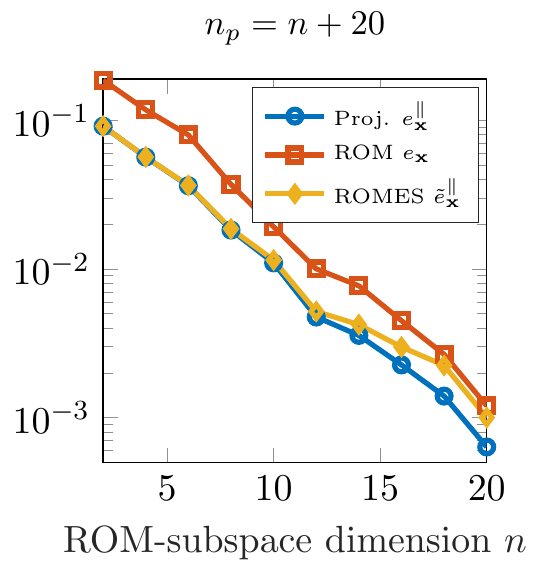}
  \includegraphics{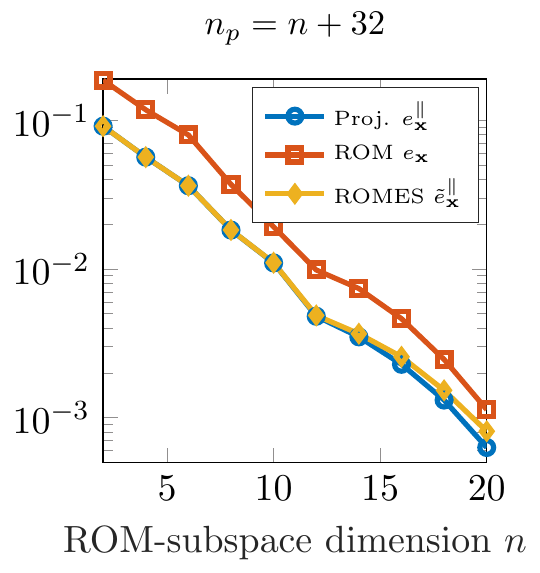}
 }
 \caption{\label{fig:inplane_plot}
 \textit{Test case 1}.
 Mean relative ROM error $\errorROM$ (red),  mean relative ROM error
	after applying the in-plane ROMES correction $\errorROMESParallel$ (yellow),
	and
	mean relative projection error $\errorROMParallel$ (blue) for a varying
	reduced-subspace dimension $\nrb$ and dual-subspace dimension $\nrbDual$.
	Here, we
set $\nrbperp = 0$, $\lossArg{\hyperparams}{i}{j} =
\lossPredictionArg{\hyperparams}{i}{j}$,
$\card{\paramROMES} = 1000$ and $\card{\paramEval}=1500$. \vspace{-0.1cm}
	}
\end{figure}

We set $\nrbperp = 0$, $\lossArg{\hyperparams}{i}{j} =
\lossPredictionArg{\hyperparams}{i}{j}$, and
$\card{\paramROMES} = 1000$ training-parameter instances, and we
vary the reduced-subspace dimension $\nrb$ and dual-basis dimension
$\nrbDual$.  Figure \ref{fig:inplane_plot} reports the results obtained for
error measures $\errorROM$, $\errorROMESParallel$, and $\errorROMParallel$ over a
range of values for $\nrb$ and $\nrbDual = \nrb + i$ for $i=8,14,20$.  These
figures illustrate that the ROMES models for the in-plane ROMES correction $\stateErrorInplaneRedModel$ enable the mean state error to be
significantly reduced with respect to the ROM error, as $\errorROMESParallel$ is smaller than
$\errorROM$ in all cases; moreover, the mean relative ROM error after applying
the in-plane ROMES correction state error $\errorROMESParallel$
approaches the optimal value defined by the mean relative projection error $\errorROMParallel$ as the
dual-basis dimension $\nrbDual$ increases. Indeed, because equality in
\eqref{eq:stateErrorInplaneRedExpressionTwoCheap} holds when the ROM
approximation of the duals is exact (i.e.,
$\dualArg{i}\approx\dualApproxArg{i}$) and the residual is linear in its first
argument (i.e., Eq.~\eqref{eq:reslinear} holds)---the latter of which is true
for this problem---we expect the ROMES models to be extremely accurate for
this problem as the dimension of the dual reduced basis $\nrbDual$ becomes
large. Thus, we conclude that the proposed approach is indeed able to
accurately approximate the in-plane error. That is, the proposed method
constructs an accurate statistical closure model.



%

\begin{figure}[h!t]
  \vspace{-0.1cm}
 \centerline{
 \includegraphics[width=0.45\textwidth]{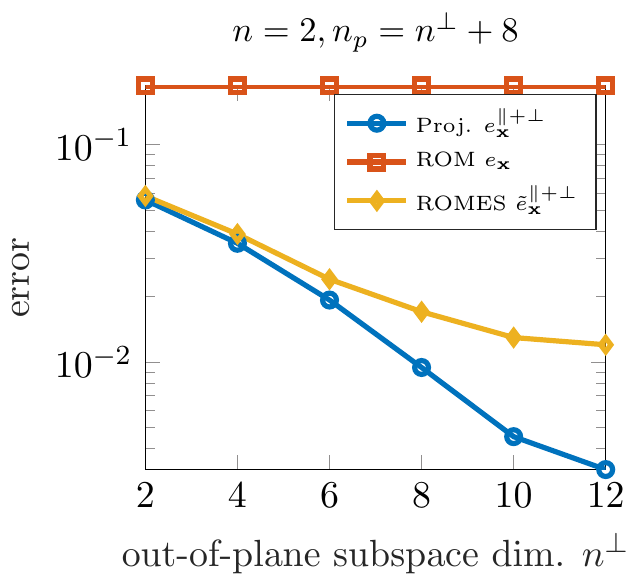}
 \includegraphics[width=0.45\textwidth]{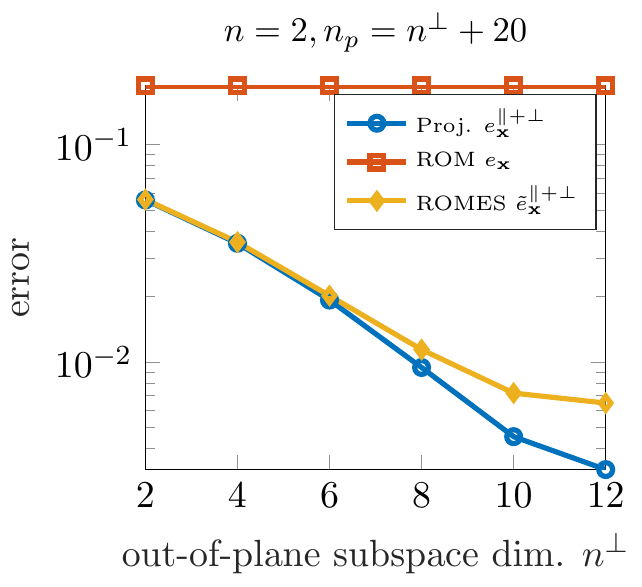}
 }
 \centerline{
  \includegraphics[width=0.45\textwidth]{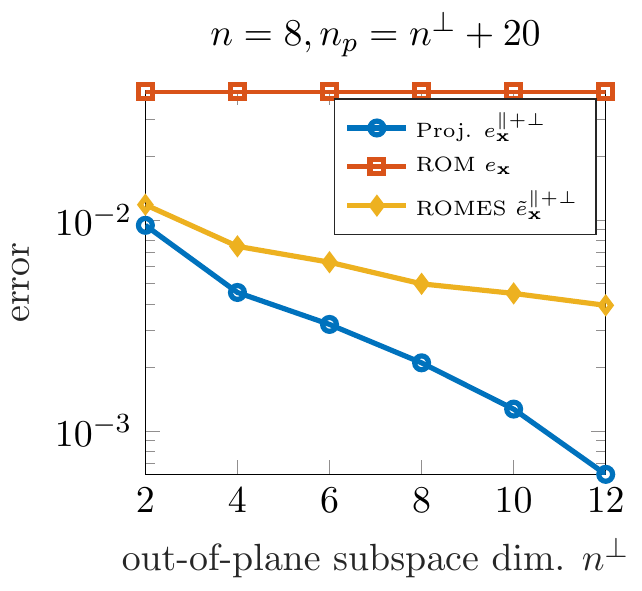}
  \includegraphics[width=0.45\textwidth]{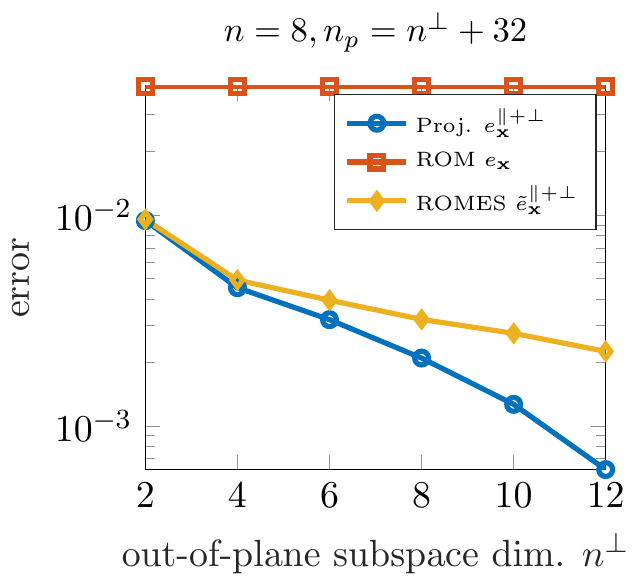}
}
 \caption{\label{fig:outofplane_plot}
\textit{Test case 1}.
 Mean relative ROM error $\errorROM$ (red),  mean relative ROM error
	after applying the in-plane and out-of-plane ROMES corrections
	$\errorROMESParallelPerp$ (yellow),
	and
	mean relative projection error $\errorROMParallelPerp$ (blue) for a varying
	reduced-subspace dimension $\nrb$, dual-subspace dimension $\nrbDual$, and
	out-of-plane subspace dimension $\nrbperp$.
	Here, we
set $\lossArg{\hyperparams}{i}{j} =
\lossPredictionArg{\hyperparams}{i}{j}$,
	$\card{\paramROMES} = 1000$, and
	$\card{\paramEval} = 1500$. \vspace{-0.1cm}}
\end{figure}

We perform a similar analysis for the out-of-plane
error by
comparing the
mean relative ROM error $\errorROM$
with the mean relative ROM error after applying both the in-plane and
out-of-plane ROMES corrections
\begin{equation} \label{eq:ROMESerrorout}
\errorROMESParallelPerp\defeq\expectationParam{ \frac{ \| \stateError(\param) -
	\rbmat\expectation{\stateErrorInplaneRedModel(\param)} - \rbmatperp
	\expectation{\stateErrorOutofplaneRedModel(\param)} \|_2 }{ \| \state(\param) \|_2 } }
\end{equation}
and
the mean relative projection error
\begin{equation} \label{eq:Perrorout}
\errorROMParallelPerp\defeq\expectationParam{ \frac{ \|
	\stateError(\param)-\stateErrorInplane(\param)
	-\stateErrorOutofplane(\param)\|_2 }{ \| \state(\param) \|_2 } },
\end{equation}
which represents the minimum value achievable by the ROM error with both
in-plane and out-of-plane ROMES correction.  Figure \ref{fig:outofplane_plot}
reports the results obtained for error measures $\errorROM$,
$\errorROMESParallelPerp$,  and $\errorROMParallelPerp$ for various values of
the reduced-subspace dimension $\nrb$, the dual-basis dimension $\nrbDual$,
and the out-of-plane subspace dimension $\nrbperp$.  We again observe that the
ROMES corrections enable significant error reduction, and performance
improves as the dual-basis dimension $\nrbDual$ increases.  Thus, we conclude
that the proposed approach is able to accurately approximate both the in-plane
and out-of-plane errors.

\subsubsection{Quantity-of-interest error
approximation}\label{sec:qoiApproxTest1}

We now consider the ability of the proposed ROMES models to construct
statistical models of quantities of interest $\qoiModel$ as proposed in
Section \ref{seq:stateqoimodels}.  To this end, we consider
$\nqoi=2$ quantities of interest
\begin{equation}
	\qoiFOM:\paramDummy\mapsto	\qoiFunc(\state(\paramDummy); \paramDummy ) =
	\begin{bmatrix}
		\qoiMatrixExperiment^{T} \state(\paramDummy)\\
\state(\paramDummy)^{T} \mathbf{M} \state(\paramDummy)
	\end{bmatrix},
\end{equation}
where $\qoiFOMArg{1}$ and $\qoiFOMArg{2}$ represent the mean value (i.e., $\int_{\spatialDomainArg{5}}u \mathrm{d}\xx= \qoiMatrixExperiment^{T} \state
$) and the mean
squared value (i.e., $\int_{\spatialDomainArg{5}}u^2 \mathrm{d}\xx=\state^{T} \mathbf{M} \state
$) of the state variable over
$\spatialDomainArg{5}$, respectively. We emphasize that these quantities of
interest were not specified during the offline stage (see Remark
\ref{rem:qoiNotNeeded}).

\begin{figure}[ht]
	 \subfigure[$\nrb=2$, $\nrbperp=0$]{
 \includegraphics[width=0.3\textwidth]{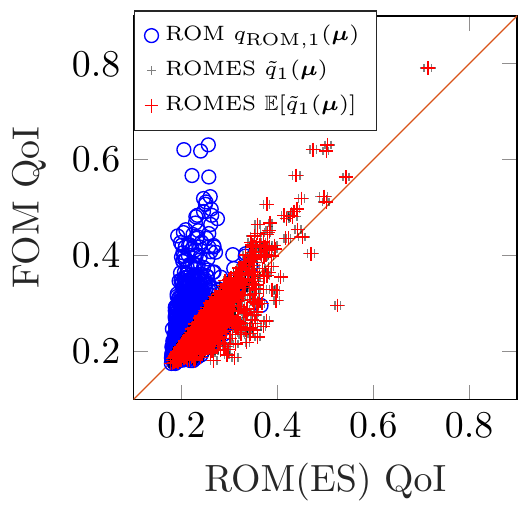}
	}
	 \subfigure[$\nrb=6$, $\nrbperp=0$]{
 \includegraphics[width=0.3\textwidth]{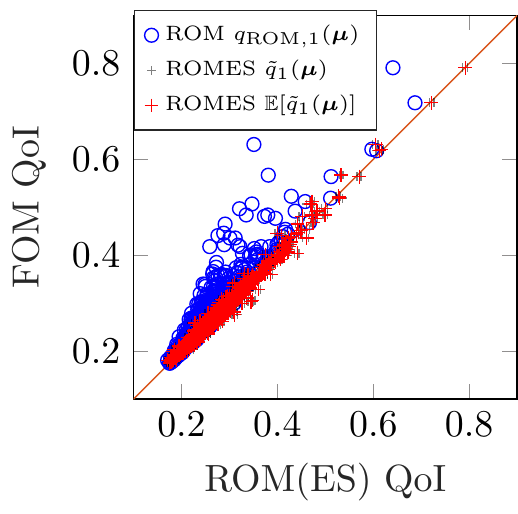}
	}
	 \subfigure[$\nrb=10$, $\nrbperp=0$]{
		 \includegraphics[width=0.3\textwidth]{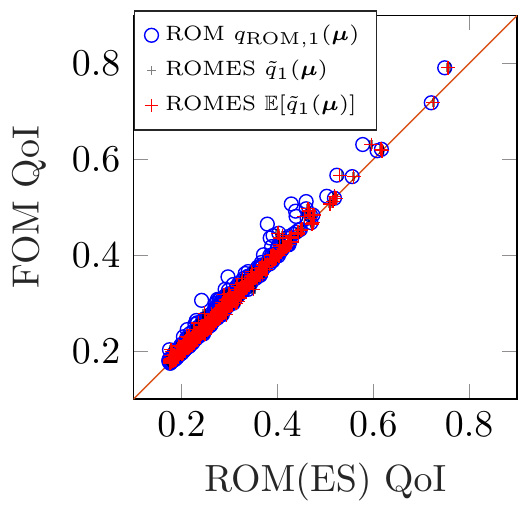}
	}
	 \subfigure[$\nrb=2$, $\nrbperp=2$]{
 \includegraphics[width=0.3\textwidth]{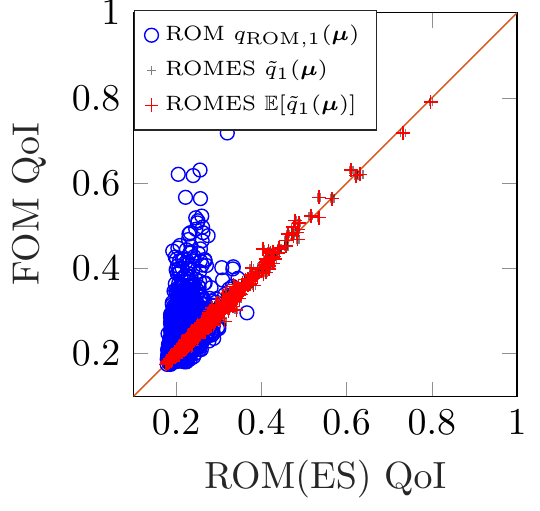}
	}
	 \subfigure[$\nrb=2$, $\nrbperp=6$]{
 \includegraphics[width=0.3\textwidth]{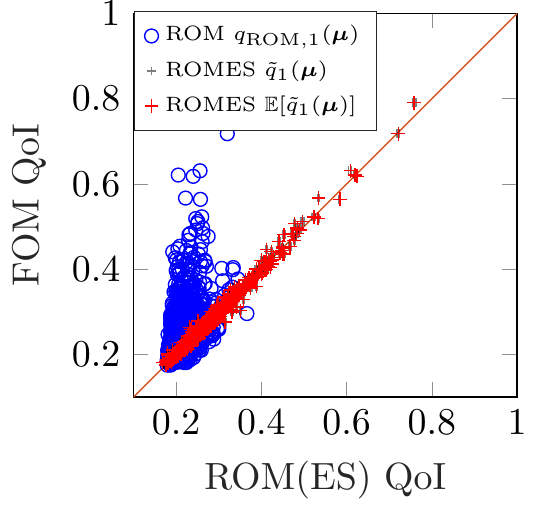}
	}
	 \subfigure[$\nrb=2$, $\nrbperp=10$]{
 \includegraphics[width=0.3\textwidth]{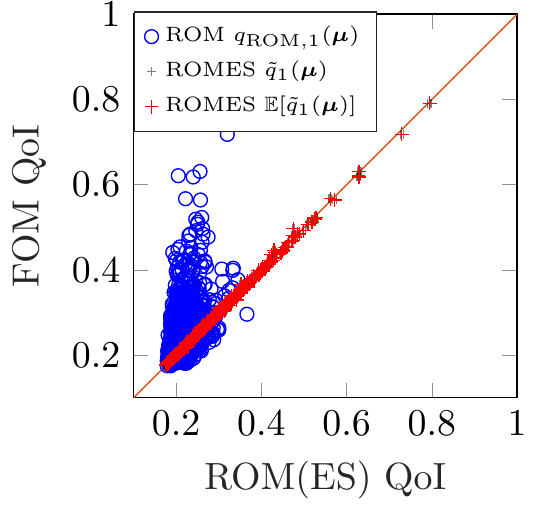}
	}
 \caption{\label{fig:Lu_output_scatter_plot}
 \textit{Test case 1}.
 Scatter plot of the FOM-computed (linear) quantity of interest $\qoiFOMArg{1}(\param)$ versus both the
	ROM-computed quantity of interest $\qoiROMArg{1}(\param)$ (blue circles) and the expected value of
	the ROMES-corrected quantity of interest $\expectation{\qoiModelArg{1}(\param)}$
	(red crosses) for
several values of the reduced-subspace dimension $\nrb$ and for
$\param\in\paramEval$. They grey crosses are computed by taking the maximum and the minimum value of the $\tilde{q}_1(\param)$ over $100$ realizations of the ROMES statistical model.
Here, we
set $\lossArg{\hyperparams}{i}{j} =
\lossPredictionArg{\hyperparams}{i}{j}$,
$\card{\paramROMES} = 1000$,
	and
	$\nrbDual=10$.
	}
\end{figure}

\begin{figure}[ht]
	 \subfigure[$\nrb=2$, $\nrbperp=0$]{
 \includegraphics[width=0.3\textwidth]{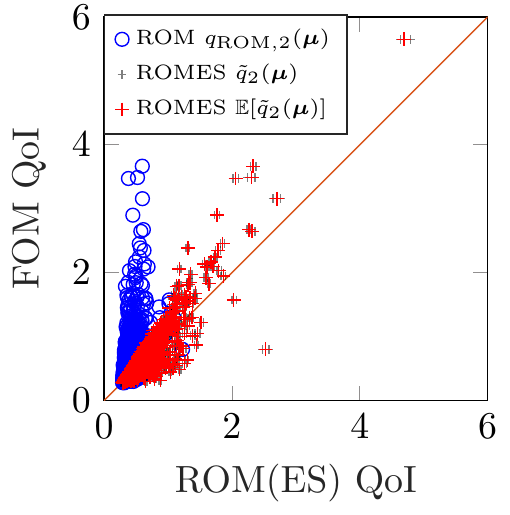}
	}
	 \subfigure[$\nrb=6$, $\nrbperp=0$]{
 \includegraphics[width=0.3\textwidth]{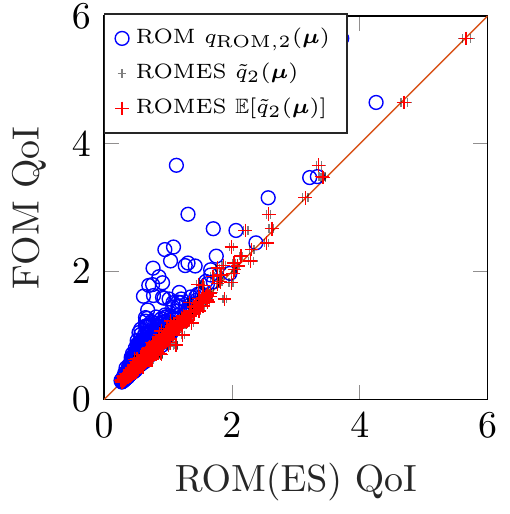}
	}
	 \subfigure[$\nrb=10$, $\nrbperp=0$]{
		 \includegraphics[width=0.3\textwidth]{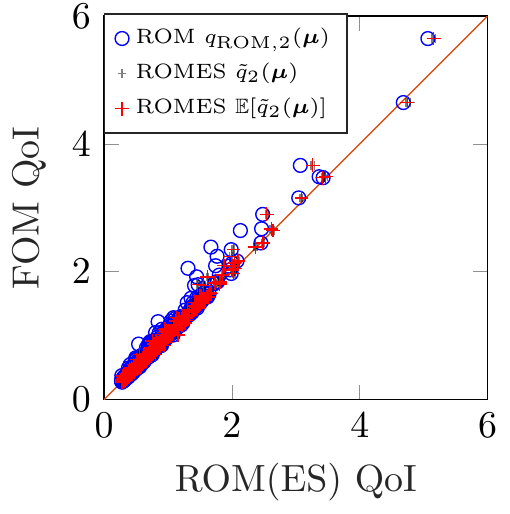}
	}
	 \subfigure[$\nrb=2$, $\nrbperp=2$]{
 \includegraphics[width=0.3\textwidth]{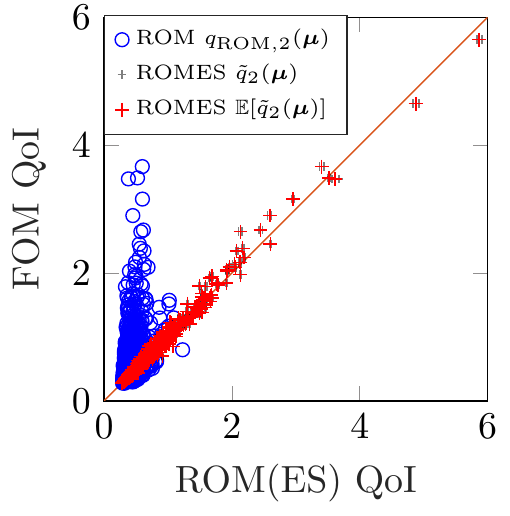}
	}
	 \subfigure[$\nrb=2$, $\nrbperp=6$]{
 \includegraphics[width=0.3\textwidth]{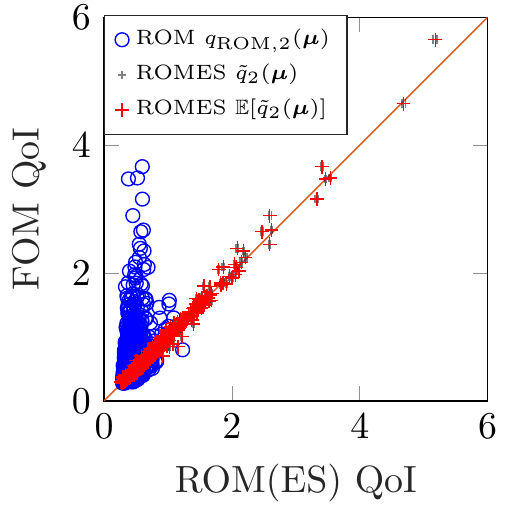}
	}
	 \subfigure[$\nrb=2$, $\nrbperp=10$]{
 \includegraphics[width=0.3\textwidth]{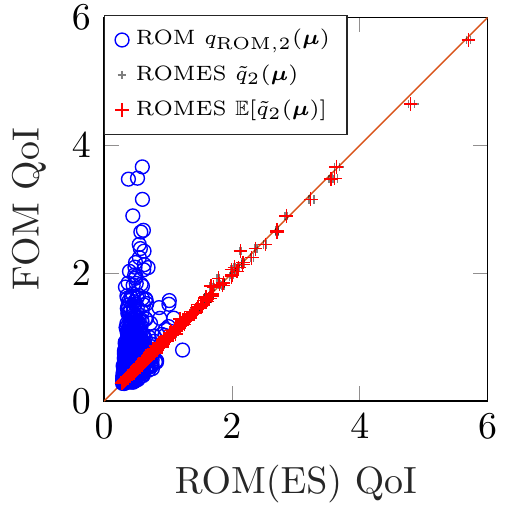}
	}
 \caption{\label{fig:Squaredu_output_scatter_plot}
 \textit{Test case 1}.
 Scatter plot of the FOM-computed (nonlinear) quantity of interest $\qoiFOMArg{2}(\param)$ versus both the
	ROM-computed quantity of interest $\qoiROMArg{2}(\param)$ (blue circles) and the expected value of
	the ROMES-corrected quantity of interest $\expectation{\qoiModelArg{2}(\param)}$
	(red crosses) for
several values of the reduced-subspace dimension $\nrb$ and for
$\param\in\paramEval$. They grey crosses are computed by taking the maximum and the minimum value of the $\tilde{q}_2(\param)$ over $100$ realizations of the ROMES statistical model.
Here, we
set $\lossArg{\hyperparams}{i}{j} =
\lossPredictionArg{\hyperparams}{i}{j}$,
$\card{\paramROMES} = 1000$,
	and
	$\nrbDual=10$.
	}
\end{figure}

We set the loss function to $\lossArg{\hyperparams}{i}{j} =
\lossPredictionArg{\hyperparams}{i}{j}$ and
number of training-parameter instances to $\card{\paramROMES} = 1000$.  Figures
\ref{fig:Lu_output_scatter_plot} and \ref{fig:Squaredu_output_scatter_plot}
plot the FOM-computed quantity of interest $\qoiFOMArg{i}(\param)$ versus both the
ROM-computed quantity of interest $\qoiROMArg{i}(\param)$ and the expected value of
the ROMES-corrected quantity of interest
$\expectation{\qoiModelArg{i}(\param)}$, $i=1,2$ for
several values of the reduced-subspace dimensions $\nrb$ and $\nrbperp$ and for
$\param\in\paramEval$.

\begin{figure}[h!t]
 \centerline{
  \includegraphics[height=0.35\textwidth]{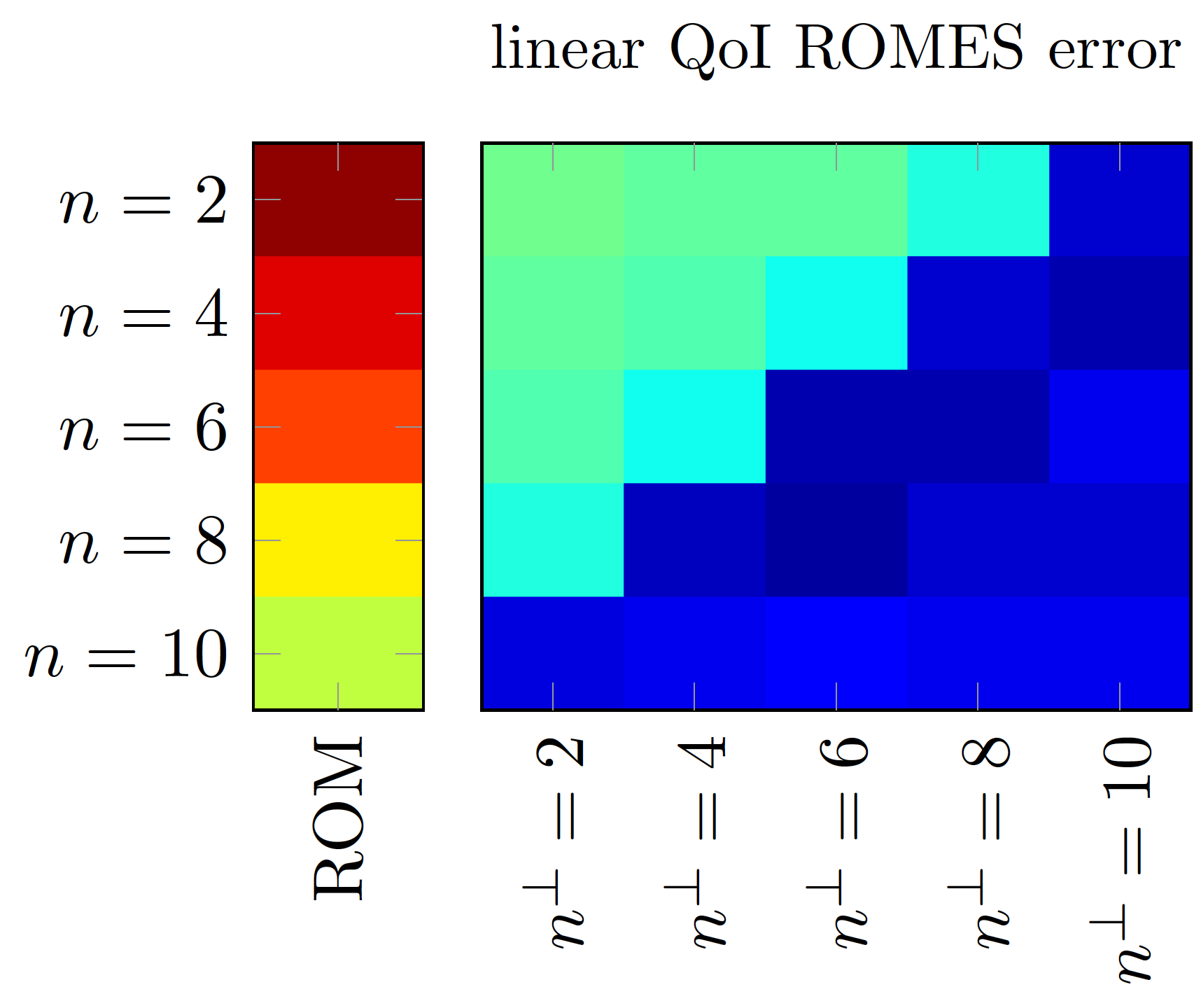}
  \includegraphics[height=0.35\textwidth]{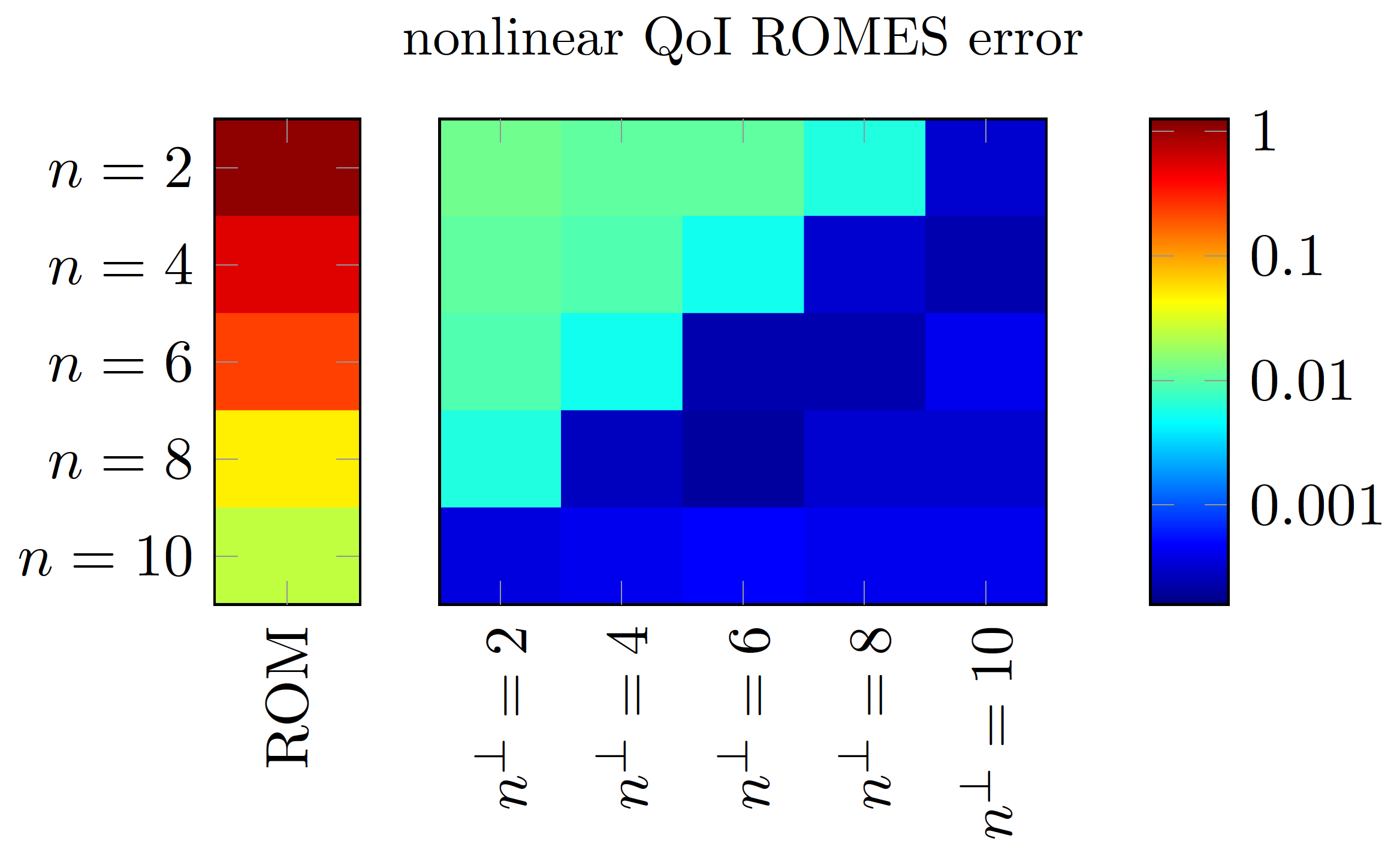}
 }
 \caption{\label{fig:output_error_plot}
 \textit{Test case 1}.
 FVU values associated with the ROM-computed quantities of interest $\qoiROMArg{1}$ and $\qoiROMArg{2}$ and the expected value of
   the ROMES-corrected quantity of interest $\expectation{\qoiModelArg{1}}$ and $\expectation{\qoiModelArg{2}}$
   for
 several values of the reduced-subspace dimension $\nrb$ and $\nrbperp$.
 Here, we
 set $\lossArg{\hyperparams}{i}{j} =
 \lossPredictionArg{\hyperparams}{i}{j}$,
   $\card{\paramROMES} = 1000$,
	 and $\nrbDual=32$.
   }
\end{figure}

Figure \ref{fig:output_error_plot} reports the
associated FVU values, with the FVU defined as
\begin{equation}
	\FVUQoiArgs{i}{q}\defeq
	\frac{\sum_{\param\in\paramEval}(\qoiFOMArg{i}(\param)-q(\param))^2}
		{\sum_{\param\in\paramEval}
		(
\qoiFOMArg{i}(\param)-
\qoiFOMAvgArg{i}
		)^2},
\end{equation}
where  $\qoiFOMAvgArg{i}$ denotes the mean value of the quantity of interest
$\qoiFOMArg{i}(\param)$ for $\param\in\paramEval$.
These plots demonstrate that the proposed method significantly reduces the
quantity-of-interest error without the need for prescribing the quantities of
interest in the offline stage (see Remark \ref{rem:qoiNotNeeded}), and
performance is improved as the dual-basis dimension $\nrbDual$ increases.

\subsubsection{Computational efficiency}\label{sec:paretoLinear}

The previous results in this section illustrate the ability of the proposed
method to reduce errors with respect to a `ROM-only' approach (i.e., a method
that executes only Step \ref{step:ROMsol} in Algorithm \ref{alg:online}).
However, the method achieves this at an increased online cost, as it
additionally executes Steps \ref{step:dual}--\ref{step:statmodels} in
Algorithm \ref{alg:online}; the dominant additional cost arises from the need
to compute the approximate dual solutoins in Step \ref{step:dual} (see Remark
\ref{rem:statmodel}).  However, we note that regardless of computational cost,
the proposed method yields a statistical model of the full-order model state
and quantity of interest, while a ROM-only approach does not. Thus, even with
increased cost, the proposed method is more amenable to integration within
uncertainty-quantification applications.  Nonetheless, we now perform an
assessment of the computational efficiency of a `ROM-only' approach and the
proposed method.

\begin{figure}[h!]
 \centerline{
 \includegraphics{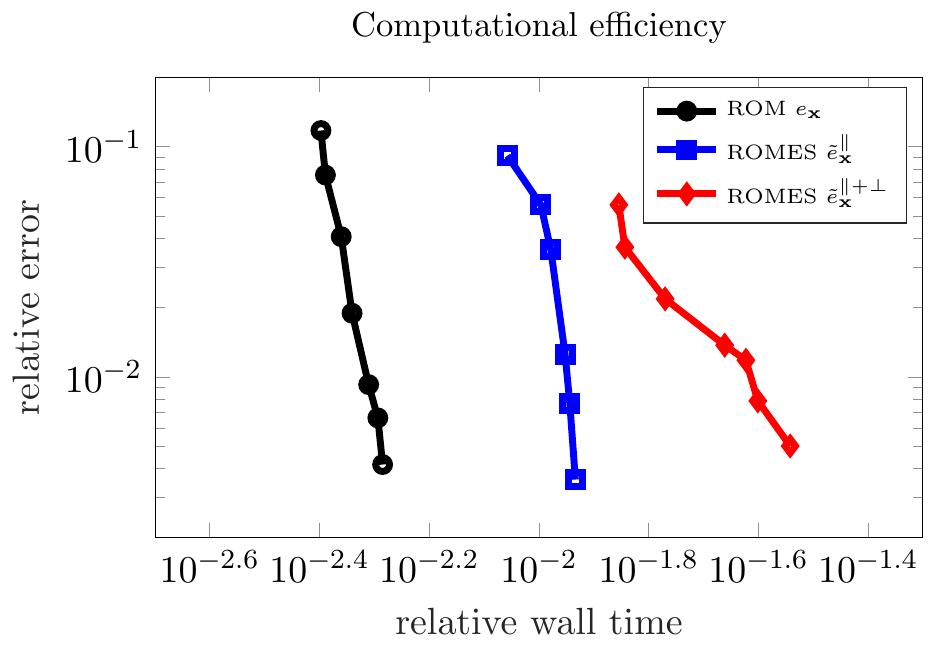}
 }
 \caption{\label{fig:comparison_performance}
	 \textit{Test case 1}. Computational efficiency of the `ROM-only' approach
	 (ROM $\errorROM$), the proposed method with a ROMES in-plane correction only (ROMES
	 $\errorROMESParallel$), and the proposed method with both an in-plane and out-of-plane
	 correction (ROMES $\errorROMESParallelPerp$) over a range of method
	 parameters, and associated Pareto fronts. Here, the relative wall time is
	 reported relateve to that incurred by the
	 full-order model.  Note that the ROM-only
	 approach is Pareto dominant due to the fact that the residual is linear in
	 its first argument for this problem (i.e., Eq.~\eqref{eq:reslinear} holds).
	 }
\end{figure}

To perform this assessment, we subject the `ROM-only' method, the proposed
method with a ROMES in-plane correction only, and the proposed method with
both an in-plane and out-of-plane correction to a wide range of parameter
values. In particular, we consider all combinations of
$\nrb\in\{1,\ldots,20\}$, $\nrbDual\in\{\nrb,\ldots,\nrb+15\}$ (not relevant
to the `ROM-only' method), and $\nrbperp\in\{1,\ldots,15\}$ (not relevant to
the `ROM-only' method or the proposed method with in-plane correction only).
Figure \ref{fig:comparison_performance} reports these results.  For each value
of these parameters, we compute both the relative error and the wall
time for the simulations relative to that incurred by the full-order model (as
averaged over all online points $\param\in\paramEval$).  The relative errors
for the ROM-only approach, the proposed method with a ROMES in-plane
correction only, and the proposed method with both an in-plane and
out-of-plane correction correspond to $\errorROM$ (Eq.~\eqref{eq:ROMerror}),
$\errorROMESParallel$ (Eq.~\eqref{eq:ROMESerror}), and
$\errorROMESParallelPerp$ (Eq.~\eqref{eq:ROMESerrorout}), respectively.  The
figure reports a Pareto front for each method, which is characterized by the
method parameters that minimize the competing objectives of relative error and
relative wall time.

These results show that the `ROM-only' approach is Pareto dominant, which is
likely due to the fact that the residual is linear in its first argument for
this problem (i.e., Eq.~\eqref{eq:reslinear} holds).
Because of this, the cost
of the dual ROM solves in Step \ref{step:dual} of Algorithm \ref{alg:online} is
similar to that of the (primal) ROM solve executed in Step \ref{step:ROMsol}
of Algorithm \ref{alg:online}. Thus, in this case, it
is always computationally more efficient to employ a larger ROM dimension
$\nrb$ than to approximate the errors according to the proposed technique if
we are only interested in minimizing the FVU.
As described in Section \ref{seq:stateqoimodels}, we expect the proposed method to be
most effective when the ROM equations are nonlinear, as the dual problems
remain linear in this case, thus allowing Steps
\ref{step:dual}--\ref{step:statmodels} to be computationally inexpensive
relative to Step \ref{step:ROMsol}. The next set of experiments will highlight
this fact.

Nonetheless, we again emphasize that the `ROM-only' approach does not generate
a statistical model of the FOM state or the FOM quantity of interest, while
the proposed approach does provide this. Thus, even in the linear case, the
proposed method may still be considered more amenable to integration with
uncertainty quantification than the `ROM-only' approach, as the proposed
method provides a mechanism to quantify the ROM-induced uncertainty.

%


\subsection{Test case 2: nonlinear mechanical response} \label{sec:5_2}

We now assess the proposed method on a problem characterized by a residual that is
nonlinear in its first argument. In particular, we consider a static,
nonlinear mechanical-response problem in three spatial dimensions.  We consider a Saint
Venant--Kirchhoff material, whose strain-energy function is given by
$$
{W}(\mathbf{E}) = \frac{\lambda_1}{2} \mathrm{tr}(\mathbf{E})^2 + \lambda_2
\mathrm{tr}(\mathbf{E}^2),
$$
where $\mathbf{E} \in \RR{3\times 3}$ denotes the Lagrangian Green strain
tensor and $\lambda_1$ and $\lambda_2$ denote Lam\'e constants
  $$
	\lambda_1 = \frac{ \nu E}{(1 + \nu)(1 - 2 \nu)} \qquad \text{and}\qquad\lambda_2 = \frac{ E}{2(1+\nu)}.
   $$
Defining the deformation gradient tensor as $\mathbf{F} = \mathbf{I} +  \nabla
\mathbf{u}\in \RR{3\times 3}$, where  $
\mathbf{u} = \mathbf{u}(\xx;\param) \in
\RR{3}$ denotes
 the deformation, we obtain the Piola tensor
$$
\mathbf{P} = \lambda_1 \mathrm{tr}(\mathbf{E}) \mathbf{F} + 2 \lambda_2 \mathbf{F} \mathbf{E}.
$$
The shear test on the domain $\spatialDomainArg{0} = [0,1]^3$
reads as follows:
 find $\mathbf{u}$ satisfying
\begin{equation} \label{eq:test2}
  \begin{cases}
    \div( \mathbf{P}(\mathbf{u} ; \param ) ) = \bf{0} & \xx \in
		\spatialDomainArg{0} \\
    \mathbf{P}(\mathbf{u} ; \param ) \mathbf{n}(\xx) = \mu_3 \mathbf{n}_z & \xx \in \Gamma_N \\
    \mathbf{P}(\mathbf{u} ; \param ) \mathbf{n}(\xx) = \bf{0} & \xx \in
		\Gamma_{N,\mathrm{free}} \\
		\mathbf{u} = \bf{0} & \xx \in \Gamma_D,
  \end{cases}
\end{equation}
where $\mathbf{n}(\xx) \equiv [n_x(\xx),n_y(\xx),n_z(\xx)]$ denotes the
outward unit normal.
We consider $\nparam=3$
parameters comprising
  the Young's modulus $\mu_1=E \in [6 \times 10^4, 8 \times 10^4]$,
  the Poisson coefficient $\mu_2=\nu \in [0.3,0.45]$, and
  the external-load magnitude $\mu_3 \in [1 \times 10^3, 2.5 \times 10^3]$.

\begin{figure}[h!]
\centerline{
\includegraphics[width=0.215\textwidth]{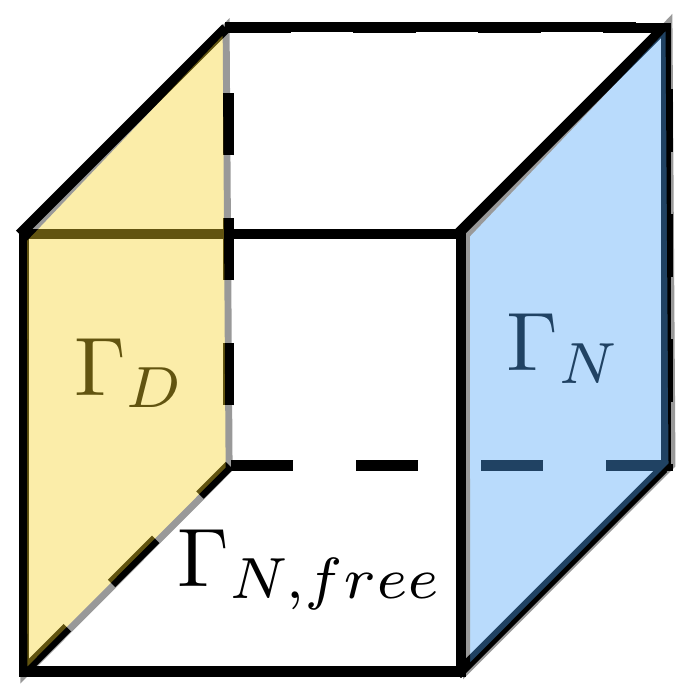} \hspace{0.85cm}
\includegraphics[width=0.215\textwidth]{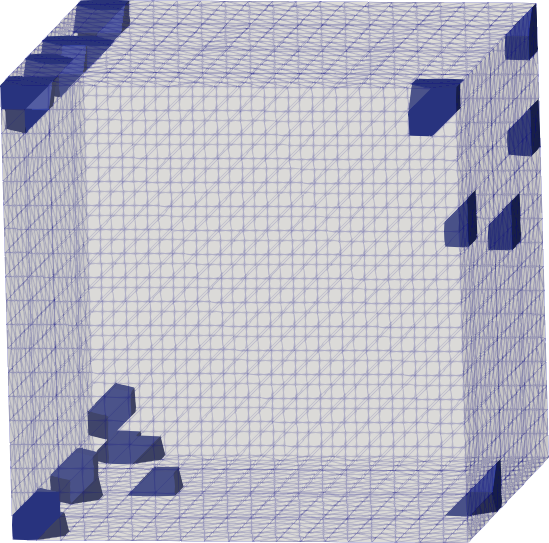}  \hspace{0.85cm}
\includegraphics[width=0.215\textwidth]{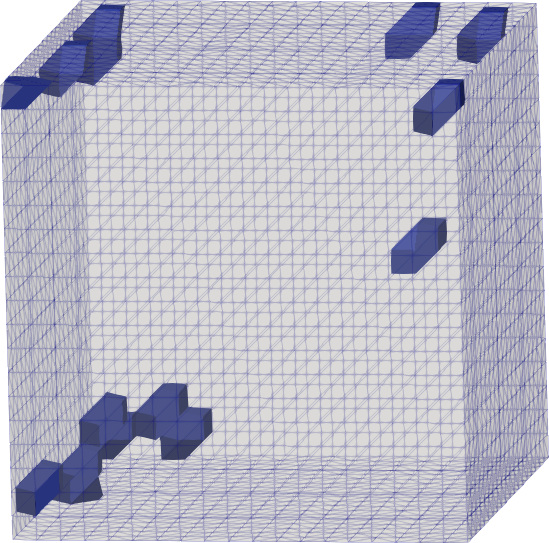}}
 \caption{
 \textit{Test case 2.}
 Geometry and boundary faces (left) and reduced meshes
for DEIM approximation of residual vectors (center) and MDEIM approximation of Jacobian matrices.}
 \label{fig:geometry}
\end{figure}

\begin{figure}[h!]
\centering
\includegraphics[width=\textwidth]{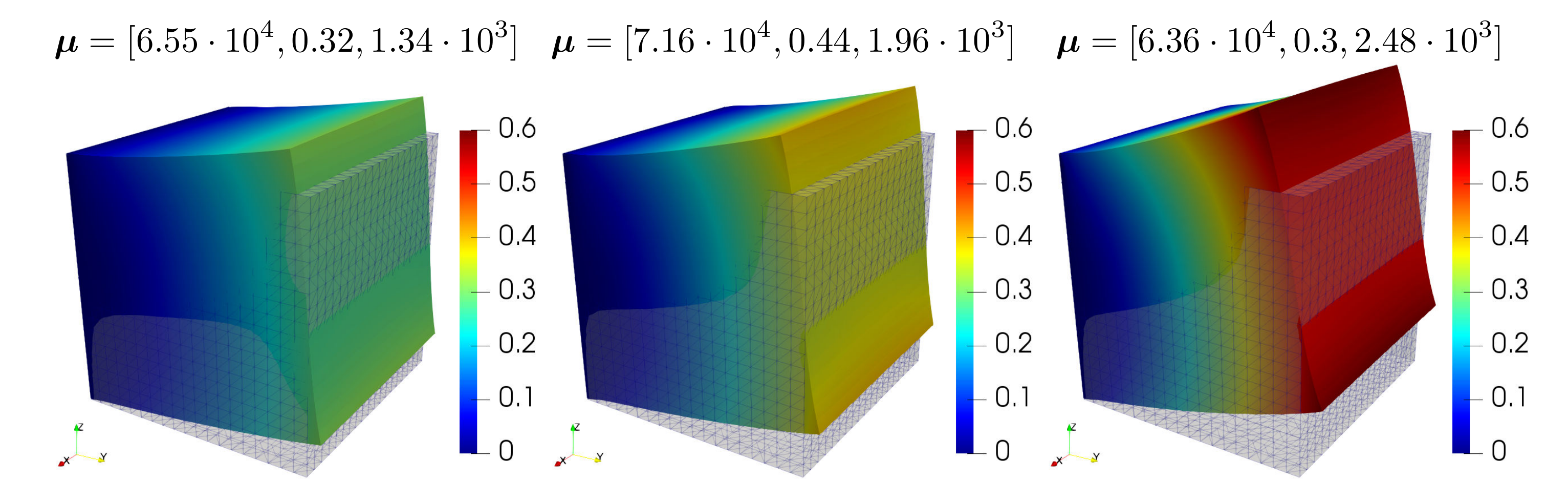}
 \caption{
\textit{Test case 2.}
 Finite element approximation of problem \eqref{eq:test2} for three different values of the parameter vector $\param$. }
 \label{fig:test2examples}
\end{figure}

We discretize the spatial domain using the finite-element method on a
conformal computational mesh given by  $41154$ tetrahedra and linear finite
elements. This yields FOM governing equations of the form
\eqref{eq:discretizedPB} with $\ndof=8000$ degrees of freedom with the
residual $\res$ nonlinear in its first argument.

We execute the offline stage using Algorithm \ref{alg:offline} as follows.
We construct the training-parameter sets by drawing uniform random samples
from the parameter domain $\paramDomain$.  We set $\card{\paramDual} = 10$,
while $\card{\paramROMES}$ varies across experiments.  In Step
\ref{step:basis}, we apply POD to FOM solutions computed at parameter
instances $\paramTrain$ with $\card{\paramTrain}=10$.  We  employ Galerkin
projection such that $\testrbmat = \rbmat$; the reduced-subspace dimension
$\nrb$ varies across experiments.
Because the residual is nonlinear in its first argument, we require a
hyper-reduction method to ensure that solving the ROM equations incurs an
$\ndof$-independent computational complexity.  For this purpose, we apply the
DEIM method \cite{Barrault2004,Chaturantabut2010} to approximate the nonlinear
component of the residual, which comprises the sum of a nonlinear component
and a linear component (the boundary conditions). For each value of the
reduced dimension $\nrb$, we collect snapshots of this nonlinear component
evaluated at the ROM solution (without hyper-reduction) at 30 parameter
instances (which includes $\paramTrain$), and we truncate the POD basis such
that it preserves $1- 1\times 10^{-8}$ of the relative statistical energy.
Step \ref{step:rbmatperp}
constructs the basis matrix $\rbmatperp$ from the discarded POD modes; the
out-of-plane subspace dimension $\nrbperp$ also depends on the particular
experiment.  In Step \ref{step:dualbasis}, we construct a single shared trial
dual basis matrix $\rbmatDual$ (i.e., $\rbmatDualArg{i} = \rbmatDual$,
$i=1,\ldots,\nrbtot$) by combining snapshots from from all $\nrbtot$ dual
solves executed at parameter instances $\param\in\paramDual$ (see Remark
\ref{rem:uniqueShared} in Section \ref{sec:romDWR}).  We also employ Galerkin
projection for the dual problem such that
$\testrbmatDualArg{i}=\testrbmatDual = \rbmatDual$, $i=1,\ldots,\nrbtot$; the
dual-basis dimension $\nrbDual$ also varies across experiments.
Because the system matrix in the dual ROM equations
\eqref{eq:dualinplaneROM} exhibits non-affine parameter dependence (but the
right-hand side is linear and parameter-independent), we apply MDEIM
\cite{negri2015efficient} to approximate the system matrix; because the
right-hand side is linear and parameter-independent, it does not require
hyper-reduction.  For MDEIM, we collect snapshots of the system matrix in an
identical way to the DEIM snapshot-collection procedure described above, and
we use the same truncation criterion.

For the online stage, we execute Algorithm \ref{alg:online} for all parameter
instances in $\paramEval$, which comprises $\card{\paramEval} = 1000$ values
drawn uniformly at random from $\paramDomain$. The remaining inputs to
Algorithm \ref{alg:online} result from the  outputs of Algorithm
\ref{alg:offline}. Note that we use DEIM when dealing with the ROM equations
\eqref{eq:romEq} in Step \ref{step:ROMsol} and MDEIM when assembling the dual ROM equations
\eqref{eq:dualinplaneROM} in Step \ref{step:dual}.

\subsubsection{ROMES model validation}

As in Section \ref{sec:GPvalTest1},
we first consider statistical validation of the ROMES models, i.e., condition
\ref{prop:validated} in Section \ref{sec:statisticalmodel}.
When constructing
the Gaussian processes in Step \ref{step:ROMESconstruction} according to the
description in Section \ref{sec:GPROMES}, we
define the set of candidate hyperparameter values $\hyperparamsSet$ as
{ uniform full-factorial sampling in each
		hyperparameter dimension characterized by 10
		equispaced values within the limits $\variance\in [0.01 \sigma_{t} , 0.25 \sigma_{t}]$,
		$\signalsd\in [ 0.1 \sigma_{t}  , \sigma_{t} ]$, and $\width\in [ 0.001
		\sigma_{t} , 0.1 \sigma_{t} ]$, with $\sigma_{t}$ denoting
		the standard deviation of the data $\{
			\stateErrorRedArg{i}(\param))\}_{\param\in\paramROMES}$.

We first consider using the negative log-likelihood loss function
$\lossArg{\hyperparams}{i}{j} = \lossloglikelihoodArg{\hyperparams}{i}{j}$
defined in Eq.~\eqref{eq:loglikelihoodloss} for hyperparameter selection.  Figure
\ref{fig:GPS} reports the resulting ROMES models constructed for the first
error generalized coordinate using a training set $\card{\paramROMES} =
400$ with two values for the dual-subspace dimension $\nrbDual$.  We note that
for $\nrbDual=4$, the data appear to be skewed and the resulting Gaussian
process exhibits large variance. By increasing the dual-subspace dimension to
$\nrbDual=12$,  the feature becomes higher quality and thus leads to a
lower-variance Gaussian process that qualitatively captures the relationship
between the error indicator and error generalized coordinate well. We now
investigate this further.

\begin{figure}[ht]
 \centerline{
 \includegraphics{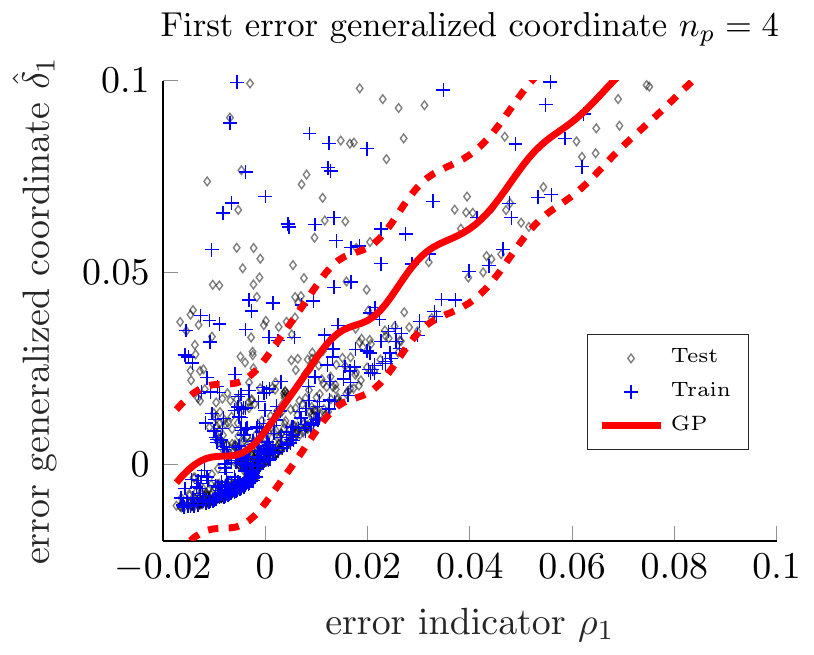}
 \includegraphics{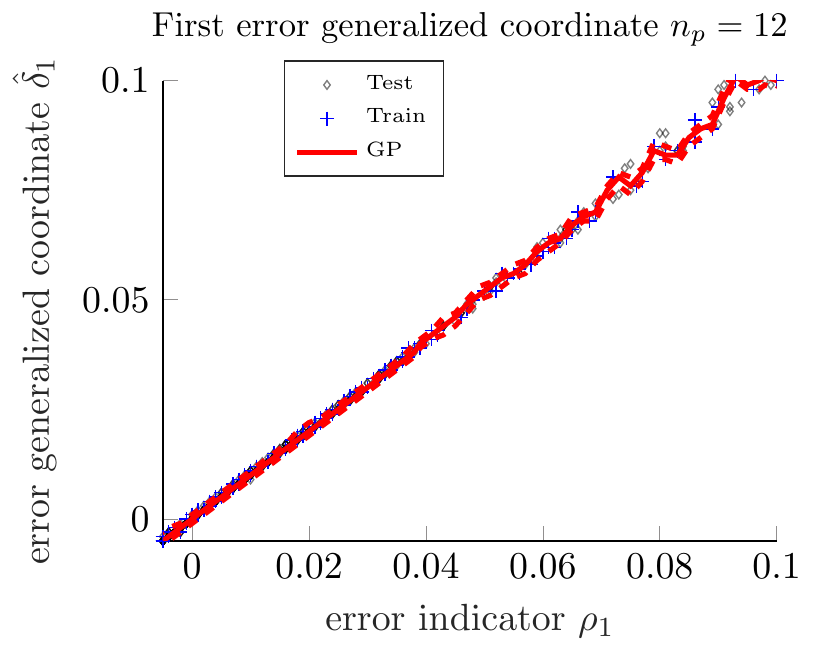}
 }
	\caption{\label{fig:GPS}  \textit{Test case 2}. ROMES models
	constructed for the first error generalized coordinate with different
	dual-subspace
dimensions $\nrbDual$.
  The solid line
  represents the GP mean; the dashed lines represent the limits of the 99\%
  prediction interval; the grey crosses represent data related to prediction
  points $\param\in\paramEval$, while the blue crosses represent training data related to prediction
  points $\param\in\paramROMES$.
	We have employed $\nrb=3$, $\nrbperp = 0$, and we have selected
	hyperparameters according to Eq.~\eqref{eq:hyperparamsOpt} with
	$\lossArg{\hyperparams}{i}{j} = \lossloglikelihoodArg{\hyperparams}{i}{j}$,
	and a training set with $\card{\paramROMES} = 400$.
 }
\end{figure}

We assess the effect of the number of training-parameter instances
$\card{\paramROMES}$ on prediction accuracy, as measured by the fraction of
variance unexplained (FVU) $\FVUArg{i}$ defined in Eq.~\eqref{eq:FVUArgDef},
the validation frequency $\validationFrequency{i}{\omega}$ defined in
Eq.~\eqref{eq:validationFrequencyDef}, and the Komolgorov--Smirnov (KS)
statistic, which quantifies the maximum discrepancy between the CDR of
$\normal{0}{1}$ and the empirical CDF of the standardized samples
$\{\stateErrorRedArg{i}(\param) -
\overallmeanArg{i}{\indicatorArg{i}(\param)})/\overallSTDArg{i}{\indicatorArg{i}(\param)}\}_{\param\in\paramEval}$.

Table \ref{Fig:ValidationFrequencyInplaneTest2} reports these results, which
show that employing $\card{\paramROMES}=400$ is sufficient for the test FVU to
have reasonably stabilized; thus, subsequent experiments in this
section set $\card{\paramROMES}=400$. As observed in Section \ref{sec:GPvalTest1}, the
converged values of FVU are quite small, but the converged prediction levels
are not all correct.  For example, $\validationFrequency{1}{0.8} = 0.9200$
even though this value should be 0.8.
This again occurs because the data are not Gaussian; Figure
\ref{fig:histGPMSEtest2} shows this.
In fact, the data do not pass Shapiro--Wilk normality test, as they yield
values of f $0.41$ for the first error generalized coordinate and of $0.46$ for the second, which again implies that it will not be possible to
achieve statistical validation in every possible metric if we employ
Gaussian-process regression.  This motivates the need for tailored loss functions as
described in Section \ref{sec:GPROMES}.

 \begin{table}[h!tb]
 \centering
 \begin{tabular}{@{} c | c c c c | c c c c @{}}
   \toprule
		error index $i$& \multicolumn{4}{c |}{$1$} & \multicolumn{4}{c}{$2$} \\
  $\card{\paramROMES}$   & $100$ & $200$ & $300$ & $400$ & $100$ & $200$ & $300$ & $400$ \\
  \midrule
		FVU &  0.0035  &  0.0036   & 0.0031  &  0.0028  &  0.0025   & 0.0032   & 0.0011   & 0.0008   \\
     $\validationFrequency{i}{0.8}$   & 0.8983 &   0.9167  &  0.9233  &  0.9200  &  0.9217   & 0.8917  &  0.9183   & 0.9117  \\
		 $\validationFrequency{i}{0.9}$   & 0.9150  &  0.9200 &   0.9317 &   0.9283  &  0.9300  &  0.9033  &  0.9350  &  0.9233 \\
		 $\validationFrequency{i}{0.95}$  & 0.9250  &  0.9300  &  0.9383  &  0.9417  &  0.9333  &  0.9117  &  0.9467  &  0.9317 \\
		 $\validationFrequency{i}{0.99}$  & 0.9467  &  0.9450  &  0.9450  &  0.9500  &  0.9567  &  0.9367   & 0.9567  &  0.9483 \\
		 KS statistic  & 0.2030  &  0.2822   & 0.2804  &  0.2388  &  0.2781  &  0.3123  &  0.3439  &  0.3438  \\
   \bottomrule
 \end{tabular}
   \caption[ValidationFrequencyInplane]{ \label{Fig:ValidationFrequencyInplaneTest2}
		 \textit{Test case 2}. Convergence of error measures associated with the ROMES models
		 constructed for the first two error generalized coordinates as the number
		 of training-parameter instances $\card{\paramEval}$ increases.
		 We have employed $\nrb=2$, $\nrbperp = 0$, $\nrbDual=12$,
		 $\card{\paramEval}=600$ and
	 , and
	 have selected hyperparameters according to Eq.~\eqref{eq:hyperparamsOpt} with
	 $\lossArg{\hyperparams}{i}{j} = \lossloglikelihoodArg{\hyperparams}{i}{j}$.
	 }
 \end{table}

\begin{figure}[h!t]
 \centerline{
 \includegraphics{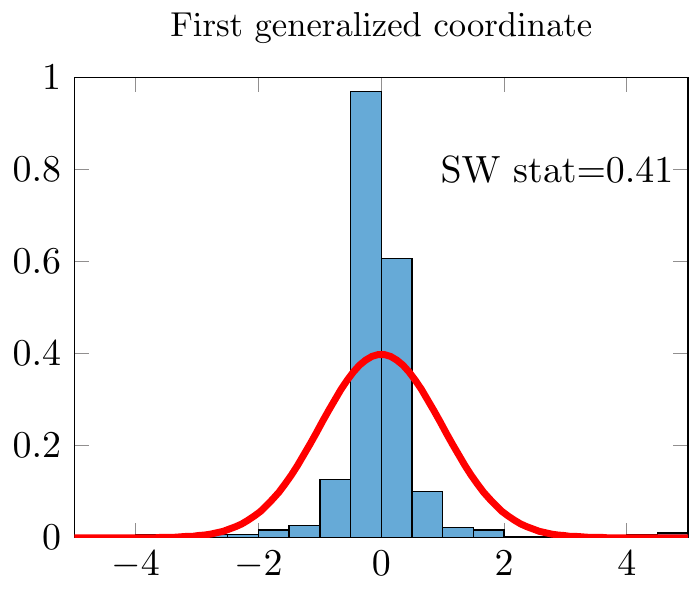}
 \includegraphics{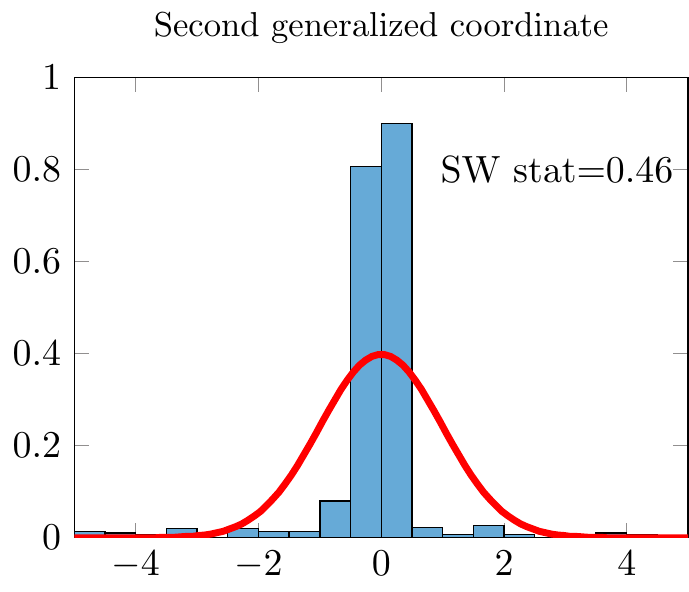}
 }
	\caption{\label{fig:histGPMSEtest2}  \textit{Test case 2}. Histogram of the
 	standardized data
 	$\{\stateErrorRedArg{i}(\param) -
 	\overallmeanArg{i}{\indicatorArg{i}(\param)})/\overallSTDArg{i}{\indicatorArg{i}(\param)}\}_{\param\in\paramEval}$,
 	$i=1,2$ (blue bar plot)
 	as compared to the PDF of the standard Gaussian distribution $\normal{0}{1}$
 	(red curve).
 	We have employed $\nrb=3$, $\nrbperp = 0$, $\nrbDual=10$, and
 	 have selected hyperparameters according to Eq.~\eqref{eq:hyperparamsOpt} with
 	 $\lossArg{\hyperparams}{i}{j} = \lossloglikelihoodArg{\hyperparams}{i}{j}$.
 	 The number of training-parameter instances is $\card{\paramROMES}=400$.
	}
\end{figure}

%
%
%



As in Section \ref{sec:GPvalTest1}, we consider five different loss
functions
$\lossArg{\hyperparams}{i}{j}$
for hyperparameter selection
(see Section \ref{sec:GPROMES}): (1) the negative log-likelihood loss
$\lossloglikelihoodNo$ (Eq.~\eqref{eq:loglikelihoodloss}), (2) the loss based on
matching the 0.80-prediction interval $\lossEightyNo$
(Eq.~\eqref{eq:lossFrequency} with $\frequency=0.80$),
(3) the loss based on
matching the 0.95-prediction interval $\lossNinetyFiveNo$
(Eq.~\eqref{eq:lossFrequency} with $\frequency=0.95$), (4) the loss based on a
linear combination of $\frequency$-prediction interval losses
$\lossPredictionNo$ (Eq.~\eqref{eq:lossPredictionNo}), and (5) the loss based
on the KS statistic $\lossKSNo$.

Table \ref{Fig:ValidationFrequencyInplaneTwoTestTwo} reports these results for
$\card{\paramROMES}=1000$.  As in Section \ref{sec:GPvalTest1}, we observe that
the loss function $\lossNo$  has a significant effect on the performance of
the resulting ROMES models according to  different statistical-validation
criteria, and can be chosen to target performence with respect to particular
criteria: $\validationFrequency{1}{0.8} = 0.81$ instead of $\validationFrequency{1}{0.8} = 0.92$  and
$\validationFrequency{1}{0.8} = 0.8250$ instead of $\validationFrequency{1}{0.8} = 0.9117$ when $\lossEightyNo$ is adopted.

We again employ
a loss function of $\lossNo=\lossPredictionNo$ in the remaining
experiments within this section
due to its favorable performance over a range of statistical-valation
criteria.

 \begin{table}[h!tb]
 \centering
 \begin{tabular}{@{} c | c  c c c c}
   \toprule
		error generalized coordinate index $i$& \multicolumn{5}{c }{$1$} \\
	 loss function $\lossArg{\hyperparams}{i}{j}$ & $\lossloglikelihoodNo$ & $\lossPredictionNo$ &
	 $\lossEightyNo$ & $\lossNinetyFiveNo$ & $\lossKSNo$ \\
  \midrule
		FVU &  0.0028  &  0.0029  &  0.0015  &  0.0017  &  0.0015 \\
           $\validationFrequency{i}{0.8}$   & 0.9200  &  0.8750  &  0.8100  &  0.9183  &  0.8133  \\
					 $\validationFrequency{i}{0.9}$   & 0.9283  &  0.9117  &  0.8567  &  0.9400  &  0.8567 \\
					 $\validationFrequency{i}{0.95}$  & 0.9417  &  0.9217  &  0.8733  &  0.9483  &  0.8767  \\
					 $\validationFrequency{i}{0.99}$  & 0.9500  &  0.9317  &  0.8950  &  0.9567  &  0.8900 \\
					 Komolgorov--Smirnov statistic  & 0.2388  &  0.1816  &  0.1179  &  0.3280  &  0.0908 \\
   \bottomrule
		error generalized coordinate index $i$& \multicolumn{5}{c }{$2$} \\
	 loss function $\lossArg{\hyperparams}{i}{j}$ & $\lossloglikelihoodNo$ & $\lossPredictionNo$ &
	 $\lossEightyNo$ & $\lossNinetyFiveNo$ & $\lossKSNo$\\
  \midrule
		FVU &  0.00082  &  0.00081  &  0.00075  &  0.00076  &  0.00076 \\
	 $\validationFrequency{i}{0.8}$   & 0.9117  &  0.8617  &  0.8250  &  0.9167  &  0.8133 \\
					 $\validationFrequency{i}{0.9}$   & 0.9233  &  0.8817  &  0.8533  &  0.9383  &  0.8383 \\
					 $\validationFrequency{i}{0.95}$  & 0.9317  &  0.8983  &  0.8667  &  0.9500  &  0.8633 \\
					 $\validationFrequency{i}{0.99}$  & 0.9483  &  0.9100  &  0.8950  &  0.9667  &  0.8867 \\
					 Komolgorov--Smirnov statistic  & 0.3438  &  0.2467  &  0.1329  &  0.3120  &  0.1135 \\
   \bottomrule
 \end{tabular}
   \caption[ValidationFrequencyInplane]{  \label{Fig:ValidationFrequencyInplaneTwoTestTwo}
 	 \textit{Test case 2}. Statistical-validation criteria for ROMES models when
	 different loss functions $\lossNo$ are employed for hyperparameter
	 selection according to Eqs.~\eqref{eq:lossi} and \eqref{eq:hyperparamsOpt}
	 in Section \ref{sec:GPROMES}.
We have employed $\nrb=2$, $\nrbperp = 0$,
	 $\nrbDual=10$, and have selected hyperparameters
	 according to Eq.~\eqref{eq:hyperparamsOpt} with
	 the specified loss function $\lossArg{\hyperparams}{i}{j}$.
	 The number of training-parameter instances is $\card{\paramROMES}=400$, while $\card{\paramEval}=600$.
	}
 \end{table}

\subsubsection{In-plane and out-of-plane error
approximation}\label{sec:test2InplaneOutofplane}

As in Section \ref{sec:inplaneoutofplane}, we now assess the ability of the
proposed method to accurately approximate the in-plane error
$\stateErrorInplane$ and the out-of-plane error $\stateErrorOutofplane$.  In
particular, we compare the mean relative ROM error $\errorROM$
(Eq.~\eqref{eq:ROMerror}) with the mean relative ROM error after applying the
in-plane ROMES correction $\errorROMESParallel$ (Eq.~\eqref{eq:ROMESerror})
and the mean relative projection error $\errorROMParallel$
(Eq.~\eqref{eq:Perror}).

Figure \ref{fig:inplane_plot_test2} reports the results obtained for
$\nrbperp = 0$, $\lossArg{\hyperparams}{i}{j} =
\lossPredictionArg{\hyperparams}{i}{j}$, and
$\card{\paramROMES} = 400$ training-parameter instances, and a
range of values for $\nrb$ and $\nrbDual$.
We first note that the
mean relative ROM error $\errorROM$ is relatively close to the (optimal)
mean relative projection error $\errorROMParallel$, implying that the in-plane
error is quite small for this particular problem; however the proposed method is
indeed able to bridge this gap, as the in-plane ROMES correction
$\stateErrorInplaneRedModel$ enables $\errorROMESParallel$ to be nearly equal
to the optimal value
$\errorROMParallel$.


\begin{figure}[h!t]
 \centerline{
 \includegraphics{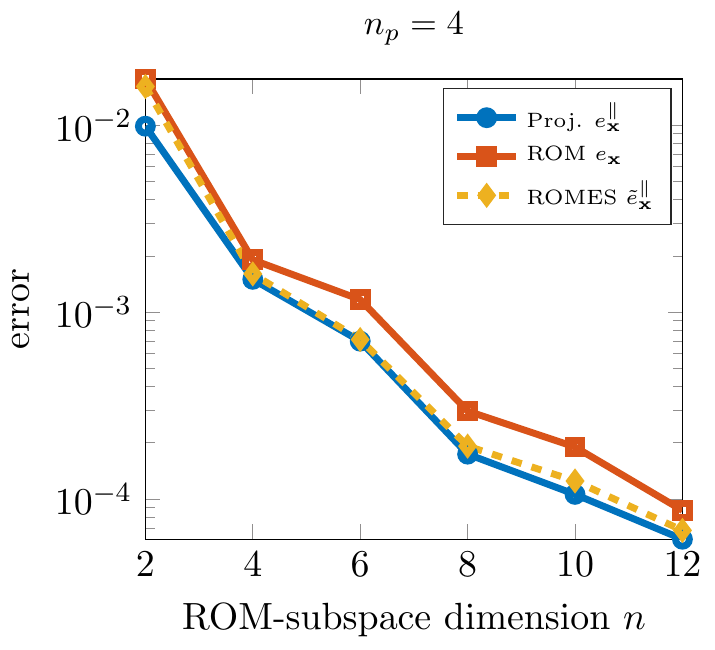}
 \includegraphics{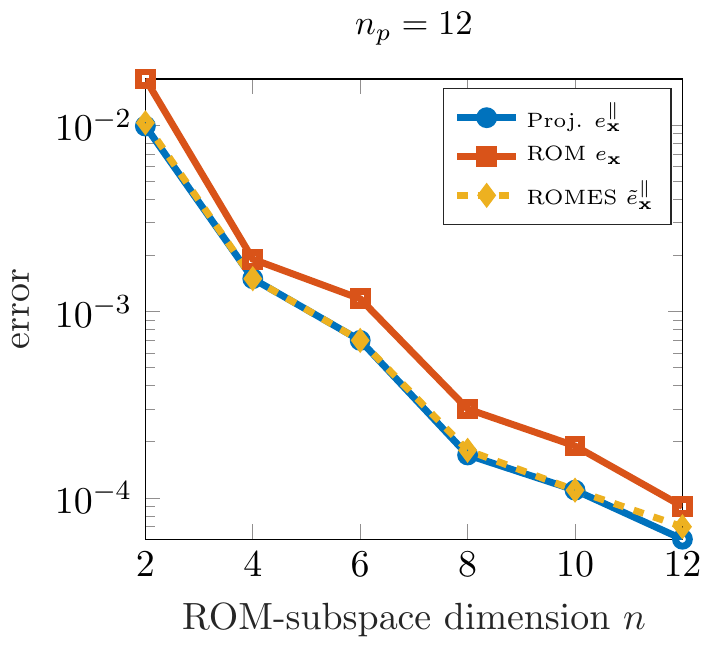}
 }
 \caption{\label{fig:inplane_plot_test2}
 \textit{Test case 2}.
 Mean relative ROM error $\errorROM$ (red),  mean relative ROM error
	after applying the in-plane ROMES correction $\errorROMESParallel$ (yellow),
	and
	mean relative projection error $\errorROMParallel$ (blue) for a varying
	reduced-subspace dimension $\nrb$ and dual-subspace dimension $\nrbDual$.
	Here, we
set $\nrbperp = 0$, $\lossArg{\hyperparams}{i}{j} =
\lossPredictionArg{\hyperparams}{i}{j}$, $\card{\paramEval}=600$ and
$\card{\paramROMES} = 400$.
	}
\end{figure}

We again perform a similar analysis for the out-of-plane error by comparing
the mean relative ROM error $\errorROM$ with the mean relative ROM error after
applying both the in-plane and out-of-plane ROMES corrections
$\errorROMESParallelPerp$ (Eq.~\eqref{eq:ROMESerrorout}) and the mean relative
projection error $\errorROMParallelPerp$ (Eq.~\eqref{eq:Perrorout}). Figure
\ref{fig:outofplane_plot_test2} reports these results for various values of
the reduced-subspace dimension $\nrb$, the dual-basis dimension $\nrbDual$,
and the out-of-plane subspace dimension $\nrbperp$. We again note that the
ROMES correction nearly eliminates both the in- and out-of-plane errors,
as $\errorROMESParallelPerp$ nearly achieves the optimal value of
$\errorROMParallelPerp$ for all considered parameters.

\begin{figure}[ht]
 \centerline{
 \includegraphics{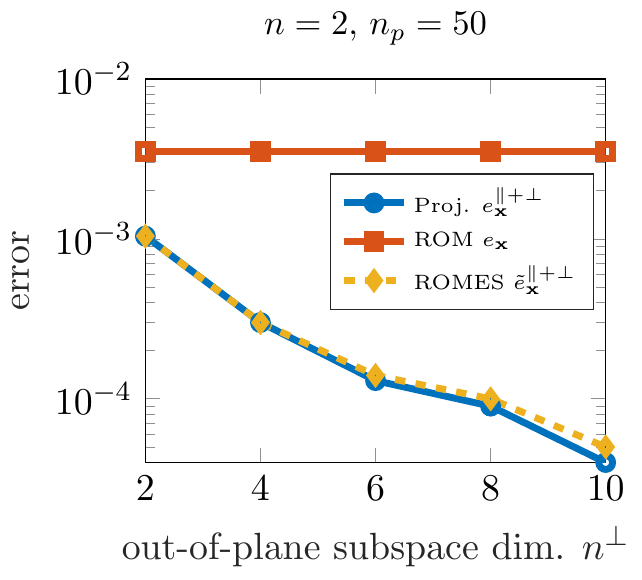}
 \includegraphics{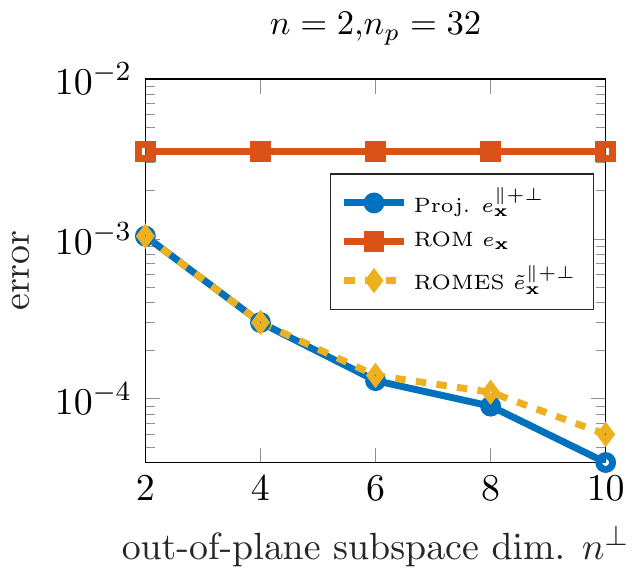}
 }
 \caption{\label{fig:outofplane_plot_test2}
 \textit{Test case 2}.
 Mean relative ROM error $\errorROM$ (red),  mean relative ROM error
	after applying the in-plane and out-of-plane ROMES corrections
	$\errorROMESParallelPerp$ (yellow),
	and
	mean relative projection error $\errorROMParallelPerp$ (blue) for a varying
	reduced-subspace dimension $\nrb$, dual-subspace dimension $\nrbDual$, and
	out-of-plane subspace dimension $\nrbperp$.
	Here, we
set $\lossArg{\hyperparams}{i}{j} =
\lossPredictionArg{\hyperparams}{i}{j}$,
$\card{\paramROMES} = 400$, $\card{\paramEval} = 600$
}
\end{figure}

\subsubsection{Quantity-of-interest error approximation}

As in Section \ref{sec:qoiApproxTest1}, we now consider the ability of the
proposed method to construct statistical models of quantities of
interest $\qoiModel$ as proposed in Section \ref{seq:stateqoimodels}.
For this purpose, we consider $\nqoi=\ndof$ quantities of interest defined by
the von Mises stress at all $\ndof$ (unconstrained) grid points in the mesh,
i.e.,
\begin{align*}
\begin{split}
	&\qoiFOM:\paramDummy\mapsto	\qoiFunc(\state(\paramDummy); \paramDummy ) =
	\\
	&\begin{bmatrix}
		\sqrt{\frac{(\sigma_{11}(\spaceDof{1}) - \sigma_{22}(\spaceDof{1}))^2 + (\sigma_{22}(\spaceDof{1}) - \sigma_{33}(\spaceDof{1}))^2 +
		(\sigma_{33}(\spaceDof{1}) - \sigma_{11}(\spaceDof{1}))^2 + 6(\sigma_{12}(\spaceDof{1})^2 + \sigma_{23}(\spaceDof{1})^2 +
		\sigma_{31}(\spaceDof{1})^2)}{2}}\\
		\vdots\\
		\sqrt{\frac{(\sigma_{11}(\spaceDof{\ndof}) -
		\sigma_{22}(\spaceDof{\ndof}))^2 + (\sigma_{22}(\spaceDof{\ndof}) -
		\sigma_{33}(\spaceDof{\ndof}))^2 +
		(\sigma_{33}(\spaceDof{\ndof}) - \sigma_{11}(\spaceDof{\ndof}))^2 +
		6(\sigma_{12}(\spaceDof{\ndof})^2 + \sigma_{23}(\spaceDof{\ndof})^2 +
		\sigma_{31}(\spaceDof{\ndof})^2)}{2}}
	\end{bmatrix},
\end{split}
\end{align*}
where $\sigma_{ij} = (\mathbf{P}(\state)(\mathbf{I}+\nabla \state)^T)_{ij}$,
$i,j=1,2,3$ and $\spaceDof{i}$ denotes the $i$th (unconstrained) grid point in the
computational mesh. This is an example of high-dimensional quantity of
interest whose error cannot be modeled tractably using the original ROMES
method \cite{kevin:romes}, as this would require constructing $\ndof$ separate
Gaussian-process models. We emphasize that these quantities of interest were
not specified during the offline stage (see Remark \ref{rem:qoiNotNeeded}).

We set the loss function to $\lossArg{\hyperparams}{i}{j} =
\lossPredictionArg{\hyperparams}{i}{j}$, the
number of training-parameter instances to $\card{\paramROMES} = 400$, the dual-basis dimension to $\nrbDual=32$.
Figure \ref{fig:MAXVM_output_scatter_plot} plots
the maximum value of the FOM-computed quantity of interest
$\max_i(\qoiFOMArg{i}(\param))$ versus both the
the
maximum value of the ROM-computed quantity of interest
$\max_i(\qoiROMArg{i}(\param))$ and the maximum value of the expected value of
the ROMES-corrected quantity of interest
$\max_i(\expectation{\qoiModelArg{i}(\param)})$ for
several values of the reduced-subspace dimension $\nrb$ and out-of-plane
subspace dimension $\nrbperp$ and for
$\param\in\paramEval$.
Figure \ref{fig:output_error_plot_chess_matrix} reports the
associated FVU values, with the FVU defined as
\begin{equation}
	\FVUQoiMaxArgs{\boldsymbol q}\defeq
	\frac{\sum_{\param\in\paramEval}(\max_{i}(\qoiFOMArg{i}(\param))-\max_i(
		q_i(\param)))^2}
		{\sum_{\param\in\paramEval}
		(
	\max_i(\qoiFOMArg{i}(\param))-
\qoiFOMAvgMax
		)^2},
\end{equation}
where  $\qoiFOMAvgMax$ denotes the mean value of the quantity of interest
$\max_i(\qoiFOMArg{i}(\param))$ for $\param\in\paramEval$.
These results show the ability of the proposed method to significantly reduce
the quantity-of-interest error without the need for prescribing the quantities
of interest in the offline stage (see Remark \ref{rem:qoiNotNeeded}). These
results also show that performance is improved as the dual-basis dimension
$\nrbDual$ increases, albeit at increased computational cost.

\begin{figure}[h!t]
 \centerline{
 \includegraphics{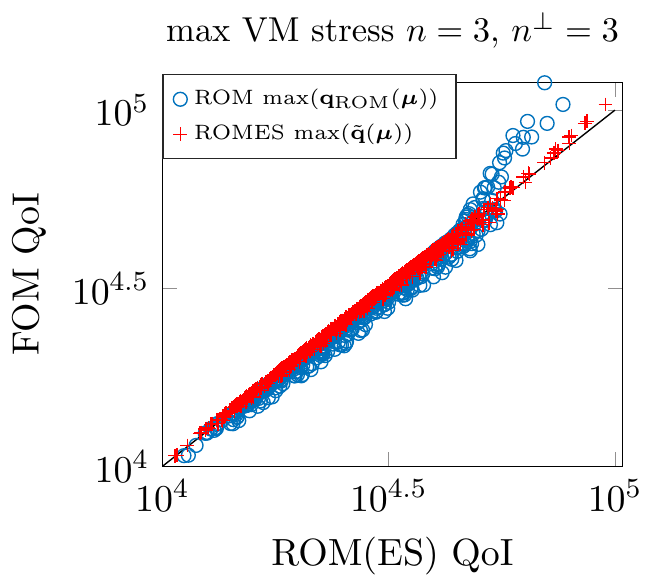}
 \includegraphics{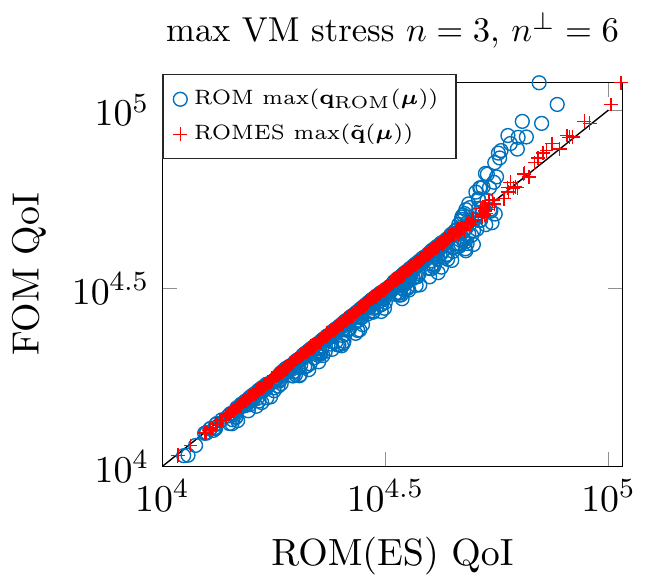}
 }
 \centerline{
 \includegraphics{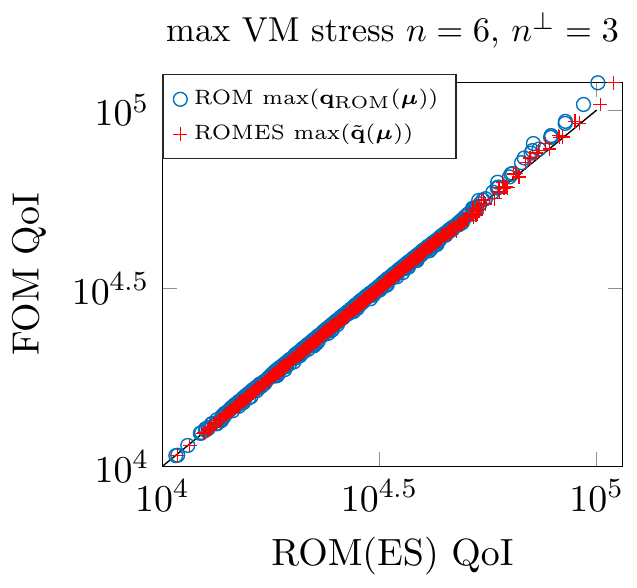}
 \includegraphics{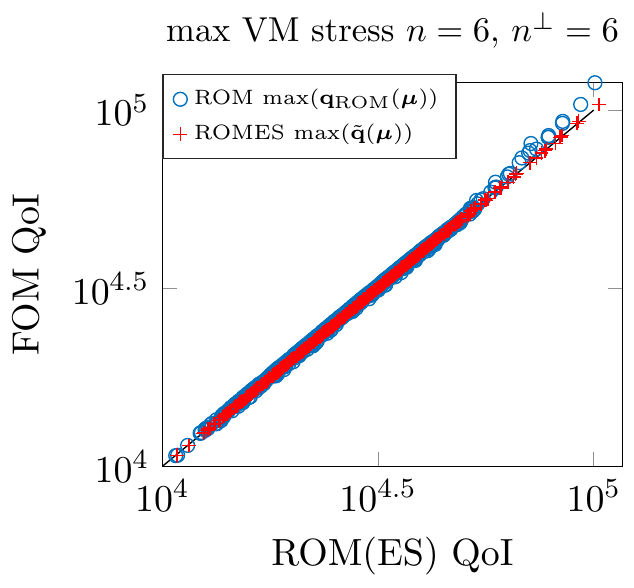}
 } \caption{\label{fig:MAXVM_output_scatter_plot}
 \textit{Test case 2}.
Scatter plot of the the maximum value of the FOM-computed quantity of interest
$\max_i(\qoiFOMArg{i}(\param))$ versus both the
the
maximum value of the ROM-computed quantity of interest
$\max_i(\qoiROMArg{i}(\param))$  (blue circles) and the maximum value of the expected value of
the ROMES-corrected quantity of interest
$\max_i(\expectation{\qoiModelArg{i}(\param)})$
	(red crosses) for
several values of the reduced-subspace dimension $\nrb$ and for
$\param\in\paramEval$. We have employed
	 $\nrbDual=32$, and have selected hyperparameters
	 according to Eq.~\eqref{eq:hyperparamsOpt} with
	 the specified loss function $\lossArg{\hyperparams}{i}{j}$.
	 The number of training-parameter instances is $\card{\paramROMES}=400$, while $\card{\paramROMES}=600$.
}
\end{figure}

\begin{figure}[h!t]
 \centerline{
  \includegraphics[height=0.35\textwidth]{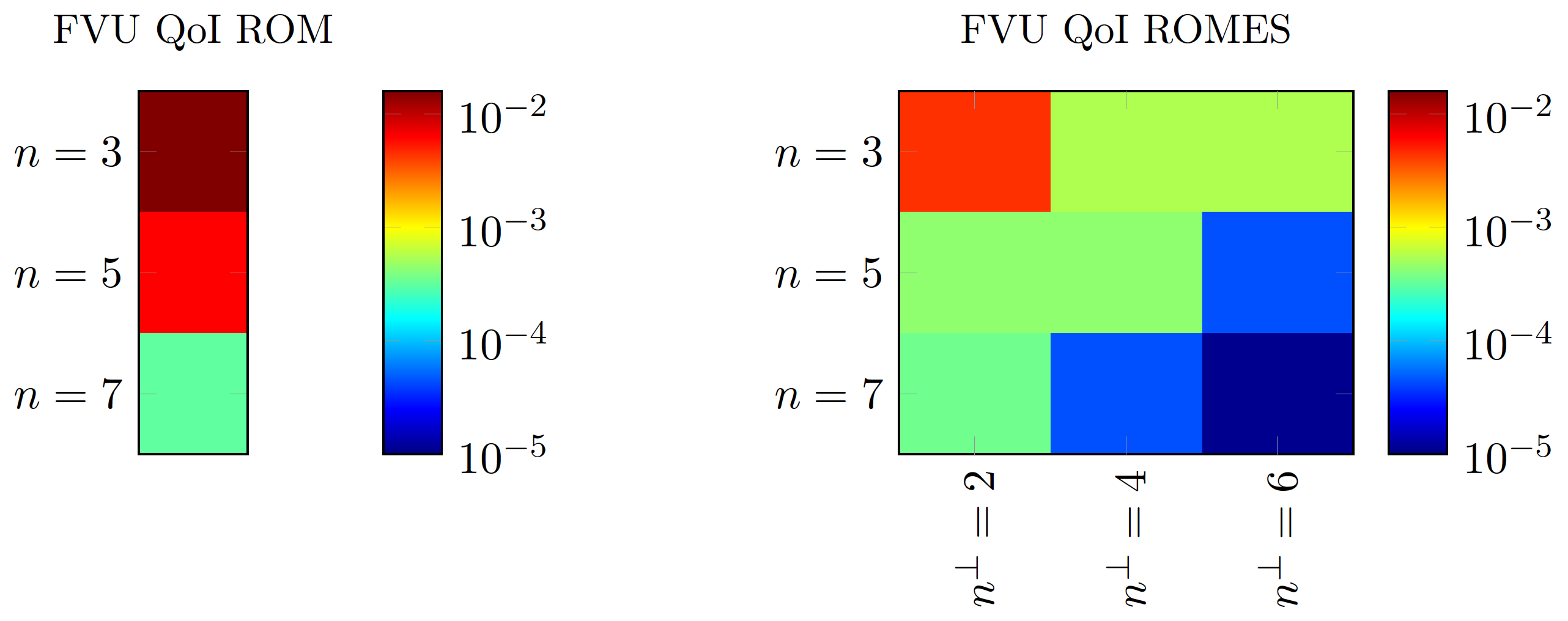}
 }
 \caption{\label{fig:output_error_plot_chess_matrix}
 \textit{Test case 2}.
 FVU values associated with the ROM-computed quantity of interest $\max_i(\qoiROMArg{i}(\param))$ and the expected value of
   the ROMES-corrected quantity of interest $\max_i(\expectation{\qoiModelArg{i}(\param)})$
   for
 several values of the reduced-subspace dimension $\nrb$ and $\nrbperp$.
 Here, we
 set $\lossArg{\hyperparams}{i}{j} =
 \lossPredictionArg{\hyperparams}{i}{j}$,
   $\card{\paramROMES} = 1000$ and $\nrbDual=32$.
   }
\end{figure}

Finally, Figure \ref{fig:VMSn1} reports the values of the the mean relative ROM error
$\errorROMQoiArg{i}$
\begin{equation} \label{eq:ROMerrorOut}
  \errorROMQoiArg{i}\defeq\expectationParam{ \frac{ | \qoiErrorArg{i}(\param)
  | }{ | \qoiFOMArg{i}( \param ) | } } ,
\end{equation}
and the mean relative ROM error
with in- and out-of-plane ROMES correction $\errorROMESParallelPerpQoiArg{i}$
\begin{equation} \label{eq:ROMESerrorOut}
  \errorROMESParallelPerpQoiArg{i}\defeq\expectationParam{ \frac{ |
  \qoiArg{i}( \param ) - \qoiFuncArg{i}(\stateROM    (\param)+
  \rbmat\expectation{\stateErrorInplaneRedModel(\param)}
+\rbmatperp\expectation{\stateErrorOutofplaneRedModel(\param)}
  ;
  \param ) | }{ |\qoiArg{i}( \param ) | } },
\end{equation}
for
$i=1,\ldots,\ndof$ as distributed over the physical domain.
We observe that applying both the in-plane and out-of-plane ROMES correction
yields a very small mean relative error, thereby illustrating the ability of
the method to accurately model the error in field quantities.

\begin{figure}[h!]
  \vspace{-0.25cm}
\centering
\includegraphics[width=0.9\textwidth]{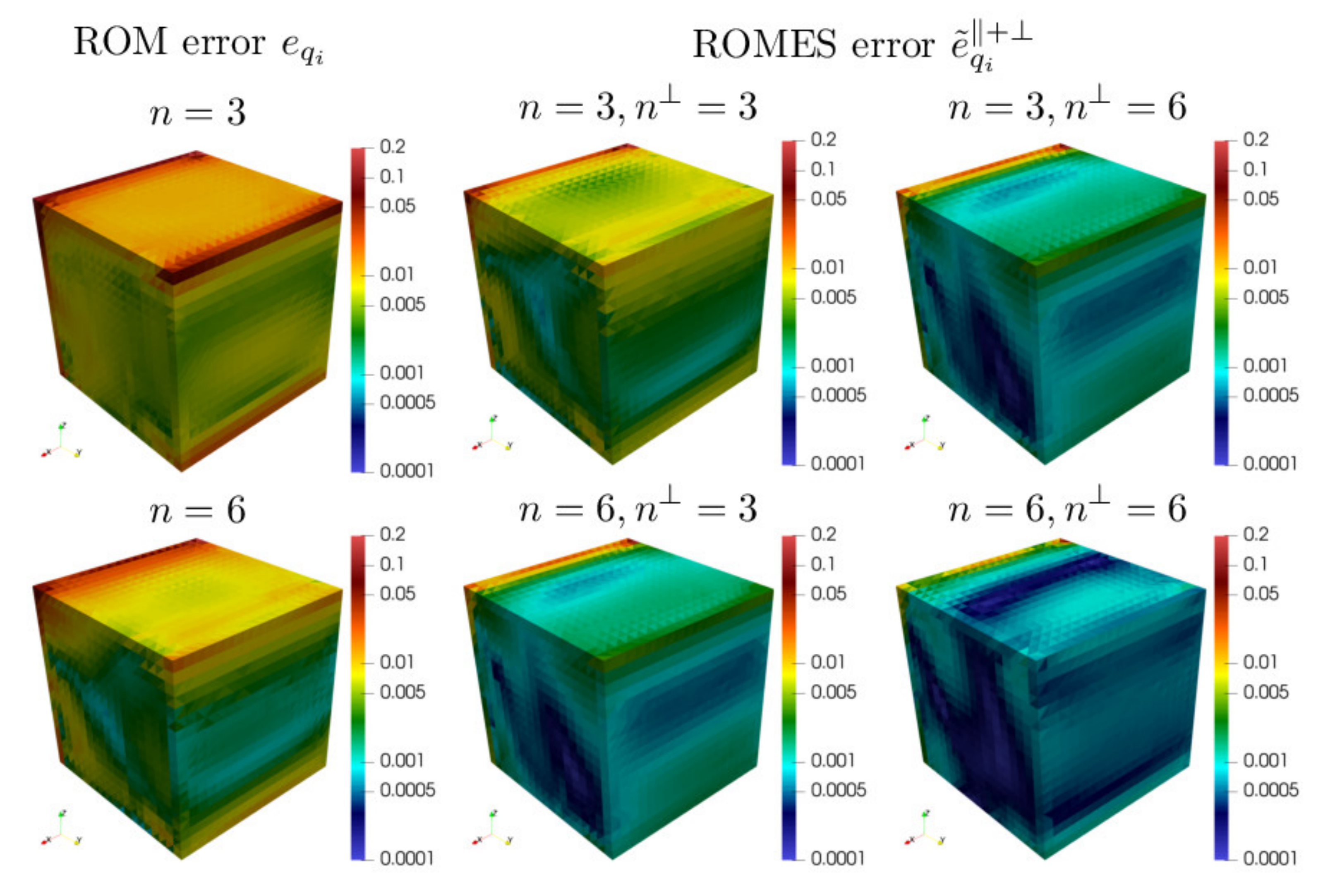}
 \caption{\textit{Test case 2.} The values of the mean relative ROM error
	$\errorROMQoiArg{i}$ (left column) and the mean relative ROM error with in- and out-of-plane ROMES
	correction $\errorROMESParallelPerpQoiArg{i}$ for $i=1,\ldots,\ndof$ (center
	and right columnes) as distributed over
the physical domain.
	Here, we
set $\lossArg{\hyperparams}{i}{j} =
\lossPredictionArg{\hyperparams}{i}{j}$,
	$\card{\paramROMES} = 400$, $\card{\paramEval} = 600$ and $\nrbDual=12$. \vspace{-0.25cm}
}
 \label{fig:VMSn1}
\end{figure}



\subsubsection{Computational efficiency}

As in Section \ref{sec:paretoLinear}, we now analyze the computational
efficiency of the proposed method with respect to a `ROM-only' approach.  As
discussed in Remark \ref{rem:statmodel}, we expect the proposed method
to yield favorable performance relative to the linear problem considered in
Section \ref{sec:5_1}, as the dual ROM equations \eqref{eq:dualinplaneROM} are always
linear in their first argument, even when the ROM equations \eqref{eq:romEq}
are nonlinear in their first argument; thus, relative to the (primal) ROM
solve, the dual solves are computationally inexpensive.

We repeat the study executed in Section \ref{sec:paretoLinear}, and subject
the `ROM-only' method, the
proposed method with a ROMES in-plane correction only, and the proposed method
with both an in-plane and out-of-plane correction
to a wide range of parameter values.
In particular, we
consider all combinations of $\nrb\in\{2,\ldots,10\}$,
$\nrbDual\in\{12,22,32\}$ (not relevant to the `ROM-only' method),
and $\nrbperp\in\{2,\ldots,10\}$ (not relevant to the `ROM-only' method or the
proposed method with in-plane correction only). Figure
\ref{fig:accuracy_vs_efficiency_test2} reports these results and associated
Pareto fronts.

These results show that the proposed method with both in-plane and
out-of-plane ROMES corrections are Pareto dominant. Specifically, for a fixed
wall time, the method yields approximately one order of magnitude in error
reduction; for a fixed error, the method yields approximately a 30\% reduction
in wall time. Furthermore, this approach provides a statistical model of the
FOM state and quantities of interest, which is not provided by the `ROM-only'
approach. Thus, the proposed method has demonstrated superior performance not
only in its computational efficiency, but also in its ability to quantify the
ROM-induced epistemic uncertainty, which is essential for rigorous integration
into uncertainty-quantification applications. We note that the proposed method
with an in-plane ROMES correction only yields similar performance to the
`ROM-only' approach. This occurs because the ROM in-plane errors are already
quite small; this was previously discussed in Section
\ref{sec:test2InplaneOutofplane} and observed in Figure
\ref{fig:inplane_plot_test2}.

%

\begin{figure}[h!t]
  \vspace{-0.25cm}
 \centerline{
 \includegraphics{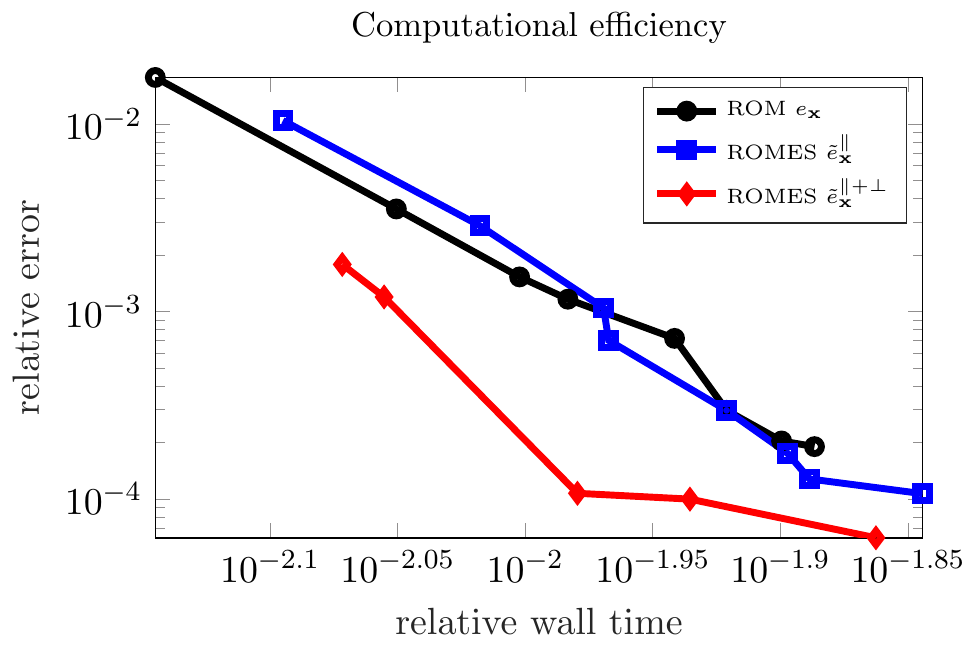}
}
 \caption{\label{fig:accuracy_vs_efficiency_test2}
 \textit{Test case 2.}
Computational efficiency of
	the `ROM-only' approach  (ROM $\errorROM$), the
	proposed method with a ROMES in-plane correction only (ROMES
  $\errorROMESParallel$), and
the proposed method
	with both an in-plane and out-of-plane correction (ROMES $\errorROMESParallelPerp$)
	 three approaches over a
	 range of method parameters, and associated Pareto fronts. Here, the
	 relative error corresponds to
	 $\errorROMESParallelPerp$ (Eq.~\eqref{eq:ROMESerrorout}),
and the wall time for the
	 simulations relative to that incurred by the full-order model.  Note that
	 the proposed method with both in-plane and out-of-plane ROMES correction is
	 Pareto dominant. \vspace{-0.25cm}
	}
\end{figure}



\section{Conclusions}

This work has proposed a technique for constructing a statistical closure
model for reduced-order models (ROMs) applied to stationary systems. The
proposed method applies the ROMES method to construct a statistical model for
the state error through constructing statistical models for the generalized
coordinates characterizing both the in-plane error (i.e., the closure model)
and a low-dimensional approximation of the out-of-plane error.
Key ingredients of the method include (1) cheaply computable error indicators
associated with a ROM-approximated dual-weighted residual (Section
\ref{sec:errorindicator}), (2) a Gaussian-process model to map these error
indicators to a random variable for the error generalized coordinates (Section
\ref{sec:GP}, (3) a cross-validation procedure for targeting specific
statistical-validation criteria (Section \ref{sec:GPROMES}),  and
(4) a way to statistically quantify the error in any quantity of interest
\textit{a posteriori} by propagating the state-error model through the
associated functional.

Numerical experiments demonstrated the ability of the method to accurately
model both the in-plane and out-of-plane errors (Figures
\ref{fig:inplane_plot}, \ref{fig:outofplane_plot},
\ref{fig:inplane_plot_test2}, \ref{fig:outofplane_plot_test2}),
quantity-of-interest errors (Figures \ref{fig:Lu_output_scatter_plot},
\ref{fig:Squaredu_output_scatter_plot}, \ref{fig:output_error_plot},
\ref{fig:MAXVM_output_scatter_plot}, \ref{fig:output_error_plot_chess_matrix},
\ref{fig:VMSn1}), and realize a more computationally efficient methodolgy than
a `ROM-only approach in the case of nonlinear stationary systems (Figure
\ref{fig:accuracy_vs_efficiency_test2}).

In both numerical experiments, it was not possible to rigorously validate the
Gaussian assumption underlying the proposed statistical model (Figures
\ref{fig:histGPMSE} and \ref{fig:histGPMSEtest2}). As such, we proposed the
use of specific loss functions
(e.g., matching $\omega$-prediction intervals, minimizing the
Komolgorov--Smirnov statistic) for hyperparameter selection that enabled the
statistical model to satisfy a subset of targeted statistical-validation criteria
(Tables \ref{Fig:ValidationFrequencyInplaneTwo} and \ref{Fig:ValidationFrequencyInplaneTwoTestTwo}).

Future work includes developing stochastic-process models associated with
different distributions; this will enable a wider range of
statistical-validation criteria to be met by the constructed model. In
addition, we aim to extend the proposed methodology to dynamical systems.

\begin{appendix}

\end{appendix}

\bibliographystyle{siam}
\bibliography{biblio}

\end{document}